\documentstyle[righttag,12pt]{amsart}

\newtheorem{theorem}{THEOREM}[section]
\newtheorem{definition}[theorem]{Definition}
\newtheorem{corollary}[theorem]{Corollary}
\newtheorem{lemma}[theorem]{Lemma}

\newtheorem{remark}[theorem]{REMARK}
\newtheorem{proposition}[theorem]{Proposition}

\def\HollowBox #1#2{{\dimen0=#1 \advance\dimen0 by -#2
       \dimen1=#1 \advance\dimen1 by #2
        \vrule height #1 depth #2 width #2
        \vrule height 0pt depth #2 width #1
        \llap{\vrule height #1 depth -\dimen0 width \dimen1} 
       \hskip -#2
       \vrule height #1 depth #2 width #2}}
 
 \def\BOX{\HollowBox{.115in}{.012in}}

\def\RR{{\Bbb R}}
\def\CC{{\Bbb C}}

\def\NN{{\Bbb N}}

\def\11{{\rm\k{.45}\vb0\k{-.142}1}}
 
\def\dbar{\overline{\partial}}

\def\ips{\langle \, \cdot \, , \, \cdot \, \rangle_s}
\def\rnp{\RR^{N+1}_+}
\def\l{\langle}
\def\r{\rangle}
\def\div{\,\text{div}\,}
\def\ra{\rightarrow}
\def\btu{\bigtriangleup}
\def\dom{\hbox{\rm dom}\,}
\def\ker{\text{ker}\,}

\def\ep{\epsilon}
\def\class{\text{class}\,}
\def\D{{\cal D}}
\def\L{{\cal L}}
\def\ss{\subseteq}
\def\supp{\text{supp}\,}
\def\til{\widetilde{\text{\,}^{\text{\,}}}}
\def\order{\text{order}\,}

\def\G{{\cal G}}
\def\p{\partial}
\def\tpx{|2\pi\xi'|}

\def\vn{\vec{n}}
\def\bw{\bigwedge\nolimits}
\def\on{\text{on}  \ }
\def\bvp{boundary value problem}
\def\pd#1#2{\frac{\partial #1}{\partial #2}}
\def\ut{\tilde{u}}

\def\Sx{\sqrt{1+\tpx^2}}
\def\e#1#2{\varepsilon^{#1}_{#2}}
\def\fe{f_{\text{e}}}
\def\ue{u_{\text{e}}}

\def\K{{\cal K}_{\Omega}}
\def\bgpf{\begin{pf*}{\bf Proof}}
\def\endpf{\end{pf*}}	

\def\grad{\text{grad}}
\font\tenrm=cmr12 scaled\magstep3
\def\bcheck{\!\!{\hbox{\tenrm\char20}}}

\begin{document}

\title[Hodge Theory in the Sobolev
Topology]{
Hodge Theory in the Sobolev
Topology for the de Rham Complex
} 

\author[L. Fontana]{Luigi Fontana}
\address{Dipartimento di Matematica\\ Via Saldini 50\\
Universit\`a di Milano\\ 
20133 Milano (Italy)}
\email{fontana@@vmimat.mat.unimi.it}
\author[S. G. Krantz]{Steven G. Krantz}
\address{Department of Mathematics\\ Washington University\\ St.
Louis, MO 63130 
(U.S.A.)}
\email{sk@@artsci.wustl.edu}
\thanks{Second author supported by 
NSF Grant DMS-9022140 during his stay at MSRI.}
\author[M. M. Peloso]{Marco M. Peloso}
\address{Dipartimento di Matematica\\ Politecnico di Torino\\ 10129 
Torino (Italy)}
\email{peloso@@polito.it}
\thanks{Third author supported in part
by the Consiglio Nazionale delle Ricerche}


\begin{abstract}The authors study the Hodge theory of the
exterior differential operator $d$ acting on $q$-forms
on a smoothly bounded domain in $\RR^{N+1}$, and on the
half space $\rnp$.  The novelty
is that the topology used is not an $L^2$ topology but
a Sobolev topology.  This strikingly alters the problem
as compared to the classical setup.  It gives rise to
a \bvp\  belonging to 
a class of problems
first introduced by Vi\v{s}ik and Eskin, and by
Boutet de Monvel.
\end{abstract}

\maketitle
  
\tableofcontents
\setcounter{section}{-1}
 
\section*{\bf {PRELIMINARIES}}

\section{Introductory Remarks}
 
Fix a smoothly bounded domain $\Omega \ss
\RR^{N+1}$.  In classical treatments of the $d$ operator
(see \cite{SWE}), one considers the complex
$$
\bigwedge\nolimits^q \stackrel{d}{\longrightarrow}
  \bigwedge\nolimits^{q+1}
\stackrel{d}{\longrightarrow} \cdots .
$$
Here $\bigwedge\nolimits^q$ denotes the $q$-forms of
Cartan and de Rham having $L^2(\Omega)$ coefficients
and the operator $d$ is understood to be densely
defined.  One considers the operator $\BOX = d d^* +
d^* d.$ Here the adjoints are calculated in the
$L^2(\Omega)$ topology.
 
Now $\BOX$ makes sense on those forms $\psi$ such that
$\psi \in \text{dom}\, d^*$ and $d\psi \in
\text{dom}\, d^*$.  One can decompose the space
$\bigwedge\nolimits^q$ into (the closure of) the image
of $\BOX$ and its orthogonal complement.  Then one
exploits this decomposition to construct a right
inverse for $\BOX$.  This inverse is easily used to show that
$\BOX$ and its accompanying boundary conditions form a
second order elliptic boundary value problem of the
classical (coercive) type.
 
In 1963, J.J. Kohn \cite{KOH1} determined how to carry out
the analog of these last calculations for the
$\overline{\partial}$ operator of complex analysis on
a strongly pseudoconvex domain $\Omega$ in $\CC\,^n.$
This is the so-called $\overline{\partial}$-Neumann problem.
Of course
this analysis, while similar in spirit, is much more complicated. 
It gave rise to the important ``Kohn canonical
solution'' to the equation $\overline{\partial} u =
f.$ That is the solution $u$ that is orthogonal to
holomorphic functions {\em in the $L^2(\Omega)$ topology}.
 
Experience in the function theory of several complex
variables has shown that it is useful to have many
different canonical solutions to the
$\overline{\partial}$-problem.  For instance,
in the strongly pseudoconvex case we
profitably study the Kohn solution by comparing it
with the Henkin solution (not canonical, but nearly so),
see \cite{HEN}, and the Phong solution (determined in the
$L^2(\partial \Omega)$ topology rather than the
$L^2(\Omega)$ topology), see \cite{PHO}.
 
The ultimate goal of the program that we are initiating
in this paper is to construct solutions to the
$\overline{\partial}$ problem that are orthogonal to
holomorphic functions in a Sobolev space $W^s$ inner
product.  There are {\em a priori} reasons for knowing
that this program is feasible. First, Sobolev space is
a Hilbert space, so there must be a minimal solution
in the Sobolev topology.  Second, Boas \cite{BOA} has
studied the space of $W^s$ holomorphic functions as a
Hilbert space with reproducing kernel.  The associated
Bergman projection operator is of course closely
related to the Neumann operator for the
$\overline{\partial}$ problem.
 
The present paper carries out the first step of the
proposed program.  We work out the Hodge theory for
the exterior differentiation operator $d$ in the inner
product induced by the Sobolev space $W^s$ topology.
Of particular interest are the boundary conditions that
arise when we calculate the adjoint $d^*$ in this
topology, and the elliptic boundary value problem that
arises when we consider $\BOX = d d^* + d^* d.$ We
calculate a complete existence and regularity theory.
 
In this paper, we restrict attention to the case $s = 1$.  This
is done both for convenience and to keep the notation relatively
simple---even in this basic case the calculations are often unduly
cumbersome.  In geometric applications, the case $s=1$ is already of
great interest.  We leave the detailed treatment of higher order $s$
to a future paper.
 
In future work, we will carry out this program of analysis in the
Sobolev space topology for the $\overline{\partial}$-Neumann problem
on a strongly pseudoconvex domain.  Not only will this give rise to
new canonical solutions for the $\overline{\partial}$ problem, but it
should give a new way to view the Sobolev regularity of the 
classical $\dbar$ problem, and of understanding the subelliptic gain 
of $1/2$ 
in regularity.

Essentially the paper is divided into two parts. In the first of these 
we study the \bvp\ that arises from our Hodge theory
on the special domain given by the half space,
and in the second one we 
deal with the problem on a smoothly bounded domain.  
{\it We would like point out that 
the second part has been written so that 
it can be read independently from the
first part.\/}  
When we use results from the half space case, we give
precise reference to them. 

In detail,
the plan of the paper is as follows:  {\bf Section 1} introduces
basic notation and definitions, while in {\bf Section 2} 
we formulate the problem and states the main results.
{\bf Sections 3} through
{\bf 6} are devoted the problem on the half space.

In {\bf Section 3} we 
calculate the operator $d^*$ on 1-forms and 
also calculate its domain when the region under
study is the upper half space $\rnp$.
{\bf Section 4} completes the detailed calculation of $d^*$.
{\bf Section 5} applies a pseudodifferential
formalism developed by Boutet de Monvel
to the study of our elliptic boundary value problem. This section 
is included essentially 
to show how the \bvp\ under investigation can be view as an example
of a very general kind of problem first introduced in
[BDM1]-[BDM3] and in \cite{ESK}.
In {\bf Section 5} we explicitly
construct the solution of the \bvp\ on 
the half space in the case of functions.
{\bf Section 6} 
studies the problem for $q$-forms.
{\bf Sections 7} through {\bf 12} deals with the problem on a bounded
smooth domain. In {\bf Section 7} we set up the problem and calculate
the domain of $d^*$ and the semi-explicit expression for $d^*$.
In {\bf Section 8} we introduce notation and some technical facts needed
in the sequel.
{\bf Section 9} contains the proof of the coercive estimate, and the
proof of our result about existence of solutions for the \bvp.
{\bf Section 10} is devoted to the proof of the {\em a priori}
estimate in the case of functions.
{\bf Section 11} 
gives the proof of the regularity result in the case
of $q$-forms.  Finally, in {\bf Section 12}, we conclude the proof of
our main result in the case of a smoothly bounded domain.

We thank G. Grubb for helpful communications regarding
this work.
 
\section{Basic Notation and Definitions}
 
We use the symbol $d$ to denote the usual operator of exterior
differentiation
acting on $q$-forms.  We let $\Omega$ denote a smoothly bounded 
domain
in $\RR^{N+1}.$  Usually, for simplicity only,
a domain is assumed to be connected.
 The symbol $\bigwedge\nolimits^q(\Omega)$ denotes the
$q$-forms on $\Omega$ with smooth coefficients.  
The symbol $\bigwedge\nolimits_0^q(\Omega)$ denotes the
forms with coefficients that are $C^\infty$ and
compactly supported in $\Omega.$ 
We let $\bigwedge\nolimits^q(\overline{\Omega})$ denote
the $q$-forms with coefficients that are smooth on
$\overline{\Omega}$, and
$\bigwedge\nolimits^q_0(\overline{\Omega})$ denote the
$q$-forms with coefficients in
$C^\infty(\overline{\Omega})$ and having compact
support in $\overline{\Omega}$ (that is, the support may
not be disjoint from the boundary of $\Omega$).

Some
of our explicit calculations will be performed on the special
domain $\RR^{N+1}_+ = \{x =
(x_0,x_1,\dots,x_N)\in \RR^{N+1}: x_0 > 0\}$.  
The half space is of course unbounded, so the function 
spaces we deal
with must take into account the integrability at infinity.  On the 
other hand the half space has the advantage of allowing explicit
calculations. 

If ${\cal L}$ is an operator on forms, then ${\cal L}'$
denotes its {\em formal adjoint}, that is, its adjoint
calculated when acting on elements of
$\bigwedge\nolimits_0^p(\Omega)$.  
 
In the discussion that follows we let $\Omega$ denote either
a smoothly bounded domain $\Omega$, or the half space
$\rnp$.  If $s$ is a
non-negative integer then the $s$ order Sobolev space
norm on functions on $\Omega$ is given by
$$
\|f\|^2_s =
\sum_{|\alpha| \leq s} \left \|
   \frac{\partial^\alpha f}{\partial x^\alpha}
    \right \|^2_{L^2(\Omega)} .
$$
The associated inner product is
$$
\left \l f, g \right \rangle_s =
\sum_{|\alpha| \leq s} \int_\Omega
     \frac{\partial^\alpha f}{\partial x^\alpha} \
   \overline{\frac{\partial^\alpha g}{\partial x^\alpha}}
               \, dV(x) .
$$
Here $dV$ stands for ordinary Lebesgue volume measure.

For $s$ a non-negative integer we define the Sobolev space $W^s
(\Omega)$ as the closure of $C^\infty _0 (\overline{\Omega})$.
When $s\in\RR_+$ we define $W^s (\Omega)$ by interpolation 
(see \cite{LIM}
for instance).  Moreover, for $s\in\RR_+$, we denote by
$\stackrel{\circ}{W^s}(\Omega)$ 
the closure of $C^\infty _0 (\Omega)$ in $W^s (\Omega)$.  
When $s<0$ we define the negative Sobolev space 
$W^s (\Omega)$ to be the
dual of $\stackrel{\circ}{W^s}(\Omega)$ 
with respect to the standard $L^2$-pairing.

On the Euclidean space $\RR^{N+1}$ we consider the Fourier transform
defined initially for a testing function $f\in C^\infty_0$ as
$$
\widehat{f}(\xi) \equiv {\cal F} f(\xi)
    \equiv \int_{\RR^{N+1}}  f(x) e^{-2\pi i x\cdot\xi}\, dx .
$$
We will also consider the tangential Fourier transform of functions
defined on the half space $\rnp$: If $f\in C^\infty_0
(\overline{\rnp})$ we set 
$$
\hat f (x_0,\xi') = \int_{\RR^N} f(x_0,x') e^{-2\pi ix'\cdot\xi'} \, dx' .
$$
We denote the inverse tangential Fourier transform of a function
$g(x_0 ,\xi')$ by $\check g (x_0 ,x')$.

For any $s\in\RR$ the Sobolev space $W^s (\RR^{N+1})$ can 
be defined via the Fourier transform.  Indeed, we set
$$
W^s ({\Bbb R}^{N+1}) = \bigl\{ f\in L^2 ({\Bbb R}^{N+1}) :
\int_{{\Bbb R}^{N+1}} (1+|\xi|^2)^s |{\cal F}f(\xi)|^2 \, d\xi
<\infty \bigr\} .
$$

We will also consider the Sobolev spaces $W^s (b\Omega)$ defined on
the boundary of our domain, $s\in \RR$. In the case $\Omega=\rnp$, 
$W^s (b\Omega)$ is just the classical Sobolev space on $\RR^N$.  In
the case of a smoothly bounded domain $\Omega$, 
the Sobolev space can be defined by fixing a smooth atlas 
$\{ \chi_j\}$ on $\partial \Omega$, 
letting $\phi_j$ be a partition of unity
subordinate to this atlas, and defining the 
Sobolev norm of a function $f$  on $b\Omega$ as 
$$
\| f \|^2_{W^s (b\Omega)} 
=   \sum_j \| \phi_j f\circ {\chi_j}^{-1} \|_{W^s (\RR^N)}^2 .
$$
Of course this norm is highly non-intrinsic, 
but a different choice of
an atlas gives rise to an equivalent norm.

On the half space $\rnp$, and on a bounded domain $\Omega$,
we consider the space of $q$-forms with
coefficients in $W^s$.  We denote such spaces by $W^s_q (\rnp)$ and
$W^s_q (\Omega)$ respectively.
Fix the standard basis for $q$-forms:
$\{ dx^I \}$,  where $I=(i_1,\dots,i_q)$ is an increasing
multi-index, i.e.\ $0\le i_1<\cdots<i_q \le N$. The Sobolev inner
product on the space of $q$-forms in the case $s$ a non-negative 
integer is given by
$$
\l \phi,\psi\r_s 
 = \sum_I \l \phi_I ,\psi_I \r_s 
 = \sum_I \sum_{|\alpha|\le s} \int
\pd{^\alpha \phi_I}{x^\alpha} 
\overline{\pd{^\alpha \psi_I}{x^\alpha}} dV(x) ,
$$
where the integral is taken over all of space, i.e.\ 
either on $\rnp$ or on $\Omega$.
In the case of general $s$, for $\phi\in\bw^q$, 
$\phi=\sum_I \phi_I dx^I$,  we
define the norm of $\phi$ by
$$
\| \phi\|^2_{W^s} = \sum_I \| \phi_I \|^2_{W^s} .
$$ 

Of course $W^s$
is a Hilbert space when equipped with the foregoing
inner product.  (In particular, $W^s$ can be identified
with its own dual in a natural way.)     
 
Until further notice, we use $D_j$ to denote the
partial derivative in the $j^{\rm th}$ variable, and
if $\alpha = (\alpha_0,\alpha_1,\dots,\alpha_N)$ is a
multi-index then we let $D^\alpha$ denote the
corresponding differential monomial.

\section{Formulation of the Problem and Statement of the Main
Results} 
\setcounter{equation}{0}
Given the operator $d$,
$$
d :  \bw^q \longrightarrow \bw^{q+1},
$$
we think of it as a densely defined operator on the space 
$W^1_q $ of 
$q$-forms with coefficients in $W^1$, both in the case of the half
space $\rnp$ and also in the case of a smoothly bounded domain $\Omega$. 
Let $d^*$ denote the $W^1$-Hilbert space adjoint of $d$.  It is a
densely defined (unbounded) operator 
$$
d^* : W^1_{q+1} \longrightarrow W^1_q .
$$ 
We shall study
the \bvp\ 
\begin{equation*}
\begin{cases}
\displaystyle{ (dd^* +d^* d) \phi= \alpha }& \on \rnp \ 
			(\text{resp.\ }\on \Omega) \\
\phi\in\dom d^* & \\
d\phi\in\dom d^* & 
\end{cases} \ ,
\end{equation*}
for $\alpha\in W^s_q (\Omega)$;
we shall prove existence and
regularity theorems, both in the case of the half space  and of a
smoothly bounded domain, for $q$-forms, $q=0,1,\dots,N+1$. 
The conditions $\phi,\, d\phi\in\dom d^*$ are ultimately 
expressed as boundary conditions. 

Our main results are the following.
\begin{proposition} \sl Let $\Omega$ denote a smoothly bounded domain, or
the special domain $\rnp$. 
Let $q=0,1,\dots,N$.  Then we have
$$
\dom d^* \cap \bw^{q+1}_0 (\overline{\Omega})
= \bigl\{ \phi\in \bw^{q+1}_0 (\overline{\Omega}) : \nabla_{\vn} 
\phi\lfloor \vn \bigl|_{b\Omega} =0 \bigr\} .
$$
Here $\nabla_{\vn}$ denotes the covariant differentiation 
of a form in
the normal direction, and ``$\lfloor$" the contraction operator
between a form and a vector field.
\end{proposition}
The next result shows a striking difference with the classical case.
We begin with the case of $\rnp$.  
Here, and in the rest of the
paper, we denote by $\btu'$ the Laplace operator defined on the
boundary of the half space:
$$
\btu'=\sum_{j=1}^{N} \pd{^2}{x_j^2} .
$$
\begin{proposition}  \sl
Let $d'$ denote the formal adjoint of $d$.  Then on $\dom d^*$ we have
$$
d^* = d' + {\cal K},
$$
where ${\cal K}$ is an operator sending $(q+1)$-forms to $q$-forms.
The operator ${\cal K}$ is the solution operator of the following \bvp
$$
\begin{cases}
\displaystyle{(-\btu+I)({\cal K}\phi)=0 }& \on \rnp\\
\displaystyle{\pd{}{x_0}({\cal K}\phi)=(-\btu'+I)
\left (\phi\lfloor\bigl( \pd{}{x_0} \bigr) \right )}& \on \RR^N 
= b\RR^{N+1} .
\end{cases} 
$$
Explicitly, if $\phi\in\bw^{q+1}(\overline{\rnp})$,
$\phi=\sum_{|I|=q+1}\phi_I \, dx^I$, then
$$
({\cal K}\phi)_I \hat{\text{\, }}(x_0,\xi')
= -\Sx e^{\Sx x_0} \hat\phi_{0I} (0,\xi') .
$$
\end{proposition}
In the case of a smoothly bounded domain we have the following:
\begin{proposition} \sl
Let $\Omega$ be a smoothly bounded domain.
Then, on $\dom d^*$,
$$
d^* = d' + \K,
$$
where $d'$ is the formal adjoint of $d$, and 
$\K$ is an operator sending $(q+1)$-forms to $q$-forms.
The components of ${\cal K}\phi$ 
are solutions  of the following \bvp s
$$
\begin{cases}
\displaystyle{(-\btu+I)(\K\phi)_I =0 }& \on \Omega \\
\displaystyle{\pd{}{n} (\K\phi)_I
		=T_2\phi\lfloor \vn }& \on b\Omega
\end{cases} ,
$$
where $T_2$ is a second order tangential differential operator on
forms (to be defined in detail below). 
\end{proposition}
At this point we fix the following notation that we will use
throughout the entire paper.
\begin{definition}\label{def-G}   \rm
On the half space $\rnp$ we set
$$
G= {\cal K}d+d{\cal K},
$$
and for a smoothly bounded domain $\Omega$ we set
$$
G_\Omega =\K d+d\K  .
$$
\end{definition}
Notice that our \bvp s become
\begin{equation}\label{bvp-half-space}
\begin{cases}
\displaystyle{ (-\btu+G) \phi= \alpha }& \on \rnp \\
\phi\in\dom d^* & \\
d\phi\in\dom d^* & 
\end{cases} ,
\end{equation}
and
\begin{equation}\label{bvp-domain}
\begin{cases}
\displaystyle{ (-\btu+G_\Omega) \phi= \alpha }& \on \Omega \\
\phi\in\dom d^* & \\
d\phi\in\dom d^* & 
\end{cases} .
\end{equation}
We are now ready to state our results about existence and regularity
of the \bvp.
\begin{theorem}  \sl
Let $s>1/2$. 
Let $\alpha\in W^s_q (\rnp)$,
$\text{supp}\, \alpha \subseteq \{ x:|x|<R \}$.  
Let $N\ge4$,
Then there exists a unique tempered distribution
$\phi$ solution of the \bvp\ (\ref{bvp-half-space})
such that
$$
\| \phi \|_{s+2} \le c \| \alpha\|_s ,
$$
where $c=c(s,R) >0$ does not depend on $\alpha$, nor on $\phi$.

If $N=2,3$ and $\int \alpha_I \, dx=0$, and if $N=1$ and $\int
x^\beta \alpha_I \, dx =0$ when $|\beta|\le1$, then the conclusion as
above still holds.
\end{theorem}
\begin{theorem}\label{MAIN-THM-DMN} \sl
Consider the \bvp\  (\ref{bvp-domain}).
Let $s>1/2$. 
Then there exists a finite dimensional subspace 
(the harmonic space) ${\cal H}_q$ 
of $\bw^q (\overline{\Omega})$
and a constant $c=c_s >0$ such
that if 
$\alpha\in
W^s_q (\Omega)$
is orthogonal (in the
$W^1$-sense) to ${\cal H}_q$, then
the \bvp\ (\ref{bvp-domain})
has a unique solution $\phi$ orthogonal to ${\cal H}_q$ 
such that
$$
\| \phi \|_{s+2} \le c \|  \alpha\|_s .
$$
\end{theorem}
{\bf Remark.} If $s<1$, by saying that $\alpha$ is orthogonal in the
$W^1$-sense we mean that $\alpha$ is $W^s$-limit of smooth forms that
are orthogonal in $W^1$ to ${\cal H}_q$.

Finally we have
\begin{theorem}  \sl
Let $\Omega$ be a smoothly bounded domain in $\RR^{N+1}$.  
Let $W^1_q (\Omega)$ 
denote the 1-Sobolev space of $q$-forms.  Then we have the
strong orthogonal decomposition
$$
W^1_q = dd^* (W^1_q) \bigoplus d^* d(W^1_q) \bigoplus {\cal H}_q \, , 
$$
where ${\cal H}_q$ is a finite dimensional subspace.
\end{theorem}

\newpage
\section*{\bf {THE PROBLEM ON THE HALF SPACE}}
\mbox{}

\section{The Operator $d^*$ on $1$-Forms and Its Domain}
\setcounter{equation}{0}
 
The aim of this section is to determine $\dom d^* \cap \bw^1
(\overline{\rnp})$, and to give an explicit expression for $d^*$.  

Easy computations show that the formal adjoint $d'$ of $d$ with
respect to the inner product of $W^1$ is the same as
the formal adjoint in the classical $L^2$-case. 
Therefore, for $\phi \in
\bigwedge_0^1(\RR^{N+1}_+)$, $\phi = \sum_{j=0}^N
\phi_j \, dx^j,$ we see that 
$$ 
d' \phi = -
\sum_{j=0}^N D_j \phi_j \equiv - \text{div}\, \phi .
$$ 
It is not hard to see that the specific form of $d'$ is independent of
which order $s$ of Sobolev inner product we use.

Next we want to compute the Hilbert space adjoint $d^*$
of $d$.  Recall that $W^s_q$ will be the
closure of the smooth $q$-forms
$\bigwedge\nolimits^q_0(\overline{\RR^{N+1}_+})$ 
with compact support in
$\overline{\RR^{N+1}_+}$ 
with respect to the norm
$$ 
\| \phi  \|_s \equiv
\sum_{|\alpha| \leq s \atop |I| = q}
\int_{\RR^{N+1}_+} |D^\alpha \phi _I|^2 \, dV < \infty .
$$
Let $u \in
\bigwedge\nolimits^0$, $\phi \in \sum_{i=0}^N \phi_i \,
dx^i$ have compact support in $\overline{\RR^{N+1}_+}$
(that is, the support is compact but not necessarily
disjoint from the boundary).  Then
\begin{eqnarray*}
  \l du, \phi \rangle _s 
& = & \sum_{j=0}^N \l D_j u, \phi_j \r_1  \\
 & = & \sum_{k=0}^{N} \left ( \sum_{j=1}^N
 \int_{\RR^{N+1}_+} D_j D_k u
\overline{D_k \phi_j} \, dV
 + \int_{\RR^{N+1}_+} D_0 (D_k u)
\overline{D_k \phi_0} \, dV \right ) .
\end{eqnarray*}
The first sum inside the parentheses in the
last line equals
$$
- \sum_{j=1}^N \int_{\RR^{N+1}_+} D_k u \,
\overline{D_j D_k \phi_j} \, dV
$$
since $u$ and $\phi_j$ have compact support and
the derivation is in the
tangential directions (i.e.\ the directions
$x_1,\dots,x_N$).  The second
term in parentheses equals (with $x' = (x_1,\dots,x_N)$)
\begin{multline*}
 \int_{\RR^N} \left ( \left. D_k u
\overline{D_k \phi_0} \right |_{x_0 = 0}^{x_0 = \infty}
 - \int_0^\infty D_k u \overline{D_0 D_k
\phi_0} \, dx_0 \right ) \, dx' \\
  =  - \int_{\RR^N} D_k u \overline{D_k \phi_0}
\biggr |_{x_0 = 0} \, dx'
 - \int_{\RR^{N+1}_+} D_k u
\overline{D_0 D_k\phi_0} \, dV(x) .
\end{multline*}
 
Therefore, combining the last several lines, we find that
\begin{eqnarray*}
\l du, \phi \r_1  & = & \sum_{k=0}^{N}
\left ( - \sum_{j=0}^N
 \int_{\RR^{N+1}_+} D_k u
\overline{D_j D_k \phi_j}
 - \int_{\RR^N} D_k u \overline{D_k \phi_0}
\biggr |_{x_0 = 0} \, dx' \right ) \\
  & = &  \l u, d'\phi \r_1  
-\int_{\RR^N} u \overline{\phi_0} \biggr|_{x_0 =0}\, dx' 
-\sum_{k=1}^{N} \int_{\RR^N} D_k u
\overline{D_k \phi_0} \biggr |_{x_0 = 0} \, dx' .
\end{eqnarray*}
 
Now we can determine the domain of $d^*$ (we 
must do this before we can calculate the operator
$d^*$ itself).
Recall that, by the definition of Hilbert space adjoint,
$\phi \in \text{dom}\, d^*$
if and only if there is a number $c_\phi$ such that for
all $u \in \bigwedge\nolimits^0$
we have
$$
\bigl | \l du, \phi \r_1  \bigr
| \leq c_\phi \|u\|_1 .
$$
\begin{proposition}   \sl
The $1$-form $\phi$ lies in $\bigl(\dom d^* \bigr) 
\cap \bw^1(\overline{\rnp})$ 
if and only if
$$
D_0 \phi_0 \biggr |_{x_0 = 0} \equiv 0 .
$$
\end{proposition}
\bgpf  
We begin by setting
$$
{\cal D} = \dom d^* \cap
\bw_0^1 (\overline{\rnp})
$$
and
$$
{\cal E} = \left\{ \phi \in \bw^1_0
(\overline{\rnp}):
  \pd{}{x_0} \phi_0(0,x')
\equiv 0 \right\} .
$$
We first show that ${\cal E} \ss {\cal D}.$
Let $\phi \in {\cal E}.$
We have already computed that 
\begin{align*}
\l du,\phi \r_1 
&  =  \l u, d'\phi \r_1 
-\int_{\RR^N} u \overline{\phi_0} \biggr|_{x_0 =0}\, dx' 
-\sum_{k=1}^{N} \int_{\RR^N} D_k u
\overline{D_k \phi_0} \biggr|_{x_0 = 0} \, dx' \\
&  =  \l u, d'\phi \r_1  + \int_{\RR^N} u \biggr|_{x_0 =0}
(\btu' \phi_0 - \phi_0)\biggr|_{x_0 =0} \, dx' .
\end{align*}
Clearly the mapping
$$
u \longmapsto \l u, d'\phi \r_1
$$
is bounded on $W^1(\rnp)$, and so
$$
| \l du,\phi \r_1 | \le c_\phi \bigl(
\| u\|_{W^1 (\rnp)} + \| u |_{x_0 =0} \|_{L^2 (\RR^N)} \bigr) . 
$$
(Note that the bound
here depends on $\phi$,
but that is acceptable.)
Now, recall the standard trace theorems for Sobolev
spaces (see either \cite{LIM} or \cite{KR1}). Then
$$
\biggl \| u |_{x_0 =0} \biggr \|_{L^2 (\RR^N)} \le \| u \|_{W^t (\rnp)} ,
$$
for $t>1/2$, and in particular for $t=1$.
This establishes the containment
${\cal E} \ss {\cal D}.$
 
Next we establish that ${\cal D} \subset
{\cal E}.$  Seeking a contradiction, we suppose that
$\phi \in {\cal D} \equiv \text{dom}\, d^*
\cap \bw^1_0(\rnp)$
and that
$D_0 \phi_0 (0,x')$
is not identically zero.  We may suppose that
$\phi_0$ is real and that
$$
D_0 \phi_0(0,x') \geq 1 \qquad
\qquad \text{when} \ \ |x'| \leq 2
$$
(just multiply by a suitable constant and scale).
 
Now set
$$
u_\epsilon(x_0,x') =  (x_0 + \epsilon)^{3/4} \chi(x_0,x') , 
$$
where $\chi$ is a cutoff function with support in
$\{|x_0| \leq 2\} \times \{|x'| \leq 2\}$
and such that
$\chi \equiv 1$ in the set $\{|x_0| \leq 1\}
\times \{|x'| \leq 1\}.$
 
We claim that there is a constant $C > 0$,
independent of $\epsilon,$
such that
$$
\|u_\epsilon\|_{W^s(\rnp)} \leq C \qquad \qquad
\forall\epsilon,\ 0 < \epsilon \leq \epsilon_0 \, ,
                    \leqno {(i)}
$$
and
$$
\int_{\RR^N} D_0 u_\epsilon
  \overline{D_0 \phi_0} \biggr |_{x_0 = 0}
   \, dx' \, \sim \, \frac{1}{\epsilon^{1/4}}
 \qquad \quad \text{as} \ \ \epsilon \rightarrow 0 .
\leqno {(ii)}
$$
 
 Once these two facts are established, it is clear that it is not
possible to find a constant such that the mapping 
$u\longmapsto\l du, \phi\r_1$ is bounded on $W^1(\rnp)$, and thus
$\phi$ cannot be in $\dom d^*$.

 Fact $(ii)$  follows since 
$$
D_0 u_\epsilon (x_0,x') \biggr|_{x_0=0} = 
\left( \frac{3}{4}
 \epsilon^{-1/4} \chi(0,x') + \epsilon^{3/4}
D_0 \chi(0,x') \right).
$$
In order to prove $(i)$ we need only 
check the size of
$$
\int_{\rnp} |D_0 u_\epsilon|^2 dV .
$$
But
$$
D_0 u_\epsilon(x_0,x') = \frac{3}{4} \
\frac{1}{(x_0 + \epsilon)^{1/4}}
 \, \chi(x_0,x') + \Phi_\epsilon(x_0,x') ,
$$
where $\Phi_\epsilon$ lies in $C_0^\infty(\rnp)$ and is
uniformly bounded in $\epsilon.$
Thus
\begin{eqnarray*}
\int_{\rnp} |D_0 u_\epsilon|^2 \, dV &
\leq & C \int_0^1 \frac{1}{(x_0 + \epsilon)^{1/2}} \, dx_0
     + C' \\
    & \leq & C ,
\end{eqnarray*}
independent of $\epsilon$ as $\epsilon \ra 0.$
This establishes $(i)$.
\endpf 
\subsection*{The operators ${\cal K}$ and $d^*$}
Our next task is to determine the
explicit expression for the operator $d^*$.  We will learn
that, contrary to the classical case, the operator $d^*$ when
restricted to $\dom (d^*)\cap\bw^1$ {\em does not} coincide with
$d'$.  
Once we have the expression for $d^*$ then we are in a position to
formulate the \bvp\ in the case $q=1$.
The fundamental
calculation is this:
 
\begin{proposition}\label{expr-d*}   \sl
If $\phi \in \hbox{\rm dom}\,
d^* \cap \bw^1_0 (\overline{\rnp})$,
$\phi = \sum_{j=0}^N \phi_j \, dx^j$, then
$$
d^* \phi = - \hbox{\rm div}\, \phi + \int_{\RR^N}
\frac{\partial}{\partial x_0}
 \left ( e^{-\sqrt{1 + 4\pi^2 |\xi'|^2} x_0} \right )
 \widehat{\phi_0}(0,\xi') e^{2\pi i x' \cdot \xi'} \, d\xi'  
$$
\end{proposition}
\bgpf  
We need to determine the
operator $d^*$ that satisfies
$$
\l df, \phi \r_1
= \l f, d^* \phi \r_1 .
$$
The spirit of the calculation that follows
is to rewrite the expressions
that involve the inner product
$\ips$ in terms of expressions that involve
only the $L^2$ inner product.
 
Now
\begin{eqnarray*}
 \l df, \phi \r_1 & = & \int_{\rnp}
\sum_{j=0}^N \pd{f}{x_j}
       \overline{\phi_j} \, dV +
  \int_{\rnp} \sum_{j,k=0}^N
\pd{^2 f}{x_j \partial x_k}
	\overline{\pd{\phi_j}{x_k}}\, dV \\
& \stackrel{\text{(parts)}}{=}& 
	\l f,-\div \phi\r_1 - \int_{\RR^N} f(0,x')
\overline{\phi_0(0,x')} \, dx' \\
&  & \qquad + \int_{\RR^N} f(0,x') \overline{\sum_{k=1}^N
\pd{^2}{x_k^2} \phi_0(0,x')} \, dx' .
\end{eqnarray*}
Of course we have used the boundary condition that occurs in our
characterization 
of $\dom d^*$. 
 
Next we write
$$
d^* \phi = -\div \phi +u,
$$
and we wish to determine the explicit
expression for $u\equiv {\cal K}\phi$,
i.e.\ the formula for the operator ${\cal K}$
mapping  $1$-forms into functions, such that
$d^* =d'+{\cal K}$.  Then we have
\begin{eqnarray*}
\l df,\phi\r_1
& =& \l f,-\div \phi\r_1 +\l f,u\r_1 \\
& = & \l f,-\div \phi\r_1 +\int_{\rnp} f\overline{u}\, dV\\
& & \qquad -\int_{\rnp} f\overline{\btu u} \, dV
-\int_{\RR^N}f(0,x')\overline{\frac{\p u}{\p x_0}(0,x')}\, dx'.
\end{eqnarray*}
Here we have integrated by parts.
 
Comparing the two last calculations we
see that the function $u$ must
satisfy the following equation:
\begin{multline*}
- \int_{\RR^N} f(0,x') \overline{\phi_0(0,x')}\, dx'
  +   \int_{\RR^N} f(0,x') \overline{\sum_{k=1}^N
\frac{\partial^2}{\partial x_k^2} \phi_0(0,x')} \, dx' \\
 =  \int_{\rnp} f\overline{u}\, dV
         -\int_{\rnp} f\overline{\btu u} \, dV
-\int_{\RR^N}f(0,x')\overline{\frac{\p u}{\p x_0}(0,x')}\, dx'.
\end{multline*}
This equality must hold for all 
$f\in C_0^\infty(\overline{\rnp})$.
By choosing functions  supported away from the boundary we
obtain a differential equation on $\rnp$; inserting this in the
equation above gives rise to a boundary condition.  Thus the
function $u\in W^1 (\rnp)$
must be a solution of the following \bvp:
$$
\begin{cases}
\displaystyle{
-\btu u + u  =  0} & \text{on} \ \rnp \\
\displaystyle{
\pd{u}{x_0}(0,\cdot)   =  -\btu' \phi_0 (0,\cdot)
        +\phi_0(0,\cdot)}  & \on\{ 0\}\times\RR^N 
\end{cases} .
$$
Using the partial Fourier transform with
respect to $x'$, we find that
the equation $\btu u - u = 0$ becomes
$$
\frac{\partial^2}{\partial x_0^2} \widehat{u} (x_0,\xi') -
 (1 + |2\pi\xi'|^2) \widehat{u}(x_0,\xi') = 0 ,
 \qquad \text{with} \ \ \widehat{u} \in L^2(\rnp) .
$$

For fixed $\xi'$, this is an ordinary differential
equation in $x_0$
whose solution is
$$
\widehat{u}(x_0,\xi') = a(\xi')
e^{\sqrt{1 + |2\pi\xi'|^2}x_0} +
 b(\xi') e^{-\sqrt{1 + |2\pi\xi'|^2}x_0} .
$$
Since $\widehat{u} \in L^2(\rnp)$, it must be
that $a(\xi') \equiv 0$.
Thus
\begin{equation}\label{u-hat}
\hat{u}(x_0,\xi') =
b(\xi') e^{-\sqrt{1 + |2\pi\xi'|^2}x_0}
\end{equation}
and the solution to the above \bvp\ is
$$
{\cal K} \phi(x_0,x') =
\int_{\RR^N} b(\xi')
e^{-\sqrt{1 + |2\pi\xi'|^2} x_0}
 e^{2\pi i \xi' \cdot x'} \, d\xi' .
$$
Now we determine $b(\xi')$ in such a
way that the boundary
conditions  are satisfied.  We compute the derivative
of $\hat u(x_0,\xi)$, as given by (\ref{u-hat}),
with respect to
$x_0$, and evaluate at $x_0 =0$.
We find that
$$
b(\xi') \sqrt{1 + |2\pi\xi'|^2}
=  - |2\pi\xi'|^2 \widehat{\phi_0} (0,\xi')
- \widehat{\phi_0} (0,\xi') ,
$$
that is,
$$
b(\xi') = -
\sqrt{1 + |2\pi\xi'|^2} \widehat{\phi_0} (0,\xi') .
$$
Therefore the Hilbert space adjoint of
$d$ (and this has been the main
point of these calculations) is
\begin{eqnarray*}
 d^* \phi
& \equiv  & -\div \phi +{\cal K}\phi \\
& = & - \div \phi - \int_{\RR^N} \sqrt{1 + |2\pi\xi'|^2}
      e^{-\sqrt{1 + |2\pi\xi'|^2} x_0}
  \widehat{\phi_0}(0,\xi') e^{2\pi i \xi' \cdot x'}
            \, dx' \\
  & \equiv & - \div \phi + \biggl (
K_{x_0} * \phi_0(0, \, \cdot \, ) \biggr ) (x') ;
\end{eqnarray*}
here  the convolution is in $\RR^N$ with respect to the variable 
$x'$, and
$$
K_{x_0} (x') = \biggl (
\frac{\partial}{\partial x_0}
e^{- \sqrt{1 + |2\pi\xi'|^2}x_0}
\biggr )\bcheck (x')
$$
is the convolution kernel.
\endpf

At this point, following the classical methodology of
Hodge theory, we want to study the  problem 
(\ref{bvp-half-space}):
$$
\begin{cases}
\displaystyle{ \bigl( d d^* + d^* d \bigr) u  =  f} & \\
                u,du  \in  \dom d^* & 
\end{cases} .
$$
According to our calculations this problem 
becomes
\begin{equation} \label{fnct-bvp-half-space}
\begin{cases}
\displaystyle{
- \btu u + \int_{\RR^N}
\left ( \frac{\partial}{\partial x_0}
e^{-\sqrt{1 + |2\pi\xi'|^2} x_0} \right )
 \widehat{\frac{\partial u}{\partial x_0}}(0,x')
     e^{2\pi i \xi' \cdot x'} \, dx'  =   f} & \on \rnp \\
\displaystyle{
\pd{^2}{x_0^2} u(0,x') = 0 }& \on \RR^N 
\end{cases}. 
\end{equation}

At the conclusion of these calculations we want 
to remark that the 
\bvp s that appear in the case of higher degree forms is
treated in Section \ref{SEC-Q-FMS}.

Recall, from Definition \ref{def-G}, that
\begin{eqnarray}\label{G-fnct}
Gu(x_0,x')
& = & {\cal K}\biggl(\pd{u}{x_0} \biggr)(x_0, x')\nonumber\\
& = &\int_{\RR^N}
\left ( \pd{}{x_0}
  e^{-\sqrt{1 + |2\pi\xi'|^2} x_0} \right )
 \widehat{\pd{u}{x_0}}(0,x')
              e^{2\pi i \xi' \cdot x'} \, dx.
\end{eqnarray}
Now we define $\tilde{\cal K}$ by setting
${\cal K}\equiv \tilde{\cal K}\circ\gamma$,
where $\gamma = \gamma_0$ denotes the standard 
{\em trace} operator of
restriction to the boundary.
 
Notice that the kernel $\widehat{K_{x_0}}(\xi')$ is 
$-\Sx e^{-\Sx}x_0$, and thus is 
quite similar to
the normal derivative of the Fourier transform of the Poisson
kernel.
But actually the operator $\tilde{\cal K}$ is  more
regular  than  the Poisson integral itself.  The higher
degree of regularity of 
$\tilde{\cal K}$ is seen immediately by noticing
that the symbol of the convolution kernel $K_{x_0}$, that is
$-\Sx e^{-\Sx x_0}$, is smooth in $\xi'$ for all
$\xi'$.  We record here a few simple properties of the operators 
$\tilde{\cal K}$ and $G$.  
We denote by ${\cal S}(\RR^N)$ the space of
Schwartz functions on $\RR^N$, and by 
${\cal S}(\overline{\rnp})$ the
space of restriction to $\rnp$ of 
Schwartz functions on $\RR^{N+1}$.
We have the following proposition.
\begin{proposition}\label{tilde-K}   \sl
We have that
$$
\widetilde{\cal K} : {\cal S}(\RR^N) \longrightarrow
{\cal S}{(\overline{\rnp})}
$$
continuously.  Moreover
$$
\widetilde{\cal K} : 
W^{r+1/2} (\RR^N)\longrightarrow W^{r} (\rnp),
$$
continuously for $r\ge 0$.
In addition,
$$
G  : W^{r+2} (\RR^{N+1})\longrightarrow W^{r} (\rnp),
$$
and
$$
G : {\cal S}(\overline{\rnp}) \longrightarrow {\cal
S}(\overline{\rnp}),
$$
continuously.
\end{proposition}
\bgpf 
Recall that $G=\widetilde{\cal K}\circ\gamma\circ(\p/\p x_0)$, so
it suffices to prove the corresponding 
statements for $\widetilde{\cal K}$.
 
Let $v\in W^{r+1/2}(\RR^N)$. Then
$$
\widehat{\widetilde{\cal K}v} (x_0,\xi') = -\Sx e^{-\Sx}\hat v (\xi').
$$
It suffices to consider the case $r$ an integer, 
and use interpolation
to obtain the general case.  For $r\in\NN$ we have
\begin{eqnarray*}
\|\widetilde{\cal K}v \|^2_{W^r (\rnp)}
& \le & \sum_{j=0}^{r} \int_{\rnp}
(1+|\xi'|^2)^{r-j} \left| \frac{\p^j 
	\widehat{\widetilde{\cal K}v}}{\p x_0^j}
        (x_0,\xi')\right|^2 \, d\xi' \, dx_0 \\
& \le & c\int_{\rnp} 
(1+|\xi'|^2)^{r+1}|\hat v(\xi')|^2 e^{-\Sx x_0}
        \, dx_0 \, d\xi' \\
& \le & c \int_{\RR^N} (1+|\xi'|^2)^{r+1/2}|\hat v (\xi')|^2
        \, d\xi'\\
& \le & c \| v\|^2_{W^{r+1/2}(\RR^N)} .
\end{eqnarray*}
To show that $\widetilde{\cal K}$ is continuous between 
the (two different) Schwartz spaces, we notice that
\begin{eqnarray*}
\lefteqn{
\left| x_0^\ell (2\pi\xi')^\alpha\pd{}{x_0^\ell}
\pd{}{{\xi'}^\alpha}
\widehat{\widetilde{\cal K}v} (x_0,\xi')\right| }\\
& \le & \left| x_0^{\ell+|\alpha|}
(1+\tpx^2 )^{\frac{\ell+1}{2}+|\alpha|} e^{-\Sx x_0} \hat
v(\xi')\right| \\
& \le & c\biggl \{ \sup_{x_0\ge 0}\left[ 
x_0^{\ell+|\alpha|}e^{-x_0} \right] \biggr \}
\sup_{|\xi'|}(1+\tpx)^{\frac{\ell+1}{2}+|\alpha'|}
|\hat v(\xi')| , 
\end{eqnarray*}
which is finite for $v\in{\cal S}(\RR^N)$.
\endpf
 
We conclude this section with some remarks.
\begin{remark}\label{simbol}  \rm
We have that
$$
- \phi_0 = \sigma(d', \, dx_0) \, \phi ,
$$
where $\sigma(d',\, dx_0)$ is the symbol of
$d'$ in the normal direction.  
With this notation,
$$
\l du, \phi \r_1 = \l
u, d' \phi \r_1  +
  \int_{b\rnp} \l u,
\sigma (d',\, dx_0)\phi \r_1 \, dx' ,
$$
where $b \RR^{N+1}$ is the boundary of the
half space.
\end{remark}
Indeed, using the
standard definition (see, for instance, \cite{KR1}) of
symbol, we select a testing function $\rho$ such that
$\rho(0,x') = 0$ and $d\rho = \, dx_0$ (for instance,
$\rho(x) = x_0$ locally will do).  Then
\begin{eqnarray*}
\sigma(d',\, dx_0) & = & d'(\rho
\phi) \biggr | _{x = (0,x')} \\
 & = & - \sum_{j=0}^N D_j ( \rho \phi_j)
\biggr |_{x_0 = 0} \\
  & = & - \phi_0 .
\end{eqnarray*}
Thus we may write, for $\phi = \sum_{|I| = q} \phi_I \,
dx^I$,
$$
\l du, \phi \r_1 = \l
u, d' \phi \r_1  +
  \int_{b\rnp} \l u,
\sigma (d',\, dx_0)\phi \r_1 \, dx' .
$$
In the next Proposition we describe
the domain of $d^*$ in the
topology of $W^s (\rnp)$ in case $s\in\NN$.  Essentially the
same computations as in the case $s=1$
prove the following result.

\begin{proposition}  \sl
Denote by $d^*$ the Hilbert space adjoint of $d$  in the
topology of $W^s (\rnp)$ for $s=1,2,\dots$.  
Let $\phi\in\bw^{1}_0 (\overline{\rnp})$.
Then $\phi$ lies in $\bigl(\dom d^* \bigr) 
\cap \bw^1(\overline{\rnp})$ 
if and only if
$$
D_0^s \phi_0 \biggr |_{x_0 = 0} 
\equiv \pd{^s}{x_0^s} \phi_0 \biggr |_{x_0 = 0} 
\equiv 0 .
$$
\end{proposition}
\bgpf
We only indicate the changes in the proof for the case $s=1$. 
We have that
$$
\langle du,\phi \rangle _s = \langle u,
d'\phi \rangle _s
- \sum_{j=0}^{s-1} \sum_{|\alpha'| \leq s - j}
       \int_{b\RR^{N+1}_+} D_0^j D^{\alpha'} u
     \overline{D_0^j D^{\alpha'}\phi_0} \, dx' .
$$
Clearly the mapping
$$
u \longmapsto \langle u, d'\phi \rangle _s
$$
is bounded on $W^s(\RR^{N+1}_+)$.
Then, using the trace theorem, we
want to show that for any $j$, $0 \leq j \leq s-1$,
and $|\alpha'| \leq s-j$,
the mapping
$$
u \longmapsto I_j \equiv \int_{b\rnp} D_0^j
D^{\alpha'} u \ \overline{D_0^j D^{\alpha'}
\phi_0} \, dx'
$$
is bounded in the $W^s(\rnp)$ norm.
This establishes the containment
${\cal E} \subseteq {\cal D}.$
 
In order to establish that ${\cal D} \subset
{\cal E}$,  suppose that
$\phi \in {\cal D} \equiv \text{dom}\, d^*
\cap \bigwedge\nolimits^1_0(\rnp)$
and that
$$
D_0^s \phi_0 (0,x') \equiv \frac{\partial^s}{\partial x_0^s}
\phi_0 (0,x')
$$
is not identically zero.  
Setting
$$
u_\epsilon(x_0,x') = x_0^{s-1} (x_0 +
\epsilon)^{3/4} \chi(x_0,x') ,
$$
and arguing as before, it is easy to see that
$$
\|u_\epsilon\|_{W^s(\rnp)} \leq C \qquad \qquad
\forall\epsilon,\ 0 < \epsilon \leq \epsilon_0
                    \leqno {(i)}
$$
and
$$
\sum_{j=0}^s \sum_{|\alpha'| \leq s-j}
\int_{\RR^N} D_0^j D^{\alpha'} u_\epsilon
  \overline{D_0^j D^{\alpha'} \phi_0} \biggr |_{x_0 = 0}
   \, dx' \, \sim \, \frac{1}{\epsilon^{1/4}}
 \qquad \quad \text{as} \ \ \epsilon \rightarrow 0 .
\leqno {(ii)}
$$
These two facts show that $\phi$ cannot be in $\dom d^*$.
\endpf

Finally, we also would like to observe that the \bvp s
(\ref{bvp-half-space})
and (\ref{bvp-domain}) are in
fact  pseudodifferential \bvp s.  
This kind of  \bvp s has been
studied by several authors (see \cite{BDM3}, \cite{ESK}, 
\cite{GRU}, and \cite{RESC}), as we indicate in the next 
section.

\section[Boutet de Monvel-Type Analysis of the
Problem]{Boutet de Monvel-Type Analysis of the
Boundary Value Problem}
 
\setcounter{equation}{0}

The \bvp s  
(\ref{bvp-half-space})
and (\ref{bvp-domain}) can be seen as an instance of a very general 
theory developed by Boutet de Monvel and several other authors.  
In this section we
give a short summary of this theory and of its applications to our
problem.  Such a general approach mainly produces results of local
character in which the data are often assumed to 
be more regular then
we can afford.  This is why we actually choose a different, 
and more direct, approach to solving our \bvp s.
 
For maximum generality, let $\Omega$ be either a
smoothly bounded domain in $\RR^{N+1}$ or a half space in
$\RR^{N+1}.$  Define ${\cal A}$ to be the system
\begin{equation}
{\cal A}= \left ( \begin{array}{lr}
                     P_\Omega + G & K \\
                        T         & S
                   \end{array}
           \right ) :
           \left ( \begin{array}{c}
                   C^\infty_0(\overline{\Omega}) \\
                   C^\infty_0(b\Omega)^M
                   \end{array}
           \right )
                     \longrightarrow
               \left ( \begin{array}{c}
                   C^\infty(\overline{\Omega}) \\
                   C^\infty(b\Omega)^{M'}
                     \end{array}
               \right ) .        \label{A}
\end{equation}
The components of this formula are as follows:
 
1)  The numbers $M, M'$ are non-negative integers, 
$C^{\infty}_0 (b\Omega)^M$
is the Cartesian product of $M$ copies of
$C^{\infty}_0 (b\Omega)$, and $C^{\infty} (b\Omega)^{M'}$
is the Cartesian product of $M'$copies of
$C^{\infty} (b\Omega)$;
 
2)  The operator $P_\Omega$ is a pseudodifferential
operator taking functions defined on $\Omega$ into
functions defined on $\Omega.$  [For more on
pseudodifferential operators, see \cite{KR1}.]  We will require
$P_\Omega$ to have the ``transmission property'' (to
be defined below);
 
3)  The operator $G$ is a ``singular Green operator.''
It is defined on functions on $\Omega$, taking values
in the set of functions defined on $\Omega;$
 
4)  The operator $K$ is a ``Poisson operator''.  It
takes ($M$-tuples of) functions on the boundary to functions on the
domain;
 
5)  The operator $T$ is a ``trace operator''.  Let
$\gamma_\ell$ be the classical trace operator of
order $\ell$,
\begin{eqnarray*}
u & \longmapsto &
\frac{\partial^\ell u}{\partial x_0^\ell} \biggr
|_{x_0=0} \\
 C_0^\infty(\overline{\Omega}) & \longrightarrow &
C_0^\infty(b\Omega) .
\end{eqnarray*}
By definition, a trace operator $T$ is the
composition of a pseudodifferential
operator $S'$ on the boundary and some
restriction operator $\gamma_\ell$: $T = S' \circ \gamma_\ell.$
 
6)  The operator $S$ (to repeat) is a
pseudodifferential operator on the boundary.
 
Notice that our problem
\begin{equation*}
\begin{cases}
\displaystyle{
(- \btu + G) u  =  f } & \on \rnp \\
\displaystyle{
\pd{^2 u}{x_0^2} (0,x')  =  0} & \on \{ 0 \} \times \RR^N 
\end{cases} .
\end{equation*}
is a special case of the above setup,
with $M = 0, M' = 1.$
We would take $K = 0$, $S = 0$ and consider
$$
{\cal A}= \left ( \begin{array}{c}
                     P_\Omega + G \\
                           T
                 \end{array}
           \right ) :
     C^\infty_0(\overline{\Omega}) \longmapsto
     \left(
     \begin{array}{c}
              C^\infty(\overline{\Omega}) \\
                C^\infty(b\Omega)
       \end{array}
       \right) .
$$
We shall construct a parametrix ${\cal R}$
for this operator
which will have the form
$$
{\cal R} = (Q \ \ K ) :
                \left(
                \begin{array}{c}
                   C^\infty_0(\overline{\Omega})  \\
                   C_0^\infty(b\Omega)
                        \end{array}
                \right)
  \longrightarrow   C^\infty(\overline{\Omega}) .
$$
Note here that $M=1$ and $M' = 0$.
 
The classical boundary value problems (such
as the Dirichlet or Neumann
problems) have the form
$$
{\cal A}= \left ( \begin{array}{c}
                         D \\
                         T
                   \end{array}
           \right ) :
 C_0^\infty(\overline{\Omega}) \longrightarrow
                       \left ( \begin{array}{c}
          C^\infty(\overline{\Omega})   \\
                    C^\infty(b\Omega)^{M'}
                               \end{array}
                       \right ) ,
$$
where $D$ is a differential operator and $M'$
is the number of boundary conditions.
 
\setcounter{theorem}{0}
 
\begin{definition}   \rm
A pseudodifferential operator  $P$ of order $d$, defined on
$\RR^{N+1}$, 
is said to have the
{\it transmission property} on $\Omega$ if the mapping
$$
u \longmapsto RP E u \stackrel{\text{def}}{=} P_\Omega u
$$
maps $W^s(\Omega)$ continuously to $W^{s-d}(\Omega)$
for all $s > - 1/2$.  Here $E$ is
the operator given
by extending a function on $\overline{\Omega}$
to all of $\RR^N$ by
setting it equal to zero on the complement.
Also $R$ is the operator
of restriction to $\Omega.$
\end{definition}
 
For our purposes, it is enough to note that
any partial differential
operator with smooth coefficients has the
transmission property.  This
follows by inspection.
 
\begin{definition}   \rm
A {\it Poisson operator} $K$ of order $d$ is
an operator defined by the formula
$$
(Kv) (x_0,x') =  \int_{\RR^N} e^{2\pi ix'\cdot \xi'}
        \tilde{k}(x_0,x',\xi')
    \widehat{v}(\xi') \, d\xi' ,
$$
where the symbol $\tilde{k}$ satisfies the estimates for 
$|\xi'|\ge 1$ given by 
$$
\left \| x_0^\ell D_0^{\ell'}
D_{x'}^\beta D_{\xi'}^\alpha \tilde{k}
 (\, \cdot \, , x' , \xi') \right
\|_{L^2(\RR_+)} \leq c(x')
    \bigl ( 1 + |\xi'|^2
\bigr )^{(1/2)(d - 1/2 - \ell - \ell' - |\alpha|)} .
$$
Here $c(x')$ is a continuous function on
$\RR^N$ that depends on
$\ell, \ell', \alpha, \beta.$
\end{definition}
 
One can show that a Poisson operator $K$ of order $d$
is continuous as a mapping from $W^s_{\text{comp}}(\RR^N)$ 
to $W^{s-d +1/2}_{\text{loc}}(\overline{\rnp}).$
By $W^s_{\text{comp}}(\Omega)$ we mean those elements
of $W^s(\RR^N)$ with compact support in $\Omega.$ 
 
It is a standard fact that the Poisson
integral for the upper half space is a Poisson operator,
according to the above definition, of order 0; that is it
maps
$W^s_{\text{comp}}(\RR^N)$ (where $\RR^N = \partial \RR^{N+1}$)
to $W^{s+1/2}_{\text{loc}}
(\overline{\rnp}).$
 
\begin{definition}   \rm
A {\it trace operator} of order
$d \in \RR$ and of {\bf class} $r \in \NN$
is an operator of the form
$$
Tu = \biggl( \sum_{0 \leq \ell
\leq r-1} S_\ell \gamma_\ell \biggr) u ,
$$
where $\gamma_\ell$ is the standard
trace operator of order $\ell$ and
each $S_\ell$ is a pseudodifferential operator  on $\RR^N$
of order $d - \ell.$
\end{definition}
 
The notion of ``class'' is a natural
artifact of dealing with restriction
operators in the context of Sobolev
spaces. Recall that the classical
restriction theorems sending elements
of $W^s(\rnp)$ to elements
of $W^{s-1/2}(\RR^N)$ are only valid
when $s > 1/2.$  The idea of class
addresses the necessary lower bound for
$s$ in theorems such as this.
 
\begin{definition} \rm
A {\it singular Green operator} of {\it order}
$d \in \RR$ and {\it class}
$r \in \NN$ is an operator of the form
$$
G = \sum_{0\leq \ell \leq r-1} K_\ell \gamma_\ell ,
$$
where each $K_\ell$ is a Poisson
operator of order $d-\ell$ and each
$\gamma_\ell$ is a standard trace operator
of order $\ell.$
\end{definition}
 
It can be shown that a singular Green
operator of order $d$ and class
$r$ is a continuous operator
$$
G: W^s_{\text{comp}}(\overline{\rnp}) \longrightarrow
W^{s-d}_{\text{loc}}(\overline{\rnp})
$$
as long as $s > r - 1/2.$
 
\begin{definition} \rm
Let ${\cal A}$ be a system as in (\ref{A}).
We define the associated {\it boundary
system} by
$$
{\bf A} = \left ( \begin{array}{lr}
 p_\Omega(0,x',D_0,\xi') + g(x',D_0,\xi')
& k(x',D_0,\xi') \\
   t(x',D_0,\xi') & s(x',\xi')
                  \end{array}
          \right ) .
$$
\end{definition}
Notice that, in this definition,
we fix $x', \xi' \in \RR^N$, and
the symbol $p_\Omega$
of $P_\Omega$ is taken at the point
$(0,x') \in b \rnp.$
Finally, all the operators in the
display act only in the $D_0$ slot.
 
Further observe that
$$
p_\Omega(0,x',D_0,\xi') + g(x',D_0,\xi'):
{\cal S}(\overline{\RR_+})
 \longrightarrow {\cal S}(\overline{\rnp})
$$
and that
\begin{eqnarray*}
  k(x',D_0,\xi'): \CC & \longrightarrow &
{\cal S}(\overline{\rnp}) \\
 a & \longmapsto & a k(x', x_0,\xi') .
\end{eqnarray*}
Furthermore, the boundary symbol operator of a
trace operator is an operator
$$
t(x',D_0,\xi') u = \sum_{0 \leq \ell \leq r-1}
s_\ell(x',\xi') \gamma_\ell u
$$
that maps ${\cal S}(\overline{\RR_+})$ into $\CC.$

A system of the type (\ref{A}) is called {\em elliptic}
if there exists a second system
${\cal R}$ of type (\ref{A}) such that
$$
{\cal R}{\cal A}+ I \qquad \quad
\text{and} \qquad \quad {\cal A}{\cal R} + I
$$
are negligible operators; here
``negligible'' means that all terms
in the system send distributions with compact
support into $C^\infty$
functions.  In the case of the half space this definition must be  
considered as ``local", that is fixing a compact set in
$\overline{\rnp}$. 
 
Now we have
 
\begin{theorem}[\text{[BDM3, Theorem 5.1]}]\label{BdM}  \sl
 \sl A system $a$ of type (\ref{A}) is elliptic if and
only if both its boundary and its interior symbols are
invertible.
\end{theorem}

In the case of (\ref{bvp-half-space}),
$$
{\cal A}= \left ( \begin{array}{c}
                  - \btu + G \\
                         T
                   \end{array}
           \right) ,
$$
where
$T \stackrel{\text{def}}{=} (I - \btu')^{-1/2} \gamma_2$.
Then
we return to our analysis of the problem
\begin{equation}\label{BP}
\begin{cases}
\displaystyle{
- \btu u + G u  =  f  } & \on \rnp \\
\displaystyle{
  Tu  =  0 } & \on \{0\} \times \RR^N 
\end{cases}
\end{equation}
Here, for convenience, we have replaced our
usual boundary
operator $(\p^2 /\p x_0^2 )$
with the (equivalent) operator $T$. 

Our aim is to show that this is a
pseudodifferential boundary value
problem of type (\ref{A}) and 
to apply Theorem \ref{BdM}
to our system.  
 
First observe that $- \btu$, being a partial
differential operator
of order 2 with smooth coefficients, certainly
possesses the transmission
property.
 
Next, let us analyze $G$.  
Recall that $G = \widetilde{\cal K} \gamma_1$,
where $\widetilde{\cal K}$
has symbol
$$
K_{x_0} (\xi')\equiv
k(x_0,\xi') = -\Sx  e^{-\Sx x_0} .
$$
In order to see that $G$ is a singular
Green operator, we need to estimate
\begin{multline*}
\bigl \| x_0^\ell D_0^{\ell'}
D_{\xi'}^\alpha k(\, \cdot \, , \xi')
  \bigr \|_{L^2_{x_0} (\RR_+)}^2 \\
   =  \int_0^\infty \biggl |
x_0^\ell D_{\xi'}^\alpha \left [
\bigl ( - \sqrt{1 + 4\pi^2 |\xi'|^2}
\bigr )^{\ell' + 1}
 e^{-\sqrt{1 + 4\pi^2|\xi'|^2} x_0}
\right ] \biggr |^2 \, dx_0
\end{multline*}
 
The calculations made in Proposition \ref{tilde-K} 
showed that $\widetilde{\cal K}$ has order 1.
In fact $\widetilde{\cal K}$ is essentially
a single derivative of the Poisson operator
(which has order zero) so the
assertion is plausible.
It follows that
$$
\text{order}\,(G) = \text{order}\,(\widetilde{\cal K}) +
\text{order}\,(\gamma_1) = 2
$$
and
$$
\text{class}\,(G) = \text{order}\,(\gamma_1) + 1 = 2 .
$$
 
Next we observe that $T$ is a trace operator
of order 1 and class 3.
Indeed, $T \stackrel{\text{def}}{=}
(I - \btu'')^{-1/2} \gamma_2$ so that
$$
\order(T) = \order(I - \btu')^{-1/2}
+ \order(\gamma_2) = - 1 + 2 = 1
$$
and
$$
\class(T) = \order(\gamma_2) + 1 = 3 .
$$
 
The interior symbol of this system is
just the symbol of $- \btu$
which is plainly invertible.  The boundary
symbol is
$$
\left ( \begin{array}{c}
 \left ( \tpx^2 -
\frac{\partial^2}{\partial x_0^2} \right )
 + \left ( \frac{\partial}{\partial x_0}
e^{-\sqrt{1 + \tpx^2} x_0} \right ) \gamma_1 \\
          (1 + |\xi'|^2)^{-1/2} \gamma_2
        \end{array}
\right ) : {\cal S}(\overline{\RR_+})
\longrightarrow
  \left ( \begin{array}{c}
     {\cal S}(\overline{\RR_+}) \\
              \CC
          \end{array}
  \right ) .
$$
The demand that this system be invertible
means that the pseudodifferential
system
$$
\begin{cases}
\displaystyle{
v'' (x_0) - \left ( \frac{\partial}{\partial x_0}
      e^{-\sqrt{1 + 4\pi^2|\xi'|^2} x_0} \right )
          v' (0) - \tpx^2 v(x_0)  =  \psi} \qquad & \on \RR \\
\displaystyle{
  (1 + \tpx^2)^{-1/2}
v''(0)  =  \alpha } 
\end{cases}
$$
has one and only one solution
$v \in {\cal S}(\overline{\RR_+})$
for all $\psi \in {\cal S}(\overline{\RR_+})$
and $\alpha \in \CC$
and for $|\xi'| \geq c > 0$ fixed.
 
Straightforward and unenlightening calculations now prove the 
following result.
\begin{proposition}   \sl
For $|\xi'| \geq c > 0$ fixed,
the preceding system of equations
has a unique solution for every choice of
$(\psi, \alpha)
   \in {\cal S}(\overline{\RR_+}) \times \CC.$
\end{proposition}
As a consequence we have:
\begin{theorem}  \sl
Consider the \bvp\ given by (\ref{BP}):
\begin{equation*}
\begin{cases}
\displaystyle{
 (-\btu + G) u  =  f } & \on \RR_+^{N+1} \\
\displaystyle{
Tu  =  g } & \on \{0\} \times \RR^N 
\end{cases} .
\end{equation*}
Let $s >5/2$, and suppose that $\supp f$, $\supp g \ss B(0,R)$, for
some $R>0$.  Moreover, let 
$ \eta \in C^{\infty}_0
(\overline{\rnp})$. 
Then there exists a
finite dimensional subspace ${\cal L}$
of $W^s_{\text{comp}}(\rnp) \times 
W^{s-1/2}_{\text{comp}}(\RR^N)$
such that if $(f,g) \in {\cal L}^\perp$
then the problem (\ref{BP}) has a unique
solution $u$ satisfying
$$
\|\eta u\|_{W^{s+2}(\rnp)} \leq C \bigl (
\| f\|_{W^s(\rnp)} +
\| g\|_{W^{s-1/2}(\RR^N)} \bigr ) .
$$
The constant $C$ in the above estimate depends on $R,s$,
the test 
function 
$\eta$, and on the spatial dimension $N$.
\end{theorem}
Our aim is to refine this theorem in order to allow $s>1/2$ (instead of
$s>5/2$) and to be more precise about the local-global aspects of the
problem. Moreover, we do so by
determining the solution explicitly, and by describing the
compatibility conditions.  
 
\section{The Explicit Solution in the Case of Functions}
\setcounter{equation}{0}
 
Our goal in the present section is to prove theorems about existence and
regularity of the solution of the \bvp .
We shall determine the solution of our boundary value problem
explicitly. This will allow us to give precise estimates.
 
Thus we concentrate on the problem
\begin{equation} 
\begin{cases}
{\displaystyle
-\btu u + G u  =  f }\qquad & \on  \rnp \\
{\displaystyle
\pd{^2 u}{x_0^2}(0,\cdot)  =  0 }
		& \on \RR^N 
\end{cases} 
\label{BV2}
\end{equation}
for $u  \in W^{r}(\RR_+^{N+1})$, where $r>5/2$.
Notice that we need this
limitation on $r$ in order to
guarantee that $\p^2  u/\p x_0^2(0,\cdot)$
makes sense and belongs to $L^2 (\RR^N)$ (at least).
Here, for $x_0 >
0$ and $\xi' \in \RR^N$, we have
$$
\widehat{Gu}(x_0,\xi') = -
\sqrt{1 + |2\pi\xi'|^2} e^{-\sqrt{1 +
|2\pi\xi'|^2}x_0}
 \frac{\partial \widehat{u}}{\partial x_0} (0,\xi').
$$
Recall that we  write
$\widehat{w}(\xi')$ to denote
the
partial Fourier transform of $w$ in $\RR^N$, with respect to
the $x'$-variables.

For greater flexibility we will solve the \bvp\
\begin{equation}
\begin{cases}
{\displaystyle
-\btu u + G u } =  f &\qquad \on  \rnp \\
{\displaystyle
\frac{\partial^2 u}{\partial x_0^2}(0,\cdot)}  =  h
		&\qquad\on \RR^N 
\end{cases} 
\label{BV3}
\end{equation}
 
In what follows 
we shall denote by ${\cal G}$  the classical Green's
function for the Laplacian
on $\rnp$; also ${\cal P}$ denotes the standard
Poisson operator on $\rnp$.
Recall (see \cite{GAR}) that, for $N\ge2$,
\begin{multline*}
{\cal G}(x,y)  = 
\frac{1}{(N-1)\omega_{N+1}} \biggl[ 
\frac{1}{\bigl( (x_0 - y_0)^2
     	+ |x' - y'|^2 \bigr)^{(N-1)/2}} \\
-\frac{1}{\bigl( (x_0 +y_0)^2 +|x'-y'|^2 \bigr)^{(N-1)/2}}
\biggr] ,
\end{multline*}
where $\omega_{N+1}$ denotes the surface measure of the 
$N$-dimensional unit sphere in 
$\RR^{N+1}$.
If $N=1$ then
$$
{\cal G}(x,y)  = \frac{1}{2\pi} \log
\frac{(x_0 +y_0 )^2 +(x_1 -y_1 )^2}{(x_0 -y_0 )^2 +(x_1 -y_1 )^2} .
$$
Notice that
\begin{equation}\label{tang-FT-calG}
\widehat{{\cal G}f} (x_0,\xi')
=\frac12 \frac{1}{\tpx} \int_0^\infty \biggl(
e^{-\tpx\cdot|x_0 -y_0|} -e^{-\tpx(x_0 +y_0 )} 
\biggr) \hat f(y_0,\xi')\, dy_0 .
\end{equation}
Moreover,
we shall adopt the following notation:
\begin{equation}\label{m_1}
m_1(\xi') = \frac{1+\tpx^2 +\tpx\Sx}{1+2\tpx^2 +\tpx\Sx},
\end{equation}
and
\begin{equation}\label{m_2}
m_2 (\xi')= \left[ \frac{\tpx +\Sx}{1+2\tpx^2 +\tpx\Sx}\right] 
= \frac{1}{\sqrt{1 + |2\pi \xi'|^2}} m_1(\xi') .
\end{equation}
We now construct the solution to the \bvp\ explicitly.
\begin{proposition}\label{exp-sol}    \sl
If $u$
is a solution of the \bvp\ (\ref{BV3})
in $W^{r}(\rnp)$ with $r>5/2$, then $u$ is given by
the following formula
\begin{eqnarray}
u(x_0,x')
&  = &  \G f(x_0,x') +\G\left( 
	\int_{\RR^N} e^{-\Sx x_0}m_1 (\xi') \right. \nonumber \\ 
&    & \quad \left. \times \left[
\tpx\int_0^{+\infty} e^{-\tpx y_0}\hat f(y_0,\xi')\, dy_0 
-\hat h(\xi') -\hat f(0,\xi') \right] 
e^{2\pi ix'\cdot\xi' }\, d\xi' \right)\nonumber \\
& &\quad +{\cal P}\left( 
\frac{\hat h(\xi')+\hat f (0,\xi')}{1+2\tpx^2 +\tpx\Sx} 
\right)\bcheck \nonumber \\
& & \quad + {\cal P}\left( \frac{m_1(\xi')}{\tpx} \bigl( 
\hat h(\xi')+\int_0^{+\infty} \hat f (y_0,\xi') 
e^{-\tpx y_0} \, dy_0 \bigr) \right)\bcheck  
\label{u-statement}
\end{eqnarray}
\end{proposition}
\bgpf
Our first task is to reduce the order of the operator appearing in
the boundary condition by using the equation on the domain.
Thus the problem (\ref{BV3}) is equivalent to
\begin{equation*}
	\begin{cases}
-\btu u + Gu  =  f &   \on\RR_+^{N+1}\\
-\btu' u(0,\cdot)  =  f(0,\cdot)+h -Gu(0,\cdot)&  \on\RR^N 
	\end{cases}
\end{equation*}
for $u\in  W^{r}(\rnp)$ with 
$r>5/2$,
where $\btu' \stackrel{\text{def}}{=}
\sum_{j=1}^N \partial^2/\partial x_j^2 .$
The problem (\ref{BV3}) can be rewritten as
\begin{equation*}
	\begin{cases}
-\btu u(x_0,x')   =  F(x_0,x')\qquad & \on\rnp\\
 u(0,x')  =  g(x') &\on\RR^N 
 	\end{cases} 
\end{equation*}
(again for 
$u\in W^{r}(\rnp),\, r>5/2$)
where we have set
\begin{equation}
F(x_0,x')
 =   f(x_0,x')        
+ \int_{\RR^N}\Sx
e^{-\Sx x_0}\frac{\p\hat u}{\p x_0}(0,\xi')e^{2\pi ix'\cdot\xi'}\,
                d\xi' , \label{F}
\end{equation}
and
\begin{equation}\label{g}
g(x')  = {\cal N} \bigl( h+ f(0,\cdot)\bigr) (x')
+{\cal N} \left(
\bigl( \Sx \frac{\p\hat u}{\p x_0}(0,\cdot)\bigr)
\bcheck  (0,x') \right) ,
\end{equation}
with ${\cal N}$ being the Newtonian potential in $\RR^N$.
By the classical theory we then obtain that the solution we seek 
can be written as 
\begin{eqnarray*}
u(x_0,x')
& = & \int_{\rnp} \G (x_0,y) F(y)\, dy
        +\int_{\RR^N}P_{x_0}(x',y')g(y')\, dy' \\
&\equiv & (\G F) (x_0,x')+ {\cal P}g(x_0,x') .
\end{eqnarray*}
Here $P_{x_0}$ denote the Poisson kernel in $\RR^N$.
Now we need to
isolate $u$ on the left hand side (note that
$F$, $g$ are defined in terms of $u$).  In order to do this we
compute 
the derivative with respect to $x_0$ of the above equation, take
the partial Fourier transform, and evaluate
at $x_0 =0$.  Using definitions (\ref{F}) and (\ref{g}) 
we have
\begin{eqnarray*}
\left. \frac{\p\hat u}{\p x_0}\right|_{x_0 =0}
& = & \left. \left[ \biggl( \frac{\p (\G F)}{\p x_0}\biggr) \,
    {}^{\widehat{}}  +\biggl(\frac{\p {\cal P}g}{\p x_0} \biggr) 
     \, {}^{\widehat{}} \,
      \right]\right|_{x_0 =0} \\
& = & \int_{0}^{+\infty}({\cal P}F)\, \widehat{}\, (y_0,\xi')\,
                dy_0 -\tpx\hat g(\xi') \\
& = & \int_0^{+\infty}e^{-\tpx y_0}\, \widehat F(y_0 ,\xi')\, dy_0
        -\tpx\hat g(\xi) \\
& = & \int_0^{+\infty}e^{-\tpx y_0} \hat f(y_0 ,\xi')\, dy_0
                                                 \\
&   & \qquad  \Sx +\frac{\p\hat u}{\p x_0}(0,\xi')
	\int_0^{+\infty} e^{-(\Sx +\tpx )y_0}\, dy_0  \\
&   & -\tpx\left[ \frac{\hat h (\xi')}{\tpx^2} 
	+\frac{\hat f (0,\xi')}{\tpx^2}
        +\frac{\Sx (\p\hat u/\p x_0 )(0,\xi')}{\tpx^2}\right]\\
& = & \int_0^{+\infty}e^{-\tpx y_0}\hat f(y_0,\xi')\, dy_0
        -\frac{\hat h(\xi')}{\tpx} 
	-\frac{\hat f(0,\xi')}{\tpx} \\
&  & \qquad + 
\left[\frac{\Sx}{\Sx +\tpx}-\frac{\Sx}{\tpx}\right]
        \frac{\p\hat u}{\p x_0}(0,\xi') .
\end{eqnarray*}
 
Hence
\begin{multline*}
\frac{\p\hat u}{\p x_0}(0,\xi')
 \left[ 1 -\frac{\Sx}{\Sx +\tpx}  + \frac{\Sx}{\tpx}\right] \\
= \int_0^{+\infty}e^{-\tpx y_0}\hat f(y_0,\xi')\, dy_0
   -\frac{\hat h(\xi')}{\tpx} -\frac{\hat f(0,\xi')}{\tpx} ,
\end{multline*}
that is,
\begin{equation} \label{part-u(0,)}
 \frac{\p\hat u}{\p x_0}(0,\xi')
 = \tpx m_2 (\xi') \left\{
\int_0^{+\infty}e^{-\tpx y_0}\hat f(y_0,\xi')\, dy_0
 -\frac{\hat h(\xi')}{\tpx} -\frac{\hat f(0,\xi')}{\tpx} \right\} .
\end{equation}
Now we have the expression for $\p u/\p x_0$ at $x_0 =0$.
Notice that, for $t>0$,
$$
\left\| \frac{\p u}{\p x_0}(0,\cdot)\right\|_{W^{t}(\RR^N )}
\leq C\bigl( \| f\|_{W^{t+1/2}(\rnp)} + \| h\|_{W^{t-1}(\RR^N)}
\bigr) ,
$$
so that our calculations make sense so far.
Substituting equation (\ref{part-u(0,)}) into (\ref{F}) and 
(\ref{g}), and recalling that $m_1 =\Sx m_2$, we find that
\begin{multline*}
F(x_0,x') =  f(x_0,x') \\
+\int_{\RR^N} \tpx m_1(\xi')e^{-\Sx}
\biggl( \int_0^{+\infty} e^{-\tpx y_0} \hat f(y_0,\xi') \, dy_0 
-\frac{\hat h(\xi')}{\tpx} -\frac{\hat f (x_0,\xi')}{\tpx} \biggr) ,
\end{multline*}
and
\begin{align*}
g(x')
& = {\cal N}\bigl( h+f(0,\cdot) \bigr)\\
& \ +{\cal N} \biggl( \bigl(m_1 (\xi')\tpx \bigl\{ 		
	\int_0^{+\infty} e^{-\tpx y_0} \hat f(y_0,\xi')\, dy_0 
	-\frac{\hat h(\xi')}{\tpx}
	-\frac{\hat f(0,\xi')}{\tpx}\bigr\}\bigr)\bcheck \biggr) \\
& = {\cal N}\bigl( h+f(0,\cdot) \bigr)
 + {\cal N} \biggl( \bigl(\tpx m_1 (\xi') \int_0^{+\infty}
	e^{-\tpx y_0} \hat f(y_0,\xi')\, dy_0 \bigr) \bcheck \biggr)  \\
& \quad -{\cal N}\biggl( \bigl( m_1 \hat h +m_1 \hat f(0,\cdot) \bigr)\bcheck
\biggr) .
\end{align*}

Finally,
\begin{eqnarray}
u(x_0,x')
& = & \G F(x_0,x') + {\cal P}g(x')\nonumber \\
& = & \G f(x_0,x') +\G\left( \int_{\RR^N} e^{-\Sx x_0}m_1 (\xi')
                     \right.\nonumber \\
& & \quad\times  \left. 
\left[
\tpx\int_0^{+\infty} e^{-\tpx y_0}\hat f(y_0,\xi')\, dy_0 
-\hat h (\xi') -\hat f(0,\xi') 
\right] e^{2\pi ix'\cdot\xi' }\, d\xi' \right)\nonumber \\
& & \qquad +{\cal P}\left( 
\frac{\hat h (\xi')+ \hat f (0,\xi')}{1+2\tpx^2 +\tpx\Sx}
\right) \bcheck   \nonumber \\
& & \qquad +{\cal P}\left( \frac{m_1(\xi')}{\tpx} 
\int_0^{+\infty} 
\hat f (y_0,\xi') e^{-\tpx y_0} \, dy_0 \right)\bcheck
. \label{u}
\end{eqnarray}
This gives the explicit expression for the solution $u$.
\endpf
 
\subsection*{The {\em a priori} Estimate}
In this part we are going to prove the
{\em a priori} estimate for 
the solution $u$; then we will turn to
the question of existence.
Now we need the following two lemmas.
\begin{lemma}\label{FACT-1}     \sl
There exists a constant $c$ such that,
for $|\alpha|=2$, and $t\ge 0$,
$$
\| D^\alpha \G w \|_{W^t (\rnp)} \le c \| w\|_{W^t (\rnp)} .
$$
\end{lemma}
\bgpf  For a function $w$ defined on $\rnp$,
let $w_{\text{odd}}$ denote the odd
extension (in the $x_0$ variable) to $\RR^{N+1}$, and 
${\cal N}_{N+1}$
the Newtonian potential in $\RR^{N+1}$.
Let  ${\cal F}$ denote the Fourier transform in $\RR^{N+1}$.
Observe that
$$
(\G w)_{\text{odd}} ={\cal N}_{N+1}  w_{\text{odd}} ,
$$
so that we  have
$$
{\cal F} \left( \pd{^2 (\G w)_{\text{odd}}}{x_i x_j} \right) (\xi) 
= \frac{\xi_i \xi_j}{|\xi|^2} {\cal F}(w_{\text{odd}}) (\xi), 
$$
and the estimate for $t=0$ follows.  
When $t=1,2,\dots$ notice that,
for $x \in \rnp$,
$$
\biggl(\frac{\partial^2 (\G w)}{\partial x_0^2} 
\biggr)_{\text{odd}} (x)
= -w-\btu' \G w (x) .
$$
Then it suffices to consider the case when the derivative is of
degree 
at most $1$ in the $x_0$-direction.  This takes care 
of the fact that
$w_{\text{odd}}$ 
is not differentiable in the $x_0$-direction.  Now the estimate
follows easily for $t$ an integer.  Interpolation gives the result
for all $t\ge 0$.
\endpf 
 
The proof of the next lemma is easy: 
\begin{lemma}\label{FACT-2}     \sl
Let $w$ be defined on $\RR^N$,
$|\alpha|=2$, and $t\ge 0$.  Then there exists a constant $C>0$ such
that
$$
\| D^\alpha {\cal P} w\|_{W^{t}(\rnp)}^2
\le C \int_{\RR^N} \tpx^3 (1+\tpx)^{2+t} |\hat w(\xi')|^2 \, d\xi' . 
$$
\end{lemma}
 
\begin{theorem}\label{apr-est}  \sl
Let $f\in W^s (\rnp)$, with $s>1/2$.
If the \bvp\  (\ref{BV3})
admits the solution $u$ given by Proposition
\ref{exp-sol}, then it satisfies the estimate
\begin{equation}\label{a priori-est}
\| u\|_{W^{s+2} (\rnp)} \leq C\left\{ \| f\|_{W^s (\rnp)}
 +\| h\|_{W^{s-1/2} (\RR^N)}+\| u\|_{W^{1} (\rnp)} \right\} . 
\end{equation}
\end{theorem}
\bgpf
In order to obtain {\em a priori} estimates we use the following
standard fact (see \cite{LIM} Theorem 9.7).  If $s\ge 0$, then
$$
\| u\|_{W^{s+2}(\rnp)}
\le c\left( \sum_{|\alpha|=2}
\left\| \pd{^{|\alpha|}u}{x^\alpha} \right\|_{W^s(\rnp)}
+\| u\|_{W^1 (\rnp)} \right) .
$$
 
Now recall the expression for the solution given by (\ref{u})
\begin{eqnarray*}
u(x_0,x')
& = & \G f(x_0,x') +\G\left( \int_{\RR^N} e^{-\Sx x_0}m_1 (\xi') 
                                        \right. \\
& & \times  \left. \left[
\tpx\int_0^{+\infty} e^{-\tpx y_0}\hat f(y_0,\xi')\, dy_0 
-\hat h (\xi') -\hat f(0,\xi') \right] 
e^{2\pi ix'\cdot\xi' }\, d\xi' \right) \\
& & \qquad +{\cal P}\left( \frac{\hat h(\xi') 
	+\hat f (0,\xi')}{1+2\tpx^2 +\tpx\Sx}
\right) \bcheck  \\
& & \qquad +{\cal P}\left( \frac{m_1(\xi')}{\tpx} \int_0^{+\infty} 
\hat f (y_0,\xi') e^{-\tpx y_0} \, dy_0 \right)
\bcheck  \\
& \equiv & \G f +\G f_1 +{\cal P}g_1 +{\cal P}g_2 .
\end{eqnarray*}
It follows that
\begin{multline*} 
\| u\|_{W^{s+2}(\rnp)}
 \le  c \| u\|_{W^1 (\rnp)}  
+ c \sum_{|\alpha|=2}
\biggl( \left\|\pd{^{|\alpha|}\G f }{x^\alpha} 
		\right\|_{W^s(\rnp)} \\
+ \left\| \pd{^{|\alpha|}\G f_1 }{x^\alpha} \right\|_{W^s(\rnp)}
+ \left\| 
\pd{^{|\alpha|}{\cal P}g_1}{x^\alpha}\right\|_{W^s(\rnp)} 
 + \left\|\pd{^{|\alpha|}{\cal P} g_2}{x^\alpha}
        \right\|_{W^s(\rnp)}\biggr) 
\end{multline*}
 
Lemma \ref{FACT-1} implies that, for $|\alpha|=2$,
$$
\| D^\alpha \G f \|_{W^{t}(\rnp)} \le c \| f \|_{W^t (\rnp)} ,
$$
and
$$
\| D^\alpha \G f_1 \|_{W^{t}(\rnp)} \le c \| f_1 \|_{W^t (\rnp)} \, .
$$
Notice that $f_1=\widetilde{\cal K}(h_1)$, where
$$
\widehat{h_1} (\xi') =  m_2 (\xi')
\left[
\tpx\int_0^{+\infty} e^{-\tpx y_0}\hat f(y_0,\xi')\, dy_0 -\hat
h (\xi') -\hat f(0,\xi') \right] .
$$
Then, for $s >1/2$,
\begin{eqnarray*}
\| f_1 \|_{W^s (\rnp)}
& \le & c \| h_1 \|_{W^{s+1/2}(\RR^N )} \\
& = & c\biggl\{ \int_{\RR^N} (1+|\xi'|^2)^{s+1/2} |m_2(\xi')|^2 \\ 
& &  \times \left| \tpx\int_0^{+\infty} e^{-\tpx y_0}
\hat f(y_0,\xi')\, dy_0 
-\hat h(\xi') -\hat f(0,\xi') \right|^2 \, d\xi'
\biggr \}^{1/2} \\
& \le & c\left( \| h \|_{W^{s-1/2} (\RR^N)}  
	+ \| f(0,\cdot)\|_{W^{s-1/2}(\RR^N)}
        + \| f\|_{W^{s}(\rnp)}\right) \\
& \le & c \bigl( \| f\|_{W^s (\rnp)} +\| h \|_{W^{s-1/2} (\RR^N)} 
	\bigr) .
\end{eqnarray*}
 
Hence it suffices to consider the $W^t (\rnp)$ norm  of
any second derivative of ${\cal P}g_1$
and ${\cal P}g_2$.
For $s>1/2$, by Lemma \ref{FACT-2} we have
\begin{eqnarray*}
\lefteqn{\left\| D^\alpha {\cal P}
        \biggl( \frac{ \hat h (\xi')
	+\hat f (0,\xi')}{1+2\tpx^2 +\tpx\Sx} \biggr)
                \right\|_{W^s (\rnp)} \text{ \ \ \ \ }} \\
& \le & c\left\{ \int_{\RR^N} \tpx^{3} (1+\tpx^2)^{s-2}
	\bigl( |\hat h (\xi')|^2 + |\hat f (0,\xi')|^2 \bigr)
	\, d\xi' \right\}^{1/2} \\
& \le & c \bigl\{ \| f(0,\cdot)\|_{W^{s-1/2}(\RR^N)} 
+ \| h \|_{W^{s-1/2} (\RR^N)} \bigr\} \\
& \le & c\bigl\{ \| f\|_{W^s (\rnp)} 
\| h \|_{W^{s-1/2}(\RR^N)} \bigr\} , 
\end{eqnarray*}
Moreover, recalling the definition of $m_1(\xi')$ (see (\ref{m_1})),
we have
\begin{eqnarray*}
\lefteqn{
\left\| D^\alpha {\cal P}
\biggl( \frac{m_1(\xi')}{\tpx} \int_{0}^{+\infty} \hat f (y_0,\xi')
e^{-\tpx y_0 }\, dy_0 \biggr)  \right\|_{W^s (\rnp)} } \\
& \le & c \left\{ \int_{\RR^N}
\tpx^{3} (1+\tpx^2)^{s}\frac{1}{\tpx^2} \bigg|
        \int_{0}^{+\infty}  \hat f (y_0 ,\xi')e^{-\tpx y_0} \, dy_0
        \bigg|^2 \, d\xi' \right\}^{1/2} \\
& \le & c \left\{ \int_{\RR^N} (1+\tpx^2)^{s} \int_0^{+\infty}
|\hat  f(y_0 ,\xi')|^2 \, dy_0 \, d\xi' \right\}^{1/2} \\
& = & c  \| f\|_{W^s (\rnp)} .
\end{eqnarray*}
This proves the {\em a priori} estimate.
\endpf 
\subsection*{The existence Theorem}
Now we turn to the existence statement.  If we prove that the
function $u$ defined by (\ref{u}) belongs to $W^1 (\rnp)$ then, 
by the estimate (\ref{a priori-est}),
the equation (\ref{u}) defines the unique
solution of the \bvp.  The solution belongs to $W^{r} (\rnp)$ with $r>5/2$. 
\begin{theorem}\label{exist}     \sl
If $f \in W^s \cap L^1 (\rnp)$, with $s>1/2$, then the
\bvp\ (\ref{BV2})
has a unique solution $u$, and this satisfies 
\begin{eqnarray*}\label{existence}
\| u\|_{W^{s+2}(\rnp)} \leq C\left\{ \| f\|_{W^s (\rnp)}
        +\| f\|_{L^1 (\rnp)} 
+ \| h \|_{W^{s-1/2}(\RR^N )} 
\right\} 
\end{eqnarray*}
if $N\ge4$. 
If $N=2,3$ and $f$ satisfies
$$
\int f(x)\,  dx =0 , \qquad \qquad f\in L^1 (|x|dx,\rnp) ,
$$
then there exists a unique solution $u$ that satisfies 
the estimate
$$
\| u \|_{W^{s+2} (\rnp)} \le C \left\{ \| f\|_{W^s (\rnp)} 
+ \| f\|_{L^1 (|x|dx,\rnp)} + \| h \|_{W^{s-1/2}(\RR^N )} 
\right\} .
$$
If $N=1$, we assume that $f$ is such that
$$
\int f(x)\, dx =0 , \qquad\int x_i f(x)\, dx =0 ,\ i=0,1
 \qquad f \in L^1 (|x|^2 dx, \RR^2 ).
$$
Then there exists a unique solution $u$ such that
$$
\| u \|_{W^{s+2} (\rnp)} \le C\left\{ \| f\|_{W^s (\rnp)}
+ \| f\|_{L^1 (|x|^2 dx,\rnp)} 
+ \| h \|_{W^{s-1/2}(\RR^N )} \right\} .
$$
\end{theorem}
\bgpf By Theorem \ref{apr-est} it suffices to show that
$$
\| u\|_{W^1 (\rnp)} \le c\left\{ \| f\|_{W^s (\rnp)} +
\| f\|_{L^1 (\rnp)} + \| h \|_{W^{s-1/2}(\RR^N )} \right\}  ,
$$
when $N\ge 4$, or the corresponding estimate 
when $N=1,2,3$, with $f$ satisfying the stated conditions.
Notice that
$$
\| u\|_{W^1 (\rnp)} \le  \| u\|_{L^2 (\rnp)} +
\sum_{i=0}^{N} \biggl \| \pd{u}{x_i} \biggr \|_{L^2 (\rnp)} .
$$
Recall that (\ref{u}) gives the
function $u$ as
$$
u =\G f + \G f_1 +{\cal P} g_1 +{\cal P} g_2 .
$$
We calculate the partial Fourier transform of $u$ using 
the fact that
$$
\biggl( \G\bigl [(e^{-\Sx x_0}) \bcheck \bigr ] \biggr)\, 
\widehat{\mathstrut} 
\, (x_0, \xi') = e^{-\tpx x_0}-e^{-\Sx x_0}.
$$
[Note here that $\bcheck{\mathstrut}$ and $\widehat{\mathstrut}$
do not cancel, since the Green's potential ${\cal G}$ occurs in
between the two operations.  Reference line (5.3).]
We obtain 
\begin{eqnarray*}
\hat u(x_0,\xi')
& = & \widehat{\G f} (x_0,\xi') +
\widehat{\G f_1}(x_0, \xi') + \widehat{{\cal P} g_1 }(x_0, \xi') 
+ \widehat{{\cal P} g_2} (x_0, \xi') \\
& = & \widehat{\G f} (x_0,\xi') + m_1 (\xi')\left[
e^{-\tpx x_0} -e^{-\Sx x_0}\right] \\
& & \qquad \quad \times \,
\biggl[ \tpx \int_0^{+\infty}
e^{-\tpx y_0} \hat f (y_0,\xi') \, dy_0 -\hat h (\xi')
-\hat f (0,\xi') \biggr] \\
& & \quad  + \, \frac{e^{-\tpx x_0} \bigl( \hat h (\xi')
+\hat f (0,\xi')\bigr)}{1+2\tpx^2 +\tpx\Sx}\\
& & \quad + \, e^{-\tpx x_0} \frac{m_1(\xi')}{\tpx}
\int_0^{+\infty} e^{-\tpx y_0} \hat f (y_0,\xi') \, dy_0 \, .
\end{eqnarray*}
Moreover,
\begin{eqnarray*}
\left( \pd{\hat u}{x_0}\right) (x_0,\xi')
& = & \left( \pd{}{x_0}\widehat{\G f} \right) (x_0,\xi') \\
& & \quad +  m_1 (\xi')\left[ \Sx e^{-\Sx x_0}
-\tpx e^{-\tpx x_0}\right] \\
& & \qquad \quad \times
\biggl[ \tpx \int_0^{+\infty}
e^{-\tpx y_0} \hat f (y_0,\xi') \, dy_0 -\hat h (\xi')
-\hat f (0,\xi') \biggr] \\
& & \qquad
-\frac{\tpx e^{-\tpx x_0} \bigl( \hat h (\xi')
+\hat f (0,\xi')\bigr)}{1+2\tpx^2 +\tpx\Sx}\\
& & \qquad - e^{-\tpx x_0} m_1(\xi')
\int_0^{+\infty} e^{-\tpx y_0} \hat f (y_0,\xi') \, dy_0 \, .
\end{eqnarray*}
Then $\| u\|_{W^1 (\rnp)}$ will be estimated once we estimate
$\| \hat u \|_{L^2 (\rnp)}$,
$\bigl \| \tpx \hat u \bigr \|_{L^2 (\rnp)}$,
and $\|(\p u/\p x_0)\, \hat{}\,\|_{L^2 (\rnp)}$.
 
We begin by studying the terms that
arise from $\G f$.  
We extend $f$ to all of $\RR^{N+1}$ as an odd function
$f_{\text{odd}}$ of $x_0$ as we did in the proof of Lemma
\ref{FACT-1}.
We obtain 
$$
\| \G f\|_{L^2 (\rnp)} =\frac{1}{\sqrt 2}
\biggl \|  \frac{1}{\tpx^2} 
{\cal F} \bigl( {\cal N}f_{\text{odd}} \bigr) 
\biggr \|_{L^2 (\RR^{N+1})}
$$
and
$$
\left \| \pd{}{x_i} \G f \right \|_{L^2 (\rnp)} =\frac{1}{\sqrt 2}
\left \| \pd{}{x_i}{\cal F}f_{\text{odd}}
 \right \|_{L^2 (\RR^{N+1})} .
$$
{\bf Case $N\ge4$.}  We have
\begin{eqnarray*}
\| \G f\|_{L^2 (\rnp)} 
& = & c\cdot \| {\cal F} ({\cal N}_{N+1}
       f_{\text{odd}}) \|_{L^2 (\RR^{N+1})}  \\
& \le & c\left\{ \int_{|2\pi\xi|\ge1}
\frac{|{\cal F} f_{\text{odd}} (\xi)|^2}{|2\pi\xi|^4} \, d\xi 
\right\}^{1/2}
+ c\left\{ \int_{|2\pi\xi|\le1}
\frac{|{\cal F}f_{\text{odd}}(\xi)|^2}{|2\pi\xi|^4} \, d\xi
\right\}^{1/2} \\
& \le & c\left\{ \| f\|_{L^2 (\RR^{N+1})} +\left(
\int_{|2\pi\xi|\le1} |2\pi\xi|^{-4}\, d\xi 
\right)^{1/2}\sup_{\xi\in\RR^{N+1}} 
|{\cal F}f_{\text{odd}} (\xi)| 
\right\}\\
& \le & c \left\{ \| f\|_{L^2 (\RR^{N+1})} 
+\| f\|_{L^1 (\RR^{N+1})}
\right\} .
\end{eqnarray*}
Moreover,
\begin{eqnarray*}
\left \| \pd{\G f}{x_i} \right \|_{L^2 (\RR^{N+1})}  
& \le & \left\{ \int_{\RR^{N+1}}
\frac{|{\cal F}f_{\text{odd}}(\xi)|^2}{|2\pi\xi|^2} 
\, d\xi \right\}^{1/2} \\
& \le & c\left\{ \int_{|2\pi\xi|\le1} 
\frac{|{\cal F} f_{\text{odd}} (\xi)|^2}{|2\pi\xi|^2} 
\, d\xi \right\}^{1/2}
+ c\left\{ \int_{|2\pi\xi|\ge1} \frac{|{\cal F}f_{\text{odd}}
(\xi)|^2}{|2\pi\xi|^2} \, d\xi  \right\}^{1/2} \\
& \le & c \left\{ \| f\|_{L^2 (\RR^{N+1})} 
+\| f\|_{L^1 (\RR^{N+1})} \right\} .
\end{eqnarray*}
Now writing $\hat u (x_0,\xi') = \widehat{\G f} (x_0,\xi') 
+ B(x_0,\xi')$, i.e.\ setting
\begin{equation}\label{B}
B=\widehat{\G f_1}+\widehat{{\cal P}g_1}
+\widehat{{\cal P}g_2} \, , 
\end{equation} 
we 
must estimate the terms 
$\| B\|_{L^2 (\rnp)}$, 
$ \bigl \| \tpx B\bigr \|_{L^2 (\rnp)}$, and
$\|(\p/\p x_0) B\|_{L^2 (\rnp)}$.   We have
\begin{eqnarray*}
\| B\|_{L^2 (\rnp)}
& \le &
 c\left\{ \int_{|2\pi\xi|\le1}
\int_0^{+\infty} |B(x_0,\xi')|^2 \, dx_0 
\, d\xi' \right\}^{1/2}  \\
& & \qquad + c\left\{ \int_{|2\pi\xi|\ge1}
\int_0^{+\infty} |B(x_0,\xi')|^2 \, dx_0 
\, d\xi' \right\}^{1/2} .
\end{eqnarray*}
We begin by estimating the integral on 
the set where $\tpx\le1$.  Notice that,
if $N\ge4$, then
\begin{align*}
|B(x_0 ,\xi')|
& \le C \biggl\{
e^{-\tpx x_0} |\hat h (\xi')-\hat f(0,\xi')| 
\bigg| \frac{1}{1+2\tpx^2 +\tpx\Sx} -m_1 (\xi') \bigg| \\
& \qquad \qquad
+e^{-\tpx x_0} \tpx^{-1} \left| \int_0^{+\infty} e^{-\tpx y_0}
\hat f(y_0 ,\xi') dy_0 \right| \biggr\} \\
& \le C \biggl\{ e^{-\tpx x_0} \bigl(
\hat h (\xi')+\hat f(y_0 ,\xi') \bigr)\tpx 
+ \frac{e^{-\tpx x_0}}{\tpx}
\int_0^{+\infty}  |\hat f (y_0 ,\xi')|\, dy_0 \biggr\} . 
\end{align*}
From this,
the restriction theorem, and the fact that 
\begin{equation}\label{RL-est}
\int_0^{+\infty} | \hat f (y_0 ,\xi')|\, dy_0 
 \le \int_0^{+\infty} \int_{\RR^N} |f(y_0,x')|\, dx' \, dy_0 
= \| f\|_{L^1 (\rnp)} ,
\end{equation}
we obtain that
\begin{align*}
\lefteqn{
\biggl\{ \int_{\tpx\le1} \int_0^{+\infty} 
|B(x_0,\xi')|^2 \, dx_0 \, d\xi' \biggr\}^{1/2} }\\
& \le C \biggl\{ \| h\|_{L^2 (\RR^N)} 
		+ \| f(0,\xi')\|_{L^2 (\RR^N)}  
+\biggl( \int_{\tpx\le1} \frac{1}{\tpx^3} 
\bigl( \int_0^{+\infty} |\hat f(y_0,\xi')| \, dy_0 \bigr)^2 
\, d\xi'  \biggr)^{1/2} \biggr\} \\
& \le C \bigl\{  \| h\|_{L^2 (\RR^N)} + \| f\|_{W^s (\rnp)}
+ \| f\|_{L^1 (\rnp)} \bigr\}.
\end{align*}
To study the integral for $\tpx\ge 1$ we first notice that
\begin{align*}
\int_0^{+\infty} \left( e^{-\tpx x_0} - e^{-\Sx x_0}
\right)^2 \, dx_0 
& = \left( 2\tpx(\tpx +\Sx)^3 \right)^{-1} \\
& = O (\tpx^{-4}) \qquad \text{as}\  |\xi'|\rightarrow\infty. 
\end{align*}
Using this fact, Schwarz's inequality,
and the restriction theorem,
we see that the integral we want to estimate
is less than or equal to a constant times
\begin{eqnarray*}
\lefteqn{
\left\{ \int_{\tpx \ge1} \frac{1}{\tpx^4} \left[ \bigg| \tpx
\int_0^{+\infty}
e^{-\tpx y_0} \hat f(y_0,\xi') \, dy_0 \bigg|^2 
+ |\hat h(\xi')|^2 
+ |\hat f(0,\xi')|^2 
\right] \, d\xi' \right\}^{1/2} }\\
& \le & c\left\{ \bigl( \int_{\tpx\ge 1} \frac{1}{\tpx^3}
\int_0^{+\infty} |\hat f(y_0,\xi')|^2 \, dy_0
\, d\xi'\bigr)^{1/2}  
+ \| h\|_{L^2 (\RR^N)} +\| f(0,\cdot)\|_{L^2 (\RR^N)}
	\right\} \\
& \le & c\, \bigl\{ \| f\|_{W^s (\rnp)} + \| h\|_{L^2 (\RR^N)}\bigr\} . 
\end{eqnarray*}

\vfill 
\eject

In the same way we obtain that   
\begin{eqnarray*} 
\lefteqn{
\Bigl \| \tpx B \Bigr \|_{L^2 (\rnp)} }\\
& \le & \left\{ \int_{\tpx\le1} \int_0^{+\infty}
\tpx^2 |B(y_0,\xi')|^2 \, dy_0 \, d\xi'
\right\}^{1/2}\\
& & \quad  + \left\{ \int_{\tpx\ge1} \int_0^{+\infty}
\tpx^2 |B(y_0,\xi')|^2 \, dy_0 \,
d\xi' \right\}^{1/2} \\
& \le & c\, \Biggl\{ \| f(0,\cdot)\|_{L^2 (\RR^N )} 
+ \| h\|_{L^2 (\RR^N)} 
+\biggl(
\int_{\tpx\le1} \tpx^{-1}\, d\xi' \biggr)^{1/2} 
\| f\|_{L^1 (\RR^N)} \\
& & \quad + \biggl( \int_{\tpx\ge1} \bigg| \int_0^{+\infty}
e^{-\tpx y_0} \hat f(y_0,\xi')\, dy_0 \bigg|^2 \, d\xi'
\biggr)^{1/2} \\
& & \quad + \biggl( \int_{\tpx\ge1} \frac{|\hat h (\xi')|^2 + 
|\hat f(0,\xi')|^2}{\tpx^2} \, d\xi' \biggr)^{1/2} \Biggr\} \\
& \le & c\, \bigl\{ \| f\|_{W^s (\rnp)}+\| f\|_{L^1 (\rnp)} 
+\| h\|_{L^2 (\RR^N)} \bigr\} \, ,
\end{eqnarray*}
by the Schwarz inequality and the restriction theorem.
 
In order to estimate $\| \p B/\p x_0 \|_{L^2 (\rnp)}$
we observe that
\begin{multline*}
2 \int_0^{+\infty} \bigl( \tpx e^{-\tpx x_0} 
-\Sx e^{-\Sx x_0}\bigr)^2 \, dx_0 \\
= (\tpx +\Sx)^{-3} = O (\tpx^{-3})
\end{multline*}
as $|\xi'|\rightarrow \infty $.
Using this fact and (\ref{RL-est}) we see that
\begin{eqnarray*} 
\left\| \pd{B}{x_0} \right\|_{L^2 (\rnp)}
& \le & \left\{ \int_{\tpx\le1} \int_0^{+\infty} dy_0 \, d\xi'
\right\}^{1/2}  +
 \left\{ \int_{\tpx\ge1} \int_0^{+\infty} dy_0 \, d\xi'
\right\}^{1/2} \\
& \le & c\, \biggl\{
 \left( \int_{\RR^N} |\hat f(0,\xi')|^2 \, d\xi'
\right)^{1/2} + 
\left( \int_{\RR^N} |\hat h(\xi')|^2 \, d\xi'\right)^{1/2} \\
& & \quad 
+ \left( \int_{\tpx\le1} \tpx^{-1} \, d\xi'\right)^{1/2} 
\| f\|_{L^1 (\rnp)}\\ 
& &  + \quad\left( \int_{\tpx\ge1} \tpx^{-1} \bigg|
\int_0^{+\infty} e^{-\tpx y_0} \hat f(y_0,\xi') \, dy_0 \bigg|^2 \,
d\xi' \right)^{1/2}  \biggr\} \\
& \le & c\bigl\{ \| f\|_{W^{1/2} (\rnp)} +\| f\|_{L^1 (\rnp)}
+ \| h\|_{L^2 (\RR^N)} \bigr\} .
\end{eqnarray*}
This concludes the estimate in the case $N\ge4$.

{\bf Case $N<4$.}  
Now we must consider the lower dimensional cases.  Again, we 
first estimate $\G f$: 
\begin{eqnarray*}
\| \G f \|_{L^2 (\rnp)}
& = &  \left\{ \int_{|2\pi\xi|\le1} 
\frac{1}{\tpx^4}|{\cal F}f_{\text{odd}}
(\xi)|^2 \, d\xi
+ \int_{|2\pi\xi| > 1} \frac{1}{\tpx^4}
|{\cal F}f_{\text{odd}} (\xi)|^2 \, d\xi
\right\}^{1/2} \\
& \le & c\, \left\{ \left(
\int_{|2\pi\xi|\le1} \tpx^{-4}
|{\cal F}f_{\text{odd}} (\xi)|^2 \, d\xi
\right)^{1/2} + \| f\|_{L^2 (\RR^{N+1})} \right\} .
\end{eqnarray*}
If $N<4$ then  we need some conditions on $f$.  
Assume that $\int_{\rnp} f(x)\, dx =0$. Then 
${\cal F}f_{\text{odd}} (0) =0$, and if $N=2,3$ then 
\begin{eqnarray*}
\lefteqn{
\left( \int_{|2\pi\xi|\le1} \frac{1}{|2\pi\xi|^{4}} 
|{\cal F}f_{\text{odd}}(\xi)|^2
\, d\xi \right)^{1/2}} \\
& = & \left( \int_{|2\pi\xi|\le1} \frac{1}{|2\pi\xi|^2}
\frac{|{\cal F}f_{\text{odd}}(\xi) 
	-{\cal F}f_{\text{odd}} (0)|^2}{|2\pi\xi|^2}
\, d\xi \right)^{1/2} \\
& \le & \left( \int_{|2\pi\xi|\le1} |2\pi\xi|^{-2} \, d\xi
\right)^{1/2}
\sup_{\xi\in\RR^{N+1}} |\grad ({\cal F}f_{\text{odd}}) (\xi)| \\ 
& \le & c\int_{\rnp} |x|\cdot|f(x)|\, dx .
\end{eqnarray*}
 
If $N=1$ we need even stronger conditions
(because the logarithmic potential is poorly behaved
at infinity),
that is ${\cal F}f_{\text{odd}}(0)=0$, and 
$\grad ({\cal F}f_{\text{odd}})(0)=0$,
corresponding to
$$
\int_{\RR^2_{+}} f(x)\, dx =0 ;\quad
\int_{\RR^2_{+}} x_i f(x)\, dx =0,\ i=1,2; \quad
\int_{\RR^2_{+}} |x|^2 f(x)\, dx <\infty .
$$
Then, for $N=1$,
$$
\| \G f \|_{\RR^2_{+}}
\le c\left\{  \| f\|_{L^2 (\RR^2_{+})}  +\int_{\RR^2_{+}} |x|^2
f(x)\, dx \right\}.
$$
The estimate for $B$ has to be modified only in the part relative to 
the set where $\tpx\le1$:
\begin{eqnarray*}
\| B\|_{L^2 (\rnp)}
& \le & \left( \int_{\tpx\le1} \int_0^{+\infty} |B(y_0,\xi')|^2
\, dy_0 \, d\xi' \right)^{1/2} + \| f\|_{W^s (\rnp)} 
+ \| h\|_{L^2 (\RR^N)} \\
& \le & c \, \biggl\{  \| f\|_{W^s (\rnp)} +\| h\|_{L^2 (\RR^N)} \\
& & + \int_{\tpx\le1}
\tpx^{-3} \bigg| \int_0^{+\infty} \hat f(y_0,\xi') e^{-\tpx y_0} \,
dy_0 \bigg|^2 \biggr\}^{1/2} .
\end{eqnarray*}
With the above assumption we can estimate the last integral  with
\begin{multline*}
\biggl\{ \int_{\tpx\le1} \tpx^{-1}\bigg| \int_0^{+\infty} y_0 \hat f
(y_0,\xi') \frac{e^{-\tpx y_0} -1}{\tpx y_0} \, dy_0 \\   
+ \int_0^{+\infty} \frac{\hat f (y_0, \xi') 
-\hat f (y_0,0)}{\tpx} \,
dy_0 \bigg|^2 \, d\xi' \biggr \}^{1/2} 
 \le  c \int_{\rnp} |x| f(x)\, dx ,
\end{multline*}
if $N=2,3$.  If $N=1$ the estimate follows in a similar fashion.
\endpf 
 
 \section{Analysis of the Problem on the Half Space for $q$-Forms}
\label{SEC-Q-FMS}
\setcounter{equation}{0}
 
In this section we consider the space of $q$-forms, $q\ge 1$, with
coefficients in $W^1 (\rnp)$.
Throughout this section, we denote
this space of forms by $\bw^q$, that is, we do not write explicitly
the 
index for the Sobolev space $W^1(\rnp)$, this being fixed once and
for all.
 
On the space $\bw^q$ we select the basis $\{dx^I \}$, where
$I=(i_1,\dots i_q)$ is a strictly increasing multi-index.  For a
$q$-form $\phi$,
$$
\phi=\sum_I \phi_I dx^I
$$
we have that
\begin{eqnarray*}
d\phi
& = & \sum_I D_j \phi_I dx^j \wedge dx^I \\
& = & \sum_K \left( \sum_{I\, ,\, j}D_j \phi_I\e{K}{jI}
\right) dx^K ,
\end{eqnarray*}
where $K$ has $q+1$ entries, and
$$
\e{K}{jI}=\left\{
        \begin{array}{l c r}
        0 & \text{ if } & K\neq jI \text{ as sets}\\
        \pm 1 &  \text{ if } & K= jI \text{ as sets}
        \end{array}
\right.  
$$
and the sign is chosen according to the sign of the permutation that
puts $jI$ in increasing order.
\begin{lemma}\label{d'}     \sl
The formal adjoint $d'$ of $d$, $d':\bw^{q+1}\longrightarrow \bw^q$
with respect to the inner product in the Sobolev space, has the
following expression:
$$
d{}' \psi=\sum_{|I|=q} \left( \sum_{|K|=q+1 \atop j=0,\dots ,N}
\e{K}{jI} D_j \psi_K\right) dx^I .
$$
\end{lemma}
\bgpf Let $\phi\in\bw^q$, $\psi\in\bw^{q+1}$, $\phi=\sum_I
\phi_I dx^I$, $\psi=\sum_K \psi_k dx^K$, 
both with compact support in
$\rnp$. Then we have
\begin{eqnarray*}
\l d\phi,\psi,\r_1
& = & \sum_{|I|=q,|K|=q+1 \atop j=0,\dots,N}
\bigl\l D_j \phi_I dx^j \wedge dx^I, \psi_K dx^K \bigr\r_1 \\
& = & \sum_{I,K,j}\e{K}{jI} \bigl\l D_j \phi_I ,\psi_K
\bigr\r_1 \\ 
& = & \sum_{I,K,j}\biggl[ 
\sum_{k=0}^{N}  \e{K}{jI} \bigl\l D_k
D_j \phi_I,D_k \psi_K\bigr\r_0 
+ \e{K}{jI} \l D_j \phi_I ,\psi_K \r_0 \biggr] \\
& = & \sum_{I,K,j}\biggl[ \sum_{k=0}^{N}
- \e{K}{jI} \bigl\l D_k\phi_I,
D_k D_j \psi_K \bigr\r_0 
+ \e{K}{jI} \l D_j \phi_I ,\psi_K \r_0 \biggr] \\
& = & \sum_I \biggl[ \sum_{k=0}^{N} 
\bigl\l D_k \phi_I, D_k \bigl(
-\sum_{|K|=q+1 \atop j=0,\dots,N} \e{K}{jI}D_j \psi_K
\bigr)\bigr\r_0 \\
& & \qquad + \l \phi_I , \bigl( -\sum_{|K|=q+1\atop j=0,\dots,N} 
\e{K}{jI} D_j \psi_K \bigr) \r_0 \biggr]
. \qed 
\end{eqnarray*}
\renewcommand{\qed}{}\endpf
 
\begin{lemma}\label{dom-d*}   \sl
Let
$d^*$ be the adjoint of $d$ in the $W^1 (\rnp)$ norm.  Then 
$$
\dom d^* \cap \bw^{q+1}(\overline{\Omega}) 
=\left\{ \psi\in\bw^{q+1}(\overline{\Omega}) :
\left. D_0 \psi_{0J}\right|_{x_0 =0},\quad |J|=q \right\}.
$$
That is, if $\psi=\sum_{|K|=q+1}\psi_K dx^K$, 
then $\psi\in\dom d^*$ amounts to a condition
on the coefficients $\psi_K$ with index
$K=(0,\dots,k_N)\equiv 0J$ only, and it requires that
$$
\left. D_0 \psi_{0J}\right|_{x_0 =0} = 0 \, .
$$  
\end{lemma}

\vfill
\eject
\bgpf Let $\phi\in\bw^q$, $\psi\in\bw^{q+1}$.  Then
\begin{eqnarray*} 
\l d\phi,\psi\r_1
& = & \sum_{|K|=q+1} \biggl\l \sum_{|I|=q \atop j=0,\dots,N}
D_j \phi_I ,\psi_K \biggr\r_1 \\
& = & \sum_K \left[ \sum_{a=0,\dots,N} \biggl\l D_a
\left( \sum_{I,j} \e{K}{I,j} D_j \phi_I\right),
D_a \psi_K \biggr\r_0  
+\biggl\l \sum_{I,j} \e{K}{I,j} D_j \phi_I,\psi_K
\biggr\r_0  \right] \\
& = & \sum_K \left[ \sum_{a=0,\dots,N} \sum_{I,j}
\left( \e{K}{I,j} \biggl\l D_a \phi_I,
D_a D_j \psi_K \biggr\r_0 -\e{K}{0I}\int_{\RR^N}
\left. D_a\phi_I \overline{D_a \psi_K}  \right|_{x_0 =0}
\right) \right. \\
& & \quad \left. +\sum_{I,j} \left( -\e{K}{jI}
\biggl\l \phi_I,D_j \psi_K \biggr\r_0 -\e{K}{0I}\int_{\RR^N}
\left. \phi_I \overline{\psi_K}\right|_{x_0 =0} \right)
\right] \\
& = & \sum_{I} \left[ \sum_{a=0,\dots,N} \biggl \l D_a \phi_I,
 D_a \left( -\sum_{K,j}\e{K}{jI}D_j \psi_K\right)
\biggr \r_0 
+ \biggl\l \phi_I ,-\sum_{K,j}\e{K}{jI}D_j \psi_K
\biggr\r_0 \right] \\
& & \quad -\sum_{K,I,j} \biggl[ \sum_{a=0,\dots,N}\e{K}{0I}
\int_{\RR^N}  D_a \phi_I \overline{D_a \psi_K}
\biggl|_{x_0 =0}  
+ \e{K}{0I} \int_{\RR^N} \phi_I
\overline{\psi_K} \biggl |_{x_0 =0} \biggr] .
\end{eqnarray*}
Following the usual pattern we see that 
the only terms that we cannot
control with the $W^1$ norm of $\psi$ are
$$
\left.
\e{K}{0I} \int_{\RR^N} D_0 \phi_I \overline{D_0 \psi_K}
        \right|_{x_0 =0} \qquad \qquad
|K|=q+1,\quad |I|=q .
$$
Since $\phi$ was arbitrary, it must be $\left. D_0 \psi_K
\right|_{x_0 =0} =0$ when $K=0I$.
\endpf
 
Next we want to determine the expression for
the Hilbert space adjoint $d^*$ when acting on
$[\dom d^* \cap W^1_{q+1}]$.  We have the following result.
\begin{proposition}\label{espr-d*-q}    \sl
The adjoint of $d$ in the $W^1$-inner 
product is given by $d^* =d'
+{\cal K}$ where, for $\psi=\sum_{|K|=q+1}\psi_K dx^K$,
\begin{eqnarray*}
{\cal K} \psi
& = & \sum_{|I|=q} ({\cal K}\psi)_I \, dx^I \\
& = & \sum_{|I|=q \atop I\not\ni 0}\biggl ( K_{x_0}*\psi_{0I}
(0, \, \cdot \, )\biggr ) (x') \, dx^I .
\end{eqnarray*}
\end{proposition}
Before proving the proposition we make a few remarks.
\begin{remark}    \rm
\noindent \hbox{ \ }
\begin{enumerate}
\item Notice that ${\cal K}$ acts diagonally, in the sense that
$({\cal K}\psi)_I$ depends only on $\psi_{0I}$.
\item Only the terms of the form $\psi_{0I}$, 
i.e.\ terms for which
$\psi \lfloor \p/\p x_0 \neq 0$, contribute to ${\cal K}\psi$.
\item ${\cal K}\psi$ 
has only purely ``tangential" components, i.e.\
$({\cal K}\psi)_I =0$ if $0\in I$.
\end{enumerate}
\end{remark}
\begin{pf*}{\bf Proof of \ref{espr-d*-q}}
Notice that we could continue the calculation in the
previous proof to obtain (for $\psi\in\dom d^*$)
\begin{equation}\label{dag}
\l d\phi,\psi\r_1
= \l \phi, d'\psi\r_1 +\sum_{|I|=q \atop I\not\ni 0} \bigg[
\int_{\RR^N} \left. \phi_I \overline{\btu'\psi_{0I}} 
\right|_{x_0 =0}
-\int_{\RR^N} \left. \phi_I \overline{\psi_{0I}} \right|_{x_0 =0}
\bigg] .
\end{equation}
Now, as in the computation in Proposition \ref{expr-d*},
we write 
$$
d^* \psi = d' \psi +\theta,
$$
where $\theta={\cal K}\phi$, and ${\cal K}$ is 
a singular
Green's operator mapping $(q+1)$-forms into $q$-forms.
Now
\begin{equation}\label{dagg}
\l \phi,d'\psi\r_1 +\l \phi, \theta\r_1 = \l \phi,d^* \psi\r_1
=\l d\phi,\psi\r_1 .
\end{equation}
Therefore, (\ref{dag}), (\ref{dagg}) and Green's theorem give that
\begin{eqnarray*}
\lefteqn{
\sum_{|I|=q \atop I\not\ni 0} \bigg[
\int_{\RR^N} \left. \phi_I \overline{\btu'\psi_{0I}} 
\right|_{x_0 =0}
-\int_{\RR^N} \left. \phi_I \overline{\psi_{0I}} \right|_{x_0 =0}
\bigg] }\\
& = & \l \phi,\theta\r_1 \\
& = & \sum_I \biggl[ \sum_{a=0}^{N}
	\l D_a \phi_I ,D_a \theta_I \r_0
        +\l \phi_I ,\theta_I\r_0 \biggr] \\
& = & -\sum_I \l \phi_I ,\btu \theta_I \r_0
-\int_{\RR^N} \left. \phi_I \overline{D_0 \theta_I}\right|_{x_0 =0}
+\l \phi_I ,\theta_I\r_0 .
\end{eqnarray*}
Therefore $\theta$ must satisfy the following conditions
$$
\begin{cases}
\displaystyle{
-\btu \theta_I +\theta_I }  & = 0 \\
\left. -D_0 \theta_I \right|_{x_0 =0}
        & =  \left\{
                \begin{array}{l l}
               0 & \text{ if } 0\in I\\
              \btu' \psi_{0I} -\psi_{0I} & \text{ if }
                                      0\not\in I 
        \end{array} \right.
\end{cases}
$$
Thus
$$
\bigl( {\cal K} \psi \bigr)_I = \theta_I =0 
\qquad \text{if } 0\in  I,
$$
since the solution of $\btu u +u =0$ with the 
boundary condition 
$D_0 u(0,x')=0$ consists of just the zero function.  
On the other hand, if
$I\not\ni 0$ we have that $({\cal K}\psi)_I =\theta_I$ 
is the solution of 
\begin{eqnarray*}
\begin{cases}
\displaystyle{
-\btu u + u   =  0}  & \on  \rnp \\
\displaystyle{
\pd{u}{x_0}(0,\cdot )  =  -\btu' \psi_{0I} +\psi_{0I} }& 
						\on \RR^N  
\end{cases} .
\end{eqnarray*}
In the proof of Proposition \ref{expr-d*} , 
we showed that $\theta$
is given by
\begin{eqnarray*}
\theta_I
& = & - \int_{\RR^N} \sqrt{1 + |2\pi\xi'|^2}
      e^{-\sqrt{1 + |2\pi\xi'|^2} x_0}
  \widehat{\psi_{0I}}(0,\xi') e^{2\pi i \xi' \cdot x'}
            \, dx' \\
& \equiv & \biggl ( K_{x_0}*\phi_0(0, \, \cdot \, )\biggr ) (x'). 
\qed
\end{eqnarray*}
\renewcommand{\qed}{}\endpf
 
We now return to the main object of our work.  
We would like to solve
the \bvp
\begin{eqnarray*}
\begin{cases}
\displaystyle{
(dd^* +d^* d)\phi  =} & \alpha  \\
\phi\in\dom d^* &    \\
d\phi\in\dom d^* & 
\end{cases} 
\end{eqnarray*}
for $\alpha\in\bw^q.$
Recall that $d^* =d' +{\cal K}$.  A simple calculation shows that
\begin{equation}\label{?}
(dd' +d' d)\phi = \sum_{|I|=q} -\btu \phi_I dx^I .
\end{equation}
Thus we need to compute $d{\cal K}\phi$ and ${\cal K}d\phi$.  The
next lemma addresses this task.
\begin{lemma}\label{K-phi}   \sl
Let $\phi, d\phi \in \dom d^*$, with 
$\phi=\sum_I \phi_I dx^I$. Then
we have
$$
d{\cal K}\phi = \sum_{|J|=q} \left[
\sum_{|I|=q-1, I\not\ni 0\atop j=0,\dots,N}\e{J}{jI} D_j
K_{x_0}*\bigl( \phi_{0I}(0,\cdot) \bigr) \right] dx^J ,
$$
and
$$
{\cal K}d\phi = \sum_{|I|=q, I\not\ni 0} \left[
\sum_{|J|=q \atop j=0,\dots,N}
\e{0I}{jJ} K_{x_0}*\left( D_j \phi_J (0,\cdot)\right) \right] dx^I .
$$
\end{lemma}
\bgpf These are just straightforward computations.  We have
\begin{eqnarray*}
d{\cal K}\phi
& = & d \left( \sum_{|I|=q-1, I\not\ni 0} K_{x_0}*\bigl(
\phi_{0I}(0,\cdot)\bigr) dx^I \right) \\
& = & \sum_{|I|=q-1, I\not\ni 0} \sum_{j=0,\dots,N} D_j
K_{x_0}*\bigl( \phi_{0I}(0,\cdot)\bigr) dx^j \wedge dx^I  \\
& = & \sum_{|J|=q} \left[
\sum_{|I|=q-1, I\not\ni 0 \atop j=0,\dots,N} \e{J}{jI} D_j
K_{x_0}*\bigl( \phi_{0I}(0,\cdot)\bigr) \right] dx^J  .
\end{eqnarray*}
This proves the statement for $d{\cal K}\phi$.  On the other hand,
\begin{eqnarray*}
{\cal K}d\phi
& = & {\cal K}\left( \sum_{|K|=q+1}\biggl(
\sum_{|J|=q\atop j=0,\dots,N} \e{K}{jJ} D_j \phi_J \biggr) dx^K
\right) \\
& = & \sum_{|I|=q, I\not\ni 0}K_{x_0}*\left(
\sum_{|J|=q\atop j=0,\dots,N} \e{0I}{jJ} D_j \phi_J (0,\cdot)
\right) dx^I \\
& = & \sum_{|I|=q, I\not\ni 0} \left[
\sum_{|J|=q\atop j=0,\dots,N}
\e{0I}{jJ} K_{x_0}* \biggl( D_j \phi_J (0,\cdot) \biggr)\right]dx^I .
\qed
\end{eqnarray*}
\renewcommand{\qed}{}\endpf
 
Next we want to compute $({\cal K}d+d{\cal K})\phi$.
[Note that, in this discussion, the letter $K$ is both a kernel and
an index; no confusion should result.]
\begin{proposition} \label{6.6}  \sl
Let $\phi\in\bw^{q}$, with $\phi,d\phi\in\dom d^*$. Set
$({\cal K}d+d{\cal K})\phi\equiv \beta=\sum_{|K|=q}\beta_K dx^K$.
Then
\begin{eqnarray*}
\beta_K & = & \pd{}{x_0} K_{x_0}*\biggl( \phi_K (0,\cdot)\biggr)
                \qquad\qquad \text{if } K\ni 0\\
\beta_K & = & K_{x_0}*\biggl( \pd{}{x_0} \phi_K (0,\cdot)\biggr)
                \qquad\qquad \text{if } K\not\ni 0 .
\end{eqnarray*}
\end{proposition}
\noindent {\bf Remark.}  
Notice that ${\cal K}d +d{\cal K}$ is a diagonal operator on the
space of $q$-forms $\bw^q$.
\begin{pf*}{\bf Proof of \ref{6.6}} 
Suppose that $K\ni 0$. 
Then $({\cal K}d\phi)_{K}=0$,
since ${\cal K}$ applied to any form has only tangential components.
Moreover, for the same reason, $(d{\cal K}\phi)_{0K'}$ can be
obtained only by differentiating $({\cal K}\phi)_{K'}$ by $D_0$, and
``wedging" by $dx^0$.  That is,
\begin{eqnarray*}
(d{\cal K}\phi)_{0K'}
& = & D_0 ({\cal K}\phi)_{K'} dx^0 \wedge dx^{K'} \\
& = & D_0 \biggl( K_{x_0}*\bigl( \phi_{0K'}(0,\cdot)\bigr) dx^0
        \wedge d^{K'} .
\end{eqnarray*}
Now suppose that $K\not\ni 0$.  We use Lemma \ref{K-phi} to
obtain that
\begin{eqnarray*}
d{\cal K}\phi +{\cal K}d \phi
& = & \sum_{|K|=q} \left[
\sum_{|I|=q-1, I\not\ni 0\atop j=0,\dots,N} \e{K}{jI} D_j
K_{x_0}*\bigl( \phi_{0I}(0,\cdot)\bigr) \right] dx^K \\
& & \qquad\qquad + \sum_{|K|=q} \left[
\sum_{|J|=q \atop j=0,\dots,N} \e{0K}{jJ}
K_{x_0}*\bigl( D_j \phi_{0J}(0,\cdot)\bigr) \right] dx^K .
\end{eqnarray*}
Thus we need to describe the term on the right hand side of this
equation when $K\not\ni 0$.  Notice that in this case $j$
cannot be $0$ in the first sum, since otherwise $\e{K}{0I}=0$.  Then,
when $K\not\ni 0$ the coefficient of $dx^K$ in
the right hand side above equals
\begin{multline*}
K_{x_0}*\bigl( D_0 \phi_{0K}(0,\cdot)\bigr)
+ \sum_{|J|=q,J\not\ni0 \atop j=1,\dots,N} \e{jK}{0J}
D_j \biggl( K_{x_0}*\bigl( \phi_{0J}(0,\cdot)\bigr)\biggr) \\
+ \sum_{|I|=q-1, I\not\ni 0\atop j=1,\dots,N} \e{K}{jI} D_j
\biggl( K_{x_0}*\bigl( \phi_{0I}(0,\cdot)\bigr) \biggr) 
\equiv  K_{x_0}*\bigl( D_0 \phi_K (0,\cdot) \bigr) ,
\end{multline*}
since $J$ must be of the form $0I$ with $|I|=q-1,\, I\not\ni0$, and
therefore 
$$
\e{0K}{j0I}= -\e{0K}{j0I}=-\e{K}{jI} .
$$
This concludes the proof.
\endpf
\begin{corollary}  \sl
Let $({\cal K}d+d{\cal K})\psi=G\psi$. Then $\beta=G\psi$ 
is a form such that its components solve the 
\bvp s 
$$
\begin{cases}
\displaystyle{
(-\btu+I)\beta_K} & = 0\\
\displaystyle{ \pd{\beta_K}{x_0} }& = (-\btu' +I)(d\psi)_{0K} 
\end{cases} 
\qquad \text{if } K\not\ni 0 ,
$$
and
$$
\begin{cases}
(-\btu+I)\beta_K & = 0\\
\beta_K & = (-\btu' +I)\psi_{K} 
\end{cases} 
\qquad \text{if } K\ni 0 .
$$
\end{corollary}
\bgpf
By Proposition \ref{espr-d*-q}, if 
$G\psi =\sum_{|K|=q} \beta_K dx^K$ 
for $K\not\ni 0$, then 
$$
\beta_K = K_{x_0}*\left( \pd{\psi_K}{x_0} (0,\cdot)\right) .
$$
By construction
$$
(-\btu+I)\beta_K =0 \qquad \text{on }\rnp,
$$
and
\begin{align*}
\pd{\beta_K}{x_0} 
& = (-\btu' +I) (d\psi)_{0K} \\
& = (-\btu' +I)\left( \pd{\psi_K}{x_0} +T_1 \psi_{K'}\right)
\qquad \text{on } b\rnp ,\ K'\ni 0,
\end{align*}
where $T_1$ is a first order tangential differential operator. 

On the other hand, if $K\ni 0$, then $({\cal K} d\psi)_K =0$. Thus, for
$K=0K'$, 
\begin{align*}
(G\psi)_K 
& = (d{\cal K}\psi)_K \\
& = \pd{}{x_0} ({\cal K}\psi)_{K'}.
\end{align*}
Then, on $\rnp$,
$$
(-\btu+I) (G\psi)_K 
= \pd{}{x_0} (-\btu+I)({\cal K}\psi)_{K'} =0 ;
$$
and if $K=0K'$,
\begin{align*}
\left. (G\psi)_K \right|_{b\Omega}
& = \left. \pd{}{x_0} ({\cal K}\psi)_{K'} \right|_{b\Omega} \\
& = \left. (-\btu' +I) \psi_K \right|_{b\Omega} . 
   \qquad \qquad \qquad \qquad \qed
\end{align*}
\renewcommand{\qed}{}\endpf 

\noindent Before analyzing the \bvp\ we need one more result;
that is, we wish to make explicit the boundary condition 
$d\phi\in\dom d^*$.
\begin{lemma}\label{2bdry-cond-q}  \sl
Let $\phi\in\bw^q$, and
let $\phi\in\dom d^*$. Then $d\phi\in\dom d^*$ if and only
if
$$
\pd{^2 \phi_K}{x_0^2} (0,\cdot) =0  \qquad \text{for }K\not\ni 0.
$$
\end{lemma}
\bgpf By Lemma \ref{dom-d*} we have that $d\phi \in\dom d^*$
if and only if, for all multi-indices $K'\not\ni0,\, |K'|=q$,
there holds the equality
$$
\pd{}{x_0}  \biggl( \sum_{I,j}
\e{0K'}{jI} \pd{\phi_I}{x_j} \biggr) (0,x')=0 .
$$
If $j\neq 0$, then $I\ni 0$ and
$$
\pd{}{x_0} \pd{\phi_I}{x_j}(0,\cdot)
= \pd{}{x_j} \left( \pd{\phi_I}{x_0}(0,\cdot) \right) =0 ,
$$
since, for $\phi\in \dom d^*$, we have $\p /\p x_0 \phi_i (0,\cdot)=0$ when
$I\ni0$. 
Thus we obtain
\begin{equation*}
\pd{^2 \phi_I}{x^2_0}(0,\cdot)
=0 \qquad\qquad\text{for all } I\not\ni0,\ |I|=q. \qed
\end{equation*}
\renewcommand{\qed}{}\endpf
 
We finally are able to formulate the \bvp.
\begin{theorem}\label{q-bvp}      \sl
Consider the \bvp\
\begin{equation}\label{bp-q}
\begin{cases}
\displaystyle{
(dd^* +d^* d)\phi  =  \alpha} \qquad 
		& \alpha\in\bw^q (\overline{\Omega}) \\
\phi\in\dom d^*  & \\
d\phi\in\dom d^* & 
\end{cases}
\end{equation}
Let $\alpha$ be a $q$-form with coefficients in
$W^r (\rnp) \cap L^1 (\rnp)$, with $r>1/2$.  
Then, if $N\ge 4$, there
is a unique $q$-form $\phi$ with coefficients in $W^{r+2} (\rnp)$
solving the \bvp\ and satisfying the estimate
$$
\| \phi \|_{W^{r+2}} \le c\cdot \left\{ \| \alpha\|_{W^r (\rnp)}
+ \| \alpha\|_{L^1 (\rnp)} \right\} .
$$
If $N=2,3$ we need to further require that
$$
\int \alpha (x)\, dx =0\ , \qquad \qquad \alpha\in L^1 (|x|dx,\rnp), 
$$
so that the solution $\phi$ satisfies the estimate
$$
\| \phi \|_{W^{r+2}} \le c\cdot \left\{ \| \alpha\|_{W^r (\rnp)}
+ \| \alpha\|_{L^1 (|x|dx,\rnp)} \right\} .
$$
If $N=1$ we need to require that
$$
\int \alpha (x)\, dx =0\ , \qquad\int x_i \alpha(x)\, dx =0 ,\ i=0,1
 \qquad \alpha\in L^1 (|x|^2 dx, \RR^2 ) ;
$$
in this case similar estimates hold for the solution $\phi$, 
with the addition of the term
$\| \alpha\|_{L^1 (|x|^2 dx, \RR^2)} $ to the right hand side.
\end{theorem}

\noindent By the discussion preceding 
the statement of the theorem we see that
the \bvp\ (\ref{bp-q}) is equivalent to the two scalar problems
\begin{equation}\label{0inK}
\begin{cases}
\displaystyle{
 -\btu \phi_K +\pd{}{x_0} \big( K_{x_0}*\phi_K |_{x_0 = 0}
\big)  =  \alpha_K  } \qquad & \\
\displaystyle{
\pd{\phi_K}{x_0} \big|_{x_0 =0}  =  0} & 
\end{cases} 
\qquad \qquad 0 \in K, 
\end{equation}
and 
\begin{equation}\label{0not-inK}
\begin{cases}
\displaystyle{
-\btu \phi_K + \big( K_{x_0}*\pd{\phi_K}{x_0} \big|_{x_0 = 0} 
\big)  =  \alpha_K }\qquad &  \\
\displaystyle{
\pd{^2 \phi_K}{x^2_0} \big|_{x_0 =0}  =  0}  &
\end{cases} 
 \qquad  0 \not\in K.  
\end{equation}
Thus the problem is reduced to solving two different \bvp s.  The
\bvp\ (\ref{0not-inK}) 
has been solved already in the case $q=0$.
The solution of the second \bvp\ is contained in the following
theorems.  As in the case of functions, for greater 
flexibility we
solve the \bvp\ with non-zero boundary data.
\begin{theorem}\label{apriori-est-q-hs}  \sl
Consider the \bvp\
\begin{equation}\label{non-0-bndry-data}
\begin{cases}
\displaystyle{
-\btu u +\pd{}{x_0} \big( K_{x_0}* u(0,\cdot)
\big)  =  f} \qquad & \on\rnp \\
\displaystyle{
\pd{u}{x_0} (0,\cdot)  =  h} & \on\RR^N 
\end{cases} 
\end{equation}
If there exists a solution $u$, then it 
satisfies the {\em a priori} estimate
$$
\| u \|_{W^{r+2}} \le c\cdot \left\{ \| f \|_{W^r (\rnp)}
+ \| h \|_{W^{r+ 1/2} (\RR^N)}
+ \| u\|_{L^2 (\rnp)} \right\} ,
$$
for $r> 1/2$.
\end{theorem}
\bgpf We break the proof into two steps.  We first estimate the
function $u$ which is the solution of the \bvp\ with zero boundary 
data.  Then we show how to reduce the general case to this 
special one. 

In this proof we let ${\cal J}$ denote the Green's
function for the Neumann problem for the 
Laplacian on the half space. 
Suppose now that $u$ is the solution of the \bvp\ with $h=0$. 
Then  $u$ can be written as
$$
u = {\cal J} \biggl( f -\pd{}{x_0} 
\bigl (K_{x_0}*u(0,\cdot) \bigr ) \biggr) .
$$
Recall that, for $N\ge 2$,
\begin{multline*}
({\cal J} f )(x_0, x')
=  \frac{1}{(N-1)\omega_{N+1}}
\int_{\rnp} \biggl\{ \frac{1}{(x_0 -y_0 )^2
        + |x' -y'|^2)^{(N-1)/2}} \\
+ \frac{1}{(x_0 +y_0 )^2 + |x' -y'|^2)^{(N-1)/2}} 
\biggr\}  f(y_0,y')\, dy_0 \, dy' ,
\end{multline*}
and
$$
\widehat{{\cal J}f} (x_0,\xi') = \frac{1}{2\tpx} \int_0^{+\infty}
\bigl( e^{-\tpx |x_0 - y_0|} + e^{-\tpx (x_0 +y_0 )} \bigr) \hat f
(y_0,\xi') \, dy_0 .
$$   
Then
\begin{eqnarray*}
\hat u(0,\xi')
&= & \widehat{{\cal J}f} (0,\xi') -\tpx^{-1}\int_0^{+\infty} 
	e^{-(\tpx+\Sx)y_0} \, dy_0 \, (1+\tpx^2)\hat u(0,\xi') \\
& = & \widehat{{\cal J}f} (0,\xi') -\frac{(1+\tpx^2)}{\tpx\bigl (\tpx+\Sx\bigr )}
\hat u(0,\xi') .
\end{eqnarray*} 
Thus, recalling the definition (\ref{m_2}) of $m_2$, we have that
$$
\hat u(0,\xi') = \tpx m_2 (\xi') \widehat{{\cal J}f}(\xi') .
$$
Set $F=f-(\p/\p x_0)\bigl(K_{x_0}*u(0,\cdot)\bigr)$.
Since $u={\cal J}F$, in order to estimate
$\|u\|_{W^{3/2} (\rnp)}$ we can proceed as in the proof of Theorem
\ref{apr-est} and, using the fact that $\p^2 /\p x_0^2 =\btu-\btu'$,
we can reduce to estimating derivatives of $u$ of order not exceeding
$2$ in the $x_0$-variable.

Let 
$g_{\text{e}}$ denote 
the even extension (in the $x_0$ variable)
to all of $\RR^{N+1}$ of the function $g$
defined initially on $\rnp$.  
We have that
$$
\biggl[ \bigl({\cal J}g\bigr)(x_0 ,x') \biggr]_{\text{e}} 
= \bigl( {\cal N} g_{\text{e}} \bigr) (x_0, x') \ ,
$$
and
$$
{\cal F} \left ( \pd{}{x_0} K_{x_0} \right )_{\text{e}} 
= \frac{2(1+\tpx^2)^{3/2}}{1+|2\pi\xi|^2} .
$$
Combining all of these facts we obtain:
\begin{equation}\label{luigi1}
{\cal F}u_{\text{e}} =
|2\pi\xi|^{-2} \biggl( {\cal F}\fe (\xi)
-\frac{2(1+\tpx^2)^{3/2}}{1+|2\pi\xi|^2} 
\tpx m_2 (\xi') \widehat{{\cal J}f} (0,\xi') \biggr) .
\end{equation}
Thus
\begin{align*}
\| u\|_{W^{2+r}(\rnp)} 
& \le c \biggl( \int_{\RR^{N+1}} \sum_{i,j=0}^{N} (1+|2\pi\xi_i|^2
|2\pi\xi_j|^2 )|{\cal F}u(\xi)|^2 (1+\tpx^2)^r \, d\xi \biggr) \\
& \le c \biggl( \| u\|_{L^2 (\rnp)} + \int_{-\infty}^{+\infty}  
\int_{\tpx\ge1} 
|{\cal F}F_{\text{e}}(\xi)|^2 \tpx^{2r}\, d\xi' \, d\xi_0 \biggr),
\end{align*}
where
\begin{align*}
{\cal F}F_{\text{e}} (\xi)
& = {\cal F}\fe (\xi) - \frac{2(1+\tpx^2)^{3/2}}{1+|2\pi\xi|^2} 
m_2 (\xi')\int_{0}^{+\infty}e^{-\tpx y_0}\hat f(y_0,\xi')\, dy_0 \\
& = {\cal F}\fe (\xi) - \frac{2(1+\tpx^2)}{1+|2\pi\xi|^2} 
m_1 (\xi')\int_{0}^{+\infty}e^{-\tpx y_0}\hat f(y_0,\xi')\, dy_0 .
\end{align*}
Therefore, by obvious calculations and the Schwarz inequality,
\begin{align*}
\| u\|_{W^{2+r}(\rnp)}^2 
& \le c \biggl( \|u \|_{L^2 (\rnp)}^2 + \| f\|_{W^r (\rnp)}^2 
+\int_{\tpx\ge1} \frac{(1+\tpx^2)^{2+r}}{|2\pi\xi|}  \\
& \qquad 
\times \bigl( 
\int_0^{+\infty} |\hat f(y_0,\xi')|^2 \, dy_0 \bigr)\bigl( 
\int_0^{+\infty} \frac{1}{(1+\tpx^2)^2 +|2\pi\xi_0 |^4}\, d\xi_0 
\bigr) \, d\xi' \biggr) \\
& \le c \bigl( \| u\|_{L^2 (\rnp)}^2 + \|f\|_{W^{r}(\rnp)}^2 \bigr) .
\end{align*}
This proves the theorem
in the case that the boundary data $h \equiv 0$.

In the general case, let ${\cal Q}$ be the operator initially defined
on $C^\infty_0 (\RR^N)$ by
$$
({\cal Q}g)\widehat{\mathstrut} (x_0,\xi')
= - \, \frac{e^{-\Sx x_0}}{\Sx} \, \hat g (\xi').
$$
Notice that $(\p/\p x_0)({\cal Q}g)|_{x_0 =0} =g$, and that 
$$
{\cal Q}: W^s (\RR^N)\longrightarrow W^{s+1/2} (\rnp).  
$$
We seek 
an {\em a priori} estimate for a function $u$ that solves of
(\ref{non-0-bndry-data}).  Let $G'$ denote the operator 
$u\longmapsto(\p/\p x_0)\big(K_{x_0}*u(0,\cdot)\big)$. 
Set $v=u-{\cal Q}h$. Then $v$ solves the 
\bvp\
$$
\begin{cases}
\displaystyle{
(-\btu+G' )v  =  f+(-\btu +G'){\cal Q}h }\qquad & \on\rnp \\
\displaystyle{
\pd{v}{x_0} (0,\cdot)  =  0} & \on\RR^N 
\end{cases} .
$$
Thus for such $v$ we have the usual {\em a priori} estimate. Notice that
$$
\bigl[ (-\btu+G'){\cal Q}h\bigr]\hat{\text{\,}} (\xi')
= -\frac{\tpx^2 e^{-\Sx x_0}}{\Sx}\hat h(\xi') ,
$$
so that the operator $(-\btu +G'){\cal Q}$ 
has the same behavior as 
the operator $\widetilde{\cal K}$ studied in Proposition
\ref{tilde-K}, so it maps $W^s (\RR^N)$ into $W^{s+1/2}(\rnp)$
continuously. 
Therefore
\begin{align*}
\| u\|_{W^{r+2}(\rnp)}
& \le   \| v\|_{W^{r+2}(\rnp)} 
	+ \| {\cal Q}h \|_{W^{r+2} (\rnp)} \\
& \le c\, \bigl\{ \| f\|_{W^r (\rnp)} 
	+ \|(-\btu+G'){\cal Q}h\|_{W^r (\rnp)}  
	+\| v\|_{L^2 (\rnp)} \\
&  \qquad \qquad 
	+ \| {\cal Q}h\|_{W^r (\rnp)} \bigr\} \\
& \le c\, \bigl\{ \| f\|_{W^r (\rnp)} +\| h\|_{W^{r+1/2} (\RR^N)}
	+\| u\|_{L^2 (\rnp)} \bigr\} ,
\end{align*}
since $v=u-{\cal Q}h$.
This concludes the proof.
\end{pf*} 
\begin{theorem}\label{exist-q}     \sl
Let  $N\ge 4$. Then for any $f\in L^2 (\rnp)\cap L^2(\rnp)$,  
there exists a unique function $u$
solving the \bvp\ (\ref{0inK}) and such that
$$
\| u \|_{L^2 (\rnp)} \le c\cdot \left\{ \| f\|_{L^2 (\rnp)}
+ \| f \|_{L^1 (\rnp)} +\|h \|_{W^{1/2}(\RR^N)} \right\} .
$$
If $N=2,3$ we suppose in addition that $f$ satisfies
$$
\int f(x)\,  dx =0 , \qquad \qquad f\in L^1 (|x|dx,\rnp) .
$$
Then there exists a unique solution $f$ that satisfies the estimate
$$
\| u \|_{L^2 (\rnp)} \le c\cdot \left\{ \| f\|_{L^2 (\rnp)}
+ \| f\|_{L^1 (|x|dx,\rnp)} +\| h\|_{W^{1/2}(\RR^N)} \right\} .
$$
If $N=1$ we suppose that $f$ is such that
$$
\int f(x)\, dx =0 , \qquad\int x_i f(x)\, dx =0 ,\ i=0,1
 \qquad f \in L^1 (|x|^2 dx, \RR^2 ).
$$
Then there exists a unique solution $u$ such that
$$
\| u \|_{L^2 (\rnp)} \le c\cdot \left\{ \| f\|_{L^2 (\rnp)}
+ \| f\|_{L^1 (|x|^2 dx,\rnp)} +\| h\|_{W^{1/2}(\RR^N)} \right\} . 
$$
\end{theorem}
\bgpf
Recall that we have found that if a  solution $u$ exists, 
then it must
be given by the formula (\ref{luigi1}), i.e.
$$
{\cal F} \ue (\xi)
 =  |2\pi \xi|^{-2}\biggl( {\cal F}\fe (\xi)
-2\frac{(1+\tpx^2)^{3/2}}{1+|2\pi \xi|^2} m_2 (\xi') 
\widehat F_1(\xi') \biggr) ,
$$
where $F_1$ is given by 
$$
\widehat F_1 (\xi')= \tpx \widehat{{\cal J}f} (0,\xi') 
= \int_{0}^{+\infty} e^{-\tpx y_0}
     \hat f(y_0 ,\xi') \, dy_0 . 
$$
In order to show that $u$ given by the above formula is indeed
a solution, by the estimate in Theorem \ref{apriori-est-q-hs} it
suffices to show that $u\in L^2 (\rnp)$. 

Suppose that $N\ge 4$.  We have
\begin{eqnarray*}
\lefteqn{\| u \|_{L^2 (\rnp)} }\\
& = & \frac12 \int_{\RR^{N+1}} |{\cal F} \ue (\xi)|^2 \, d\xi \\ 
& \le & \left\{ \int_{\tpx\le1} \int_{-\infty}^{+\infty} 
|{\cal F} \fe(\xi) |^2 \, d\xi_0 \, d\xi' \right\}^{1/2}
+ \left\{ \int_{\tpx\ge1} \int_{-\infty}^{+\infty} 
|{\cal F} \fe(\xi) |^2 \, d\xi_0 \, d\xi' \right\}^{1/2} \\
& \le & c \biggl\{ \bigl( 
\int_{\tpx\le1} \int_{-\infty}^{+\infty} |2\pi\xi|^{-4} 
\, d\xi_0 \, d\xi'  \bigr)
\sup_{\xi\in\RR^{N+1}} \biggl( |{\cal F} \fe (\xi)| +\int_0^{+\infty}
|\hat f (y_0,\xi')| \, dy_0 \biggr)^2 \\
& & \qquad +\int_{\tpx\ge1} \int_{-\infty}^{+\infty} 
|{\cal F}\fe(\xi)|^2 \, d\xi_0 \, d\xi' \\
& & \qquad 
+ \int_{\tpx\ge1} \int_{-\infty}^{+\infty} |2\pi\xi|^{-4}  \, d\xi_0 
\big| \int_0^{+\infty} e^{-\tpx y_0} \hat f(y_0,\xi') \, dy_0 \big|^2 
\, d\xi' \biggr\} \\
& \le & c\left\{ \|f\|_{L^1 (\rnp)} + \|f\|_{L^2 (\rnp)} \right\} .
\end{eqnarray*}
This proves our statement in the case $N\ge 4$.
 
Let now $N=2,3$,
and
assume that
$$
\int_{\rnp} f(x) \, dx =0 .
$$
Notice that this implies that 
$\int_0^{+\infty} \hat f (y_0,0)\, dy_0 =0$.
Therefore,
\begin{align*}
\| u \|_{L^2 (\rnp)}^2  
& \le  \int_{|2\pi\xi|\le1} \int_{-\infty}^{+\infty} 
\frac{|{\cal F}\fe (\xi)-{\cal F}\fe (0)|^2}{|2\pi\xi|^4}
\, d\xi_0 \, d\xi'\\
& \qquad + \int_{|2\pi\xi|\le1} \int_{-\infty}^{+\infty} 
|2\pi\xi|^{-2} \big|\int_0^{+\infty} e^{-\tpx y_0} 
\frac{\hat f(y_0,\xi') -\hat f(y_0,0)}{|2\pi\xi|}\, dy_0 \big|^2 
\, d\xi_0 \, d\xi'\\
& \qquad + \int_{|2\pi\xi|\le1} \int_{-\infty}^{+\infty} 
|2\pi\xi|^{-2} \biggl( \int_0^{+\infty} 
\frac{|e^{-\tpx y_0} -1|}{\tpx} |\hat f(y_0,0)|\, dy_0 \biggr)^2 
\, d\xi_0 \, d\xi'\\ 
& \qquad + \int_{\tpx\ge1}\int_{-\infty}^{+\infty} \biggl[ 
|{\cal F}\fe (\xi)|^2 + \bigl( \int_0^{+\infty} e^{-\tpx y_0} 
|\hat f(y_0,\xi')|\, dy_0 \bigr)^2 
\biggr] \, d\xi_0 \, d\xi'\\
& \le C \int_{\tpx\le1} \tpx^{-1} \, d\xi' \left\{ 
\bigl( \int_{\rnp} |x|\cdot|f(x)|\, dx \bigr)^2 
+ \bigl( \int_{\rnp} |f(x)|\, dx \bigr)^2 \right\} \\
& \qquad 
+C\| f\|_{L^2 (\rnp)} .   
\end{align*}
Here we have used the following facts:
$$
\frac{|{\cal F}\fe (\xi)-{\cal F}\fe (0)|}{|2\pi\xi|}
\le \frac{1}{2\pi} \sup_{\xi\in\RR^{N+1}} |\grad {\cal F}\fe (\xi)| 
\le \int_{\rnp} |x|\cdot |f(x)|\, dx ;
$$
$$
\int_0^{+\infty} e^{-\tpx y_0} 
|\hat f(y_0,\xi')|\, dy_0 
\le \int_0^{+\infty} y_0 \int_{\RR^N} |f(y_0,y')| \, dy'  
\, dy_0  \le \int_{\rnp} |y|\cdot |f(y)|\, dy ;  
$$
and
\begin{eqnarray*}
\lefteqn{
\int_{\tpx\ge1} \biggl(
\int_{0}^{+\infty} e^{-\tpx y_0} |\hat f(y_0,\xi')| \,
dy_0 \biggr)^2 \, d\xi'} \\
& \le &
\int_{\tpx\ge1} \frac{1}{2\tpx} \int_0^{+\infty}
|\hat f (y_0,\xi')|^2 \, dy_0 \, d\xi' \\
& \le & \int_{\rnp} |f(x)|^2 \, dx . 
\end{eqnarray*}
Finally, in the case $N=1$ we can obtain 
the same kind of estimate if
we require the data $f$ to satisfy the additional stated
compatibility conditions.  

This  proves the result in the case of the boundary data 
$h \equiv 0$.
If $h\neq0$ we proceed as in the proof of Theorem
\ref{apriori-est-q-hs}.  Consider 
$v=u-{\cal Q}h$.  Then $v$ satisfies the above estimate, and
therefore
$$
\| u\|_{L^2 (\rnp)} \le \| v\|_{L^2 (\rnp)} 
+ \|{\cal Q}h\|_{L^2 (\rnp)} 
\le c\bigl\{  \| v\|_{L^2 (\rnp)} + \| h\|_{L^2 (\rnp)} \bigr\} .
$$
In order to estimate $\| v\|_{L^2 (\rnp)}$ we observe that we have to
replace $f$ with $f+(\btu -G'){\cal Q}h$ in the computations above,
so that, for $N\ge4$ (the case $N<4$ is similar), 
\begin{align*}
\| v\|_{L^2 (\rnp)} 
& \le \frac12 \int_{\RR^{N+1}} |{\cal F}v_{\text{e}} |^2 \, d\xi \\
& \le C \bigl\{ \|f \|_{L^2 (\rnp)} + \|f\|_{L^1 (\rnp)} \bigr\} \\
& \quad +C\int_{\RR^{N+1}} \frac{1}{|2\pi\xi|^{4}}
 \bigg| {\cal F} \bigl(
(\btu -G'){\cal Q}h\bigr)_{\text{e}}
-\frac{2(1+\tpx^2)}{1+|2\pi\xi|^2} m_1 (\xi') \hat F_2 (\xi')
\bigg|^2 d\xi ,
\end{align*}
where
\begin{align*}
\hat F_2 (\xi') 
& = -\int_0^{+\infty} e^{-\tpx y_0}e^{-\Sx y_0} \frac{\tpx^2}{\Sx} \,
d\xi_0 \hat h(\xi') \\
& = -\frac{\tpx^2}{(\tpx+\Sx)\Sx} \hat h (\xi') ;
\end{align*}
and
$$
{\cal F}\bigl( (\btu-G'){\cal Q}h\bigr)_{\text{e}} 
= \frac{2\tpx^2}{1+|2\pi\xi|^2} \hat h (\xi) .
$$
Therefore we obtain
$$
\| v\|_{L^2 (\rnp)} \le C\bigl( 
\| f\|_{L^2 (\rnp)} + \| f\|_{L^2 (\rnp)} + \| h\|_{L^2 (\rnp)} .
$$
This concludes the proof.
\endpf 

\newpage 
\section*{\bf {THE PROBLEM ON A SMOOTHLY BOUNDED DOMAIN}}
\mbox{}

\section{Formulation of the 
Problem on a Smoothly Bounded Domain}
\setcounter{equation}{0}

Let $\Omega = \{ x \in \RR^{N+1}:
\rho(x) < 0\}.$  To avoid pathologies,
we assume that $\rho$ is a $C^\infty$
function on $\RR^{N+1}$ with the
property that $\nabla \rho \ne 0$
at all points of $b\Omega.$
Then $\Omega$ is a domain with $C^\infty$
boundary (see \cite{KR2} for a
protracted discussion of these matters).
We shall also assume that
$\Omega$ is bounded.

It will simplify our calculations if we assume in advance that
$(\partial \rho/\partial n)
= |\nabla \rho| = 1$
on $b\Omega.$  This is easily arranged.
 
Recall that we defined
$$
W^1(\Omega) = \biggl \{ f\in L^2 (\Omega): \sum_{j=0}^N 
\l D_j f, D_j f \r _0 + \l f, f \r _0 < \infty \biggr \} .
$$
Here the derivative is intended in  the sense of distributions,
$$
\l f,g \rangle_0
            = \int_\Omega f\overline{g} \, dV ,
$$
and $D_j = \partial/\partial x_j, j = 0,\dots,N.$
Moreover we have defined the 1-Sobolev space of 
$q$-forms $W^1_q (\Omega)$,
for $q=1,\dots,N+1$, by setting
$$
W^1_q (\Omega) =
\bigl\{  \phi=\sum_{|I|=q}\phi_I dx^I : \phi_I \in W^1 (\Omega)
\bigr\} .
$$
For $\phi,\psi\in W^1_q (\Omega)$,
their inner product in $W^1_q (\Omega)$ is given by
\begin{align*}
\l \phi,\psi\r 
& = \l \sum_I \phi_I dx^I , \sum_J \psi_J dx^J \r_1 \\
& = \sum_I \l \phi_I ,\psi_I \r_1 .
\end{align*}
Throughout the rest of 
this entire paper, we shall denote the $s$-Sobolev 
norm of a form $\psi$ by $\| \psi\|_s$.  

Let
$$
 d: \bw^q
\longrightarrow \bw^{q+1}
$$
be defined by
\begin{align*}
d\bigl( \sum_I\phi_I \, dx^I \bigr)
& = \sum_I \sum_{j=0}^N D_j \phi_I \, dx_j \wedge \, dx^I \\
& =  \sum_{|J|=q+1} \biggl( \sum_{|I|=q \atop j=0,\dots,N}  
\e{J}{jI} D_j \phi_I \biggr) dx^J \\
& \equiv \sum_{|J|=q+1} \phi_J' dx^J .
\end{align*}
 We let  $d^*$ be the operator
on $\bw^{q+1}$ defined by
$$
\l d \phi,\psi \rangle_1
= \l \phi, d^* \psi \rangle_1 .
$$
Recall that
$$
\dom d^* \cap \bw^{q+1}
= \bigl\{ \psi \in \bw^{q+1}
: | \l d \phi,\psi \rangle_1 | \leq C_\psi
\|\phi \|_1 \bigr\} .
$$

Our goal is to solve the \bvp
\begin{equation}\label{BVP-domain}
\begin{cases}
(d d^* + d^* d) \phi  =  \alpha \qquad & \on \Omega\\
 \phi \in \dom d^* 	& \\
 d\phi \in \dom d^* 	&
 \end{cases} 
\end{equation}
for $\alpha\in W^1_q (\Omega)$, 
$q=0,\dots,N$, and prove existence and
regularity theorems.

We need to analyze both the equation on $\Omega$, and the boundary 
conditions. Our first goal is
to describe $\dom d^*$. Since we have 
developed some familiarity
with this type of calculation,
we work directly with $q$-forms for all $q$ 
(recall that, in the half space case, we restricted
attention at first to functions). 
We have the following result.
\begin{proposition}\label{DOM-d*}  \sl
Let $\Omega$ be a smoothly bounded domain in $\RR^{N+1}$. 
Then the Hilbert space adjoint $d^*$ of $d$, in the $W^1$ inner
product, acting on $q$-forms,
has domain satisfying
$$
\dom d^* \cap \bw^{q+1} (\overline{\Omega}) =
 \bigl \{ \psi \in \bw^{q+1}(\overline{\Omega}):
(\nabla_{\vn} \psi) \lfloor \vn \bigr |_{b\Omega} 
 	\equiv 0 \bigr \} .
$$
\end{proposition}
Here we use the notation $\nabla_X \phi$ 
to denote the covariant 
differentiation of the form $\phi$ in the direction given by the
vector field $X$, and ``$\lfloor$" is the 
standard contraction
operation from exterior algebra. 
Recall that by definition, if $Y_1,\dots,Y_q$ are vector fields, then
$$
(\nabla_X \phi) = X\bigl( \phi(Y_1,\dots,Y_q)\bigr) 
-\sum_{i=q}^{q} \phi(Y_1,\dots,\nabla_X Y_i ,\dots,Y_q) ;
$$
also, in local coordinates $(y_0,\dots, y_{n})$,
$$
\nabla_X V =\sum_{k=0}^{n} \bigl( X(V_k) +\sum_{i,j=0}^{n} 
\Gamma^k_{ij} X_i V_j 
\bigr) \pd{}{y_k} .
$$
Clearly covariant 
differentiation preserves the {\it type}  
of a form. 

Moreover,
$$
\phi\lfloor V=\sum_{iIJ}\phi_J V_i \e{J}{iJ} dy^I .
$$
Notice that, in the standard coordinates of $\RR^{N+1}$,
$(\nabla_X \phi)_K = X(\phi_K)$ and 
$$\bigl(\phi\lfloor(\p/\p x_0)\bigr)_K = \phi_{0K}$$
if $K\not\ni0$  and $\bigl( \phi\lfloor(\p/\p x_0)\bigr)_K =0$ if
$K\ni0$.  
For these and related notions we refer the reader to 
\cite{FED}.

Observe that, if $\phi=\sum_i \phi_i dx_i$ is a $1$-form, then the 
boundary condition $\nabla_{\vn} \phi\lfloor\vn \big|_{b\Omega} =0$
can be written as
\begin{equation}\label{Jan3-1}
\sum_{j=0}^{N} \pd{\rho}{x_j} \pd{\phi_i}{n} =0 \qquad \on b\Omega .
\end{equation}  
\begin{pf*}{\bf Proof of \ref{DOM-d*}} 
Let $\phi=\sum_{|I|=q}\phi_I dx^I$, 
and $\psi=\sum_{|J|=q+1}\psi_I dx^J$. 
We shall
use the following form of Green's
theorem:  
$$
\int_\Omega D_j f \overline{g} 
= -\int_\Omega f \overline{D_j g} + \int_{b\Omega} f\bar g
\pd{\rho}{x_j} .
$$
Recall that $\vn = (D_0 \rho,\dots,D_N \rho)$ is the normal
direction. 
\vfill
\eject

Let $d\phi
 =  \sum_{|J|=q+1} \bigl( \sum_{|I|=q} \sum_{j=0,\dots,N}  
\e{J}{jI} D_j \phi_I \bigr) dx^J 
 \equiv \sum_{|J|=q+1} \phi_J' dx^J$. Then
\begin{align*}
\l d\phi, \psi \rangle_1 
& = \sum_{|J|=q+1} \left[ \sum_{k=0}^N 
\l D_k \phi_J' ,D_k \psi_J \r_0 
    + \l \phi_J' ,\psi_J \r_0 \right]  \\
& = \sum_{|J|=q+1} \left[ \sum_{k=0}^N 
	\sum_{|I|=q \atop j=0,\dots,N} \e{J}{jI} \biggl( 
- \l D_k \phi_I , D_k D_j \psi_J \r_0 +\int_{b\Omega} D_k \phi_I
\overline{\psi_J} \pd{\rho}{x_j} \biggr) \right. \\
& \qquad \qquad  
\left. + \sum_{|I|=q \atop j=0,\dots,N} \e{J}{jI} \biggl(  
- \l \phi_I ,D_j \psi_J \r_0 + \int_{b\Omega} \phi_I
\overline{\psi_J} \pd{\rho}{x_j} \biggr) \right] \\
& = \sum_{|I|=q} \left[ \sum_{k=0}^{N} \l D_k \phi_I , D_k 
\bigl( -\sum_{|J|=q+1 \atop j=0,\dots,N} 
\e{J}{jI} D_j \psi_J \bigr) \r_0  + \l \phi_I ,\bigl( 
-\sum_{|J|=q+1 \atop j=0,\dots,N} \e{J}{jI} D_j \psi_J 
\bigr) \r_0  \right] \\
& \quad + \sum_{|I|=q} \left[ \sum_{k=0}^{N} \biggl(
\int_{b\Omega} D_k \phi_I \overline{ \bigl(
\sum_{|J|=q+1 \atop j=0,\dots,N} \e{J}{jI} D_k \psi_J
\pd{\rho}{x_j} \bigr)  } \biggr)
+ \int_{b\Omega} \phi_I \overline{ \bigl( 
\sum_{|J|=q+1 \atop j=0,\dots,N} \e{J}{jI}  \psi_J 
 \pd{\rho}{x_j} \bigr)  }
\right] \\  
& = \l \phi,d'\psi \r_1  \\
& \quad +\sum_{|I|=q} \left[ \sum_{k=0}^{N} \biggl(
\int_{b\Omega} D_k \phi_I \overline{ \bigl(
\sum_{|J|=q+1 \atop j=0,\dots,N} \e{J}{jI} D_k \psi_J
\pd{\rho}{x_j} \bigr)  } \biggr)
+ \int_{b\Omega} \phi_I \overline{ \bigl( 
\sum_{|J|=q+1 \atop j=0,\dots,N} \e{J}{jI}  \psi_J 
 \pd{\rho}{x_j} \bigr)  }
\right] .
\end{align*}
 Notice that $d'$ in the last equality above 
is precisely the formal
adjoint of $d$.  Recall that $d'=-\text{div}$ on 1-forms.  
It is clear that    
$$
|\l \phi, d' \psi \r_1 | \leq \| \phi\|_1 \|\psi \|_2 ,
$$
and that
\begin{eqnarray*}
\left|
\int_{b\Omega} \phi_I \overline{ \bigl( 
-\sum_{|J|=q+1 \atop j=0,\dots,N} \e{J}{jI}  \psi_J 
 \pd{\rho}{x_j} \bigr)  } \right| 
& \leq & c\cdot \|\phi\|_{L^2(b\Omega)} \|\psi\|_{L^2(b\Omega)} \\  
        & \leq & c\cdot \|\phi\|_1 \|\psi\|_1 ,
\end{eqnarray*}
by the trace theorem.
 
Next we consider the term
$$
\sum_{k=0}^{N} 
\int_{b\Omega} D_k \phi_I \overline{ \bigl(
\sum_{|J|=q+1 \atop j=0,\dots,N} \e{J}{jI} D_k \psi_J
\pd{\rho}{x_j} \bigr)  } \, .
$$
For $k=0,\dots,N$ we decompose  the differential operator $D_k$ into  
normal and tangential components, that is, in a suitable neighborhood
of $b\Omega$ we write 
\begin{equation}\label{Yk}
D_k =Y_k  +\pd{\rho}{x_k} \pd{}{n} ,
\end{equation}
where $Y_k$ is a tangential vector field, and $\p/\p n$ is the unit
vector field in the normal 
direction.  
Therefore, for each fixed $I$, $|I|=q$, integration 
by parts yields that
\begin{align*}
\sum_{k=0}^{N} 
\int_{b\Omega} D_k \phi_I \overline{ \bigl(
\sum_{|J|=q+1 \atop j=0,\dots,N} \e{J}{jI} D_k \psi_J
\pd{\rho}{x_j} \bigr)  } 
& =  \sum_{k=0}^{N} 
\int_{b\Omega}  \phi_I \overline{ Y^*_k \bigl(
\sum_{|J|=q+1 \atop j=0,\dots,N} \e{J}{jI} D_k \psi_J
\pd{\rho}{x_j} \bigr)  }  \\
& \qquad + \int_{b\Omega} \frac{\p \phi_I}{\p n}
\overline{ \bigl(
\sum_{|J|=q+1 \atop j=0,\dots,N} \e{J}{jI} \pd{\psi_J}{n}
\pd{\rho}{x_j} \bigr)  }   \\
& \equiv I+E .
\end{align*}
Now observe that on $b\Omega$, by definition, 
\begin{equation}\label{eq-cov-diff-bdry}
\sum_{|J|=q+1 \atop j=0,\dots,N} \e{J}{jI} \pd{\psi_J}{n}
\pd{\rho}{x_j} =  \bigl(\nabla_{\vn} \psi \lfloor \vn\bigr)_I
	 \bigr |_{b\Omega} .
\end{equation}
Using the trace theorem it is easy to see that, 
for $\psi$ with $C^\infty(\overline{\Omega})$ coefficients, 
$$
|I| \leq C_\psi \|\phi\|_1 .
$$
Therefore, if $ \nabla_{\vn} \psi \lfloor \vn \bigr|_{b\Omega} 
= 0$, then $\psi \in \dom d^*.$
 
Conversely, suppose that $\nabla_{\vn} \psi \lfloor \vn
\bigr|_{b\Omega} \not \equiv 0.$
Then, repeating the construction that
we used in the case of the half
space, we can see that the mapping
$$
\phi \mapsto \sum_I
\int_{b\Omega} \pd{\phi_I}{n} \overline{ \bigl(
\nabla_{\vn} \psi \lfloor \vn \bigr)_I }
$$
cannot be continuous on $W^1_q (\Omega).$  This concludes the proof. 
\endpf

Notice that, from the previous computation, we obtain that for
$\psi\in\dom d^*$ it holds that 
\begin{multline}\label{Jan3-2}
\l d\phi,\psi\r_1
 = \l \phi,d'\psi\r_1 \\
+ \sum_{|I|=q} \biggl[ -\sum_{k=0}^{N} 
\int_{b\Omega} \phi_I 
\overline{ \nabla_{Y^*_k}\bigl(
\sum_{|J|=q+1 \atop j=0,\dots,N} \e{J}{jI} D_k \psi_J
\pd{\rho}{x_j} \bigr)  }
+ \int_{b\Omega} \phi_I \overline{(\psi\lfloor \vn)_I}
\biggr]  . 
\end{multline}

Next we shall determine an explicit
expression for $d^*.$   
We set (as in the case of the half space) 
\begin{equation}
d^*\phi = d'\phi + \K \phi ,     \label{DAGG}
\end{equation}
where $\K $ is to be determined.
\begin{proposition}\label{esp-K-Omega}  \sl
Let $\psi\in\dom d^* \cap \bw^{q+1} (\overline{\Omega})$.  Then the  
$I$-component $(\K \psi)_I$ of $\K \psi$ is the solution of the \bvp 
$$
\begin{cases}
\displaystyle{ 
-\btu v +v = 0 }&\on \Omega \\
\displaystyle{ 
\pd{v}{n} = \sum_{k=0}^{N}\biggl[
\nabla_{Y^*_k}\bigl(
(\nabla_{Y_k} \psi) \lfloor\vn \bigr)_I 
+(\psi\lfloor \vn)\biggr]_I  } 
		& \on b\Omega
\end{cases} .
$$
\end{proposition} 
\bgpf
By the relation
$$
\l d\phi, \psi \r_1 = \l \phi,d'\psi\r_1 + \l \phi,\K\psi\r_1 ,
$$
and equation (\ref{Jan3-2}) we see that for $\psi\in\dom d^* \cap
\bw^{q+1} (\overline{\Omega})$,
\begin{equation}\label{boundary-integral} 
\l \phi,\K\psi \r_1 
= \sum_{|I|=q} \biggl[ \sum_{k=0}^{N} 
\int_{b\Omega} \phi_I 
\overline{ Y^*_k \bigl(
\sum_{|J|=q+1 \atop j=0,\dots,N} \e{J}{jI} D_k \psi_J
\pd{\rho}{x_j} \bigr)  }  
+ \int_{b\Omega} \phi_I \overline{(\psi\lfloor \vn)_I}
\biggr]  . 
\end{equation}

Now, by Green's theorem, we see that 
\begin{align*}
\l \phi ,\K\psi \r_1 
& = \sum_{|I|=q} \left[ \sum_{j=0}^{N}
\int_\Omega  D_j \phi_I  
\overline{D_j (\K\psi)_I }
+\int_\Omega \phi_I \overline{(\K\psi)_I}  \right] \\
& = \sum_{|I|=q} \left[ -\int_\Omega \phi_I 
\overline{\btu (\K\psi)_I } +\int_{b\Omega} \phi_I 
\overline{ \pd{(\K\psi)_I}{n} }
+\int_\Omega \phi_I \overline{(\K\psi)_I}  \right] .
\end{align*}
Therefore, for each fixed multi-index $I$, we must have
\begin{multline*}
\int_\Omega \phi_I \overline{ \bigl[
-\btu (\K \psi)_I + (\K \psi)_I \bigr]} +\int_{b\Omega} \phi_I 
\overline{\pd{(\K \psi)_I}{n}} \\
= \sum_{k=0}^{N} 
\int_{b\Omega} \phi_I 
\overline{ Y^*_k \bigl(
\sum_{|J|=q+1 \atop j=0,\dots,N} \e{J}{jI} D_k \psi_J
\pd{\rho}{x_j} \bigr)  }  
+ \int_{b\Omega} \phi_I \overline{(\psi\lfloor \vn)_I } .
\end{multline*} 
This implies that 
$$
\bigl[ -\btu (\K \psi)_I + (\K \psi)_I \bigr] = 0 \qquad \on \Omega.
$$
Now recalling that, since $\psi\in\dom d^*$, 
$$
\sum_{|J|=q+1\atop j=0,\dots,N}
\e{J}{jI}\pd{\rho}{x_j}\pd{\psi_J}{n}
= 0 \qquad \on b\Omega \, .
$$
By writing $D_k = Y_k  
+(\p \rho /\p x_k)\cdot(\p/\p n)$, we have that 
\begin{align*}
\sum_{k=0}^{N} Y^*_k \bigl(
\sum_{|J|=q+1 \atop j=0,\dots,N} \e{J}{jI} D_k \psi_J
\pd{\rho}{x_j} \bigr)  
& = \sum_{k=0}^{N} Y^*_k \bigl(
\sum_{|J|=q+1 \atop j=0,\dots,N} \e{J}{jI} Y_k \psi_J
\pd{\rho}{x_j} \bigr)  \\
& = \sum_{k=0}^{N} Y^*_k \bigl(
\nabla_{Y_k} \psi \lfloor \vn \bigr)_I  \\
& = \sum_{k=0}^{N} \bigl[
\nabla_{Y^*_k} (\nabla_{Y_k}\psi\lfloor \vn)\bigr]_I.
\end{align*}
From this the result follows.
\endpf
Thus $\K$ is the solution operator  for the
preceding elliptic \bvp.  By standard facts of the 
theory of elliptic problems (see \cite{HOR1})
we know that  $\K\psi$ is
uniquely determined for each $\psi \in 
\bw^q(\overline{\Omega}).$  

As in the case of the half space, the operator $\K$ 
turns out to be of order 1.
In order to prove 
the following corollary,  
we need to
apply Theorem 4.2.4 in \cite{TRI}.  This result deals with  
Schauder  
estimates for elliptic \bvp s formulated in 
norms of negative order. 
The reference \cite{TRI}
happens to use the language 
of Triebel-Lizorkin spaces $F^s_{p,q}$.
Since we consider only the case $p=q=2$ 
we shall denote the
spaces $F^s_{2,2}$ simply 
by $F^s$.
We also recall that, for $s\ge 0$, $F^s = W^s (\Omega)$.
\begin{corollary}\label{K-ord-1}  \sl
Let $s>1/2$. Then there exists $C>0$ such that
for all $\psi\in\bw^q (\overline{\Omega})$ we have
$$
\| \K\psi \|_{F^{s-1}} 
\le C \| \psi\|_{s} .
$$
\end{corollary}
\bgpf
Notice that by this result the operator $\K$, initially defined on
the dense subspace $\bw^q (\overline{\Omega})$ 
can be extended as a continuous operator from 
$F^{s-1}$ to $W^s (\Omega)$.

We need only apply standard estimates for elliptic \bvp s.
In order to obtain the sharp regularity result we use the estimates
for the negative Sobolev spaces on the boundary (see \cite{TRI},
Theorem 4.2.4).  Notice that the boundary 
differential operator acting on the
data is of order 2, and tangential.  Write $T_2$ for that
operator.
Then we have
\begin{align*}
\| \K \psi \|_{F^{s-1}}
& \le C\bigl( \| T_2 \psi \|_{W^{s-5/2}(b\Omega)} 
	+ \| \psi \|_{W^{s-5/2}(b\Omega)} \bigr)\\
& \le C \| \psi\|_{W^{s-1/2} (b\Omega)}  \\
& \le C \| \psi\|_{s}. \qed
\end{align*} 
\renewcommand{\qed}{}\endpf
We could, in principle, write down an approximate
expression for $\K$ by using local coordinates to
reduce the problem to the half space situation, where
we have calculated $\K$ quite explicitly.  We forego
that option for now.  

We can now reformulate our \bvp\ (\ref{BVP-domain}), recalling that 
on forms $dd' +d'd =\btu$:
$$
\begin{cases}
\displaystyle{
(-\btu +G_\Omega)\psi=\alpha} \qquad & \on \Omega\\
\displaystyle{
\nabla_{\vn}\psi\lfloor\vn =0} & \on b\Omega \\
\displaystyle{
\nabla_{\vn}d\psi\lfloor\vn =0} & \on b\Omega 
\end{cases} 
$$
where we set $G_\Omega =d\K +\K d$.


\section{A Special Coordinate System}\label{SPECIAL}
\setcounter{equation}{0}

In  this section we introduce a system of local
coordinates near the boundary that we will use in the rest of the
paper.

It is a standard fact (see \cite{KR2} for instance) that since
$b\Omega$ is smooth and compact there exists a tubular neighborhood
$U$ of $b\Omega$ such that for each $x\in U$ there is a unique $\pi
(x)\in b\Omega$ which realizes the distance of $x$ from $b\Omega$.
The line joining $\pi(x)$ and $x$ is orthogonal to $b\Omega$.
\begin{definition}\label{FERMI}  \rm
Let $\bigl\{ (U_j,\Phi_j)\bigr\}$ be a covering of $b\Omega$ by local
coordinates patches.  The pair $(V_j,\Psi_j)$, $V_j\subset U$, and 
$\Psi:V_j \rightarrow \RR^{N+1}$ is a {\it Fermi coordinate patch} on $U$
if 
$$
\pi(x)\in U_j \qquad \text{for all} \ p\in V_j \subset U, 
$$
and 
$$
\Psi_j (p)=
\bigl( \text{dist\,}( \pi(p)),p\bigr),\Phi_j (\pi(p)) \bigr).
$$ 
Throughout the rest of the paper, $X_0 =\vn$ will denote the vector
field defined at each $p\in U$ that is given by 
$$
\vn_p = \vn_{\pi(p)}=\bigl( \pd{\rho}{x_0} (p),\dots,
\pd{\rho}{x_n}(p) \bigr).
$$
\end{definition}
We remark that the equality 
$$
\vn_{p}= \bigl( \pd{\rho}{x_0}(p),\dots,\pd{\rho}{x_n}(p)\bigr)
$$
in general does not hold for $p\in U\setminus b\Omega$.

On the fixed coordinate patch $(U_j ,\Psi_j)$ 
we also have an orthonormal frame of
vector  
fields $X_0 = (\p/\p n), X_1,\dots, X_N$.  Notice that 
$X_1,\dots,X_N$ form an orthonormal frame for the tangential vector
fields.  We also fix the dual basis of $1$-forms, $\omega_0 = dx_0,
\omega_1,\dots,\omega_N$. Given any $q$-form $\phi$ we write 
$\phi=\sum_{|I|=q}\omega^I$, where $\omega^I
=\omega_{i_1}\wedge\dots\wedge\omega_{i_q}$, where 
$I=(i_1,\dots,i_q)$ is a strictly increasing multi-index. 

In the rest of the paper we shall {\em always} assume 
that the forms
are written in terms of this given basis.
Notice that, in these coordinates,
\begin{align*}
d\psi & =d(\sum_I \psi_I \omega^I) \\
& = \sum_j X_j \psi_I \omega_j \wedge \omega^I +\{
0\text{ order terms in }\psi \} \\
& = \sum_J \bigl( \sum_{jI} X_j \psi_I \e{J}{jI} \bigr)\omega^J 
+ \{ 0\text{ order terms}\} . 
\end{align*}
Analogously,
$$
d' \psi = \sum_L \bigl( \sum_{jI} X'_j \psi_I 
\e{I}{jL} \bigl) \omega^L +\{ 0\text{ order terms}\} . 
$$
Here $X_k'$ denotes the formal adjoint of the
vector field $X_k$.
 
Recall that (see \cite{hel} for instance) we can 
define the (Christoffel symbol) coefficients 
$\Gamma^{k}_{ij}$ as follows: 
\cite{hel}:
\begin{equation}\label{def-Christoffel}
\nabla_{X_i}X_j = \sum_{k=0}^{N} \Gamma^{k}_{ij} X_k .
\end{equation}
We remark that $\nabla_{\vn}\vn=0$, 
since $\vn$ is, by definition, a unit vector field whose integral
curves are lines.
Therefore the 
coefficients $\Gamma^{k}_{ij}$ satisfy the following relations
\begin{equation}\label{Christoffel}
\Gamma^{k}_{ij} = 0 \quad 
\begin{cases}
\text{if\  } k=0 &\text{and } i, \text{\ or\ }j=0 \\
\text{if\  } i=0 &\text{and } j=0 
\end{cases} 
\end{equation}
This is so because 
$$
\Gamma_{0j}^{0} \equiv\l \nabla_{X_0} X_j ,X_0 \r
=X_0 \l X_j ,X_0 \r - \l X_j , \nabla_{X_0} X_0 \r =0 ,
$$
$$
\Gamma^k_{00} \equiv \l \nabla_{X_0} X_0, X_k \r=0,
$$
and
$$
\Gamma^0_{i0}\equiv\l \nabla_{X_i}X_0, X_0 \r 
= \frac12 X_i \l X_0 ,X_0 \r =0 .
$$ 

Now we have a 
lemma about covariant differentiation in the normal direction. 
\begin{lemma}\label{covariant-diff}   \sl
If $\phi\in\bw^q (\overline{\Omega})$ is given by $\phi_I =\sum_I
\phi_I \omega^I$, 
 then
$$
\nabla_{X_0} \phi = \sum_{|I|=q} \left[ \pd{\phi_I}{x_0} 
+ \sum_{|J|=q} \gamma_{IJ} \phi_J \right] \omega^I ,
$$
where 
$$
\gamma_{IJ} =0 \begin{cases} 
\text{for } J\ni 0 & \text{if } I\not\ni 0 \\
\text{for } J\not\ni 0 & \text{if } I\ni 0 
\end{cases} 
$$
\end{lemma}
\bgpf
Fix a Fermi coordinate chart, then for $\phi\in\bw^q
(\overline{\Omega})$ supported in a small open set we have 
\begin{align*}
\left( \nabla_{X_0} (\sum_J \phi_J \omega^J )\right)
 ( X^I ) 
& = \sum_J \left( \nabla_{X_0} (\phi_J \omega^J ) \right)
 (X^I ) \\
& = \sum_J \pd{\phi_J}{x_0} \e{I}{J} -\sum_{J,J'} \sum_{s=1,\dots,q}
\phi_J
\bigl( (\nabla_{X_0}\omega^{j_s})\wedge \omega^{J'} 
\e{J}{j_s J'}\bigr)  (X^I ) .
\end{align*}
Moreover,
\begin{align*}
\bigl( \nabla_{X_0} (\omega^j) \bigr) (X_\ell ) 
& =  -\omega^j \bigl( \nabla_{X_0} X_\ell \bigr) \\
& = -\omega^j \bigl( \sum_k \Gamma^{k}_{0\ell} \, X_k \bigr) \\
& = - \Gamma^{j}_{0\ell} .
\end{align*}
Therefore
$$
\nabla_{X_0} \omega^{j_s} = -\sum_\ell \Gamma^{j_s}_{0\ell} 
\, \omega^\ell .
$$
Hence
$$
\bigl( \nabla_{X_0} \phi \bigr) (X^I )
= \sum_{|J|=q} \left[ \pd{\phi_J}{x_0} \e{I}{J} 
- \sum_{|J'|=q-1 \atop s=1,\dots,q}\sum_{\ell=0,\dots,N} 
\e{J}{j_s J'} (-\Gamma^{j_s}_{0\ell}\, ) \e{I}{\ell J'} 
\phi_J \right] .
$$
Thus 
$$
\nabla_{X_0} \phi  =
\sum_{|I|=q} \bigl( \pd{\phi_I}{x_0} + \sum_{|J|=q} 
\gamma_{IJ}\phi_J 
\bigr) \omega^I ,
$$
where
$$
\gamma_{IJ} = \sum_{|J'|=q-1 \atop s=1,\dots,q}\sum_{\ell=0,\dots,N}  
\e{J}{j_s J'} \Gamma^{j_s}_{0\ell}\, \e{I}{\ell J'} .
$$
Recall that (see (\ref{Christoffel})) 
the symbols $\Gamma^k_{ij}$ are zero if either $k=0$ and
either $i = 0$ or $j=0$, or $i=j=0$.  
Suppose that $I\ni 0$ and that $J\not\ni 0$.  Then $\ell =0$ so that 
$\Gamma^{j_s}_{ij}\, =0$. If $I\not\ni 0$ and $J\ni 0$, 
then $j_s =0$ 
so that $\Gamma^0_{0\ell}\, =0$. This proves the lemma.
\endpf

In the sequel we will also use the following observation on the
Laplace operator acting on forms.  In general, the 
Laplacian on forms 
is defined as $dd'+d'd$.  
Having chosen the aforementioned basis on the space
of $q$-forms, we see that if $\psi=\sum_{|I|=q}\psi_I \omega^I$ then
$$
\btu \psi= \sum_{I} \biggl[ \sum_{k=0}^{N} (X_k' X_k + X_k X_k') 
\psi_I \biggr] \omega^I +\{ \text{lower order terms in } \psi \}.
$$

\subsection*{The tangential Laplacian}  Fix a Fermi coordinate patch
$\bigl(U,\Psi=(x_0,\dots,x_N)\bigr)$.  (Notice that $(x_0,\dots,x_N)$
are not the standard coordinates of $\RR^{N+1}$.)  Given a function
(or a form)
$u$ defined on $U$, we denote by $\tilde u$ the function
$u\circ\Psi^{-1}$ 
defined on $\Psi(U)$.  

In these coordinates the standard Laplacian of $\RR^{N+1}$ has the
following form:
\begin{align}
[\btu u]\til   
&  = \sum_{j,k=0}^{N} \frac{1}{\sqrt{\det g}} \pd{}{x_j} \bigl(
g^{jk}\pd{}{x_k} \tilde u \bigr) \notag \\ 
& = \sum_{j,k=1}^{N} \frac{1}{\sqrt{\det g}} \pd{}{x_j} \bigl(
g^{jk}\pd{}{x_k} \tilde u  \bigr)
+ \frac{1}{\sqrt{\det g}} \pd{}{x_0} \bigl(
\det g \pd{}{x_0} \tilde u \bigr) \notag \\
& = \btu' \tilde u + \pd{}{x_0} \log (\sqrt{\det g})
\pd{\tilde u}{x_0} + \pd{^2 \tilde u}{x_0^2} \notag \\
& = [\btu \tilde u] +\ep {\cal T}_2 \tilde u + \pd{}{x_0} 
(\log \sqrt{\det g})\pd{\tilde u}{x_0} \notag \\
& \equiv [\btu \tilde u] +\ep \L_2  \tilde u ,  \label{Delta-tilde}
\end{align}
where $g$ is the metric matrix, and $(g^{jk})$ is its inverse.  The
operator $\btu_T$ is the Laplacian on the submanifold obtained by
fixing $x_0$, so that 
\begin{equation} \label{tangential-Delta}
\btu_T = \sum_{j=1}^{N} \pd{^2}{x_j^2} + \ep {\cal T}_2 ,
\end{equation}
where ${\cal T}_2$ is a second order tangential differential operator
with $C^\infty (\overline{\Omega})$ coefficients, and $\ep$ can be
made arbitrarily small by shrinking $U$.  

Notice that we have obtained that, on $U$,
\begin{equation}\label{Delta_T}
\btu=\btu_T +X_0^2 + \bigl( X_0 \log \sqrt{\det g} \bigr) X_0 .
\end{equation}
\subsection*{The operator $G_\Omega$ on functions}  
Recall that, by
Proposition \ref{esp-K-Omega} $G_\Omega u$
is the unique solution of the elliptic \bvp 
$$
\begin{cases}
\displaystyle{ \btu v - v =0 } & \on \Omega\\
\displaystyle{ \pd{v}{n} = \sum_{k=0}^{N} Y_k^* \bigl( 
(\nabla_{Y_k} du)\lfloor\vn \bigr) + \pd{u}{n} 
					}\qquad & \on b\Omega
\end{cases} \quad .
$$
We are interested in making explicit the boundary operator in this
case.  Notice that $Y_k^* = - Y_k + (0\text{\ order terms})$, so that the
boundary equation becomes
$$
-\sum_{j=0}^{N} Y_k^2 (X_0 u) +T_1 X_0 u + T_2 u,
$$
where the $T_j$'s are tangential differential operators of order $j$.
Now notice that 
$$
\sum_{j=0}^{N} Y_k^2 = \sum_{k=0}^{N} D_k^2 -X_0^2 +L_1
$$
where $L_1$ is a first order differential operator.  This follows
easily from the fact that, by construction,
$$
Y_k = D_k -\l D_k ,\vn\r X_0 ,
$$
(see (\ref{Yk})).  In fact, for a function $f$,
\begin{align*}
\btu f 
& = \sum_{k=0}^{N} D_k^2 f \\
& = \sum_{k=0}^{N} (Y_k +\l D_k ,\vn \r)^2 f \\
& = \sum_{k=0}^{N} 
\bigl( Y_k^2 f 
+ Y_k (\l D_k ,\vn \r X_0 )f 
	+X_0 (\l D_k ,\vn \r Y_k )f \bigr) +X_0^2 f  \\
& =  \sum_{k=0}^{N} \bigl( Y_k^2 f  
+ Y_k (\l D_k ,\vn \r )X_0 f \bigr) +X_0^2 f  
	 \\
& = \sum_{k=0}^{N} Y_k^2 + X_0^2 f +aX_0 f \, ,
\end{align*}
since $\sum_{k=0}^{N}\l D_k ,\vn \r Y_k =0$.
Therefore the boundary equation in $u$ equals
\begin{align*}
\bigl( -\sum_{k=0}^{N}Y_k^2 \bigr) (X_0 u)
+T_1 X_0 u +T_2 u
& = (-\btu +X_0^2 )(X_0 u) +aX_0^2 u + T_1 X_0 u +T_2 u \\
& = (-\btu_T +bX_0 )(X_0 u) + a X_0^2 u + T_1 X_0 u + T_2 u \\
& = -\btu_T X_0 u + T_1 X_0 u +T_2 u ,
\end{align*}
where we have used formula (\ref{Delta_T}) and the fact that
$du\in\dom d^*$.  In local coordinates the boundary equation then
becomes 
\begin{equation}\label{LuigiXIV}
\pd{v}{x_0} = -\btu_T(X_0 u)\til +\ep {\cal T}_2 (X_0 u)^2 
+ \bigl[ T_1 X_0 u + T_2 u \bigr]\til .
\end{equation}

We conclude this section by introducing a convention that we shall
use consistently throughout the rest of the paper.  By $L_j$ 
we denote a generic
differential operator of order $j$
with $C^\infty (\overline{\Omega})$
coefficients, while by $T_j$ 
we denote a  differential operator of order $j$
with $C^\infty (\overline{\Omega})$ (defined in a suitable
neighborhood of $b\Omega$),
that involves only tangential derivatives.

\section{The Existence Theorem}
\setcounter{equation}{0}
In this section we study the question
of existence  for the
\bvp
$$
\begin{cases}
(d d^* + d^* d) \phi  =  \alpha \qquad & \on \Omega\\
 \phi \in \dom d^* 	& \\
 d\phi \in \dom d^* 	&
 \end{cases} 
$$
for $\alpha\in\bw^q(\overline{\Omega})$, 
$q=0,\dots,N$, 
from an abstract Hilbert space point of view.  
We remark that by Proposition \ref{DOM-d*} the
boundary conditions are given by the equations
\begin{align*}
\nabla_{\vn} \phi \lfloor \vn &= 0 \qquad  \on b \Omega \\
\intertext{and}
\nabla_{\vn} d\phi \lfloor \vn & = 0  \qquad \on b \Omega .
\end{align*}
Notice that, using the coordinates and the frame introduced in 
(\ref{FERMI}), the above equations can be rewritten as
\begin{align} \label{XXXX}
X_0 \phi_I + L_0 \phi =0 \qquad & \on b\Omega \quad \text{if }
I\ni0\\
\intertext{and} \label{YYYY}
X_0^2 \phi_I + L_1 \phi =0 \qquad & \on b\Omega \quad \text{if } 
					I\not\ni0, 
\end{align}
where $L_j$ are differential operators of order 
$j$ in the components
of $\phi$.  We remark that if $\phi$ is a function then
the first boundary
equation is empty, and the second one becomes
\begin{equation}\label{ZZZZ}
X_0^2 \phi =0 \qquad \qquad \on b\Omega ,
\end{equation}
since $\nabla_{\vn}\vn=0$.

Our development in what follows parallels classical studies such
as that which can be found in \cite{FOK}.
Let 
$$
{\cal D} \equiv \bigl \{ \phi \in \bw^q (\overline{\Omega})
: \phi \, , \, d\phi \in \dom d^*\bigr \} .
$$
For $\phi,\psi \in {\cal D}$ we define the bilinear form
$$
Q(\phi,\psi) = \l d\phi, d\psi \r_1 + \l d^* \phi, d^* \psi \r_1
                              + \l \phi, \psi \r_1 .
$$

Our first claim is that $\D$ is dense in $W_q^1$.
We argue as follows.
Let $\phi$ be any $q$-form.  We may assume that $\phi$ has
coefficients  
smooth up to the boundary.  Then it suffices
to find a $\psi \in \bw^q (\overline{\Omega})$ of
small norm such that $\phi + \psi \in \D$.
We use Fermi local coordinates in a tubular neighborhood of
$b\Omega$, so that 
$x = (x_0,x')$ with $x' \in b\Omega$ and
$x_0=\text{dist\,} (x_0,b\Omega)$ 
parametrizing the normal direction.  Set
$$
\psi_1(x') = \nabla_{\vn}
\bigl( \phi \lfloor \vn \bigr) \biggr|_{b\Omega} , 
$$
and
$$
\psi_2 (x') =  \nabla_{\vn} \bigl( d\phi \lfloor \vn
        \bigr ) \biggr |_{b\Omega} .
$$
Finally, set
$$
\psi(x) = \left(-x_0 dx_0\wedge \psi_1(x') -\frac12 x_0^2 
\psi_2 \right) \chi (x_0) ,
$$
where $\chi \in C_0^\infty[-2\epsilon,2\epsilon]$, 
$\chi= 1$ in a neighborhood
of 0, and
$\|\chi\|_1 \leq C  \epsilon^{-1/2}$. 
Then $\|\psi\|_1 < C \epsilon^{1/2}$ 
and $\phi + \psi \in \D.$  That completes the argument.

Now we let $\tilde{\D}$ be the closure of $\D$ in the topology
induced 
by $Q$.  We wish to check that $\tilde{\D}$ is still
contained in $W^1_q$.
It is easy to see that $d$ is closed in $W^1_q$; of course
 $d^*$ is closed also (since adjoints are always closed).
These facts imply that $\tilde{\D}$ is a subspace of $W^1_q$.

At this point we apply the Friedrichs extension theorem, as in
\cite{FOK}, to show that there exists a 
canonical self adjoint operator
$$
T: W^1_q \rightarrow \tilde{\D}
$$
which is bounded in the $W^1$ topology, is injective, and such that
$$
Q(T\phi,\psi) = \l \phi, \psi \r_1 .
$$
If we set $F = T^{-1}$, then
$$
Q(\phi,\psi) = \l F\phi, \psi \r_1 
\qquad \qquad \forall \phi \in \dom F, \psi \in \tilde\D .
$$
Notice that $F = (d^* d + d d^*) + I$ when restricted to $\D$.

Now we show that the $Q$-unit ball of $\tilde\D$ is compactly
embedded in $W^1_q$ by
proving that the $Q$ norm is equivalent to the $W^2$ norm. Notice
that $Q(\phi,\phi) \leq c \|\phi \|_2^2$ since $d^*$ 
is of order $1$.
The reverse inequality follows from the next theorem.
\begin{theorem}\label{coerciveness}  \sl
There exists a constant $C_0 >0$ such that, for all 
$\psi\in\tilde\D$,
$$
Q(\psi,\psi)\ge C_0 \cdot \|\psi\|_2^2 .
$$
\end{theorem}
Assume the theorem for now.  
Since the $Q$-unit ball of $\tilde\D$ is compactly embedded in 
$W^1_q (\Omega)$,
it follows that $T$ is a compact operator.
Notice that if $T\alpha\in\D$, then 
$T\alpha$ is the unique solution of  the \bvp
\begin{equation*}
\begin{cases}
\displaystyle{
(d d^* + d^* d) \phi + \phi  = }& \alpha \\
\phi \in \dom d^* & 	\\
 d\phi \in \dom d^* & 
\end{cases} 
\end{equation*}
for $\alpha \in W^1_q$, $q = 0,1, \dots, N+1.$

We now wish to establish conditions for the solvability of
\begin{equation} \label{STAR}
\begin{cases}
\displaystyle{
(d d^* + d^* d) \phi   = }& \alpha \\
\phi \in \dom d^* & 	\\
 d\phi \in \dom d^* & 
\end{cases} 
\end{equation}
Let $Q_0$ be the bilinear form on $\tilde{\D}$ defined by
$$
Q_0(\phi, \psi) 
= \l d\phi, d\psi \r_1 + \l d^* \phi, d^* \psi \r_1 .
$$
Thus $\phi$ is a solution of (\ref{STAR}) precisely when
$\phi\in\D$ and 
$$
Q_0(\phi,\psi) = \l \alpha, \psi \r_1 \qquad \qquad 
\qquad \forall \psi \in \tilde{\D} .
$$
This in turn holds if and only if 
\begin{eqnarray*}
Q(\phi, \psi) & = & Q_0(\phi,\psi) + \l \phi, \psi \r_1   \\
              & = & \l \alpha, \psi \r_1 + \l \phi, \psi \r_1 \\ 
              & = & \l \alpha + \phi, \psi \r_1 .
\end{eqnarray*}

By the Friedrichs theorem, we have reduced our situation to solving
the equation
\begin{equation}\label{sTAr}
(F-I)\phi=\alpha
\end{equation}
with $\phi\in\dom F$;
i.e., setting $\theta=F\phi$,
$$
\theta-T \theta = \alpha .
$$
We now apply the standard theory of compact operators
to obtain that the above equation  has a solution $\theta$ 
for all $\alpha$ orthogonal to the finite dimensional subspace 
${\cal H}_q \equiv\ker(I - T)$. 
It is easy to check that
$\ker(I - T)$ is exactly
the kernel of $F-I$.
Thus, for $\alpha$ orthogonal (in the $W^1$-inner product)
to $\ker(I-T)$, we obtain that $\phi=T\theta$ is the solution of the 
equation (\ref{sTAr}).
Notice that, if $\phi\in\D$, then the equation (\ref{sTAr}) reduces to the
\bvp\
(\ref{STAR}). 
Moreover, $Z\equiv (F-I)(\D)$ is  a dense subspace of 
${\cal H}_q^\perp \ss W^1_q$, i.e.
$$
Z^\perp = {\cal H}_q .
$$
Let $\beta\in\D$.  For $\phi\in\D$ we have
\begin{align*}
\l (F-I)\phi,\beta\r_1 
& = Q_0 (\phi,\beta) \\
& = Q(\phi,\beta)-\l \phi,\beta\r_1 \\
& = \l \phi,(T-I)\beta\r_1 .
\end{align*}
Thus, if $\beta\in Z^\perp$, then $\l (F-I)\phi,\beta\r_1 =0$ for all
$\phi\in\D$.  Since $\D$ is dense in $W^1_q$, $(T-I)\beta=0$, that is
$\beta\in{\cal H}_q$.  

Thus we have proved the following
theorem: 

\begin{theorem}[Existence for all $q$]\label{EXISTENCE}   \sl
The \bvp 
$$
\begin{cases}
\displaystyle{
(d d^* + d^* d)\phi  =  \alpha } \qquad & \\
\phi \in \dom d^*  	& \\
d\phi \in \dom d^* 	& 
\end{cases}
$$
has a finite dimensional kernel ${\cal H}_q^{\perp}$ 
and finite dimensional cokernel.
[Note that the space $Z$ is dense in ${\cal H}_q^{\perp}$.]
The problem has a solution $\phi\in\D$ for 
$\alpha \in Z\ss {\cal H}_q^\perp$.
\end{theorem}
\subsection*{The coercive estimate}
We now prove the fundamental estimate from below for the bilinear
form $Q$.
\begin{pf*}{\bf Proof of \ref{coerciveness}}
We begin by noticing that 
\begin{align*}
Q(\psi,\psi)
& = \l d\psi,d\psi\r_1 +\l d^* \psi,d^* \psi\r_1 
	+\l \psi,\psi\r_1 \\
& = \l  d\psi,d\psi\r_1 +\l d' \psi,d^* \psi\r_1 
	+\l \K\psi, d^* \psi \r_1 +\l  \psi,\psi\r_1 \\ 
& = \l d\psi,d\psi\r_1 +\l d' \psi,d' \psi\r_1  
	+\l\psi,\psi\r_1 
	+\l d'\psi, \K\psi\r_1 
	+\l\K\psi,d^* \psi\r_1 .
\end{align*}
Our plan is to prove that the following claims hold true.
\medskip \\
{\bf Claim 1.} 
There exists a constant $C_1 >0$ such that for all 
$\psi\in \D$ we have 
$$
\l d\psi,d\psi\r_1 +\l d' \psi,d' \psi\r_1  
	+\l\psi,\psi\r_1 \ge C_1 \|\psi\|_2^2 .
$$
\mbox{}\smallskip\\
{\bf Claim 2.}
For any $\ep >0$ there exists $C_{\ep}>0$ such that all 
$\psi\in{\cal D}$ 
$$
|\l \K \psi,d^* \psi \r_1 | \le \ep \|\psi\|_2^2 
+ C_{\ep} \|d^* \psi\|_1^2 ,
$$
and
$$
|\l d' \psi, \K \psi \r_1 | \le \ep \|\psi\|_2^2 
+ C_{\ep} (\|d\psi\|_1^2 +\|\psi\|_1^2).
$$
\mbox{}\smallskip\\
Assuming the claims for now, we shall finish the proof.  We have that
\begin{align*}
Q(\psi,\psi)
& = \l d\psi,d\psi\r_1 +\l d' \psi,d' \psi\r_1 
	+\l\psi,\psi\r_1 
	+\l d'\psi, \K\psi\r_1 
	+\l\K\psi,d^* \psi\r_1 \\
& \ge C_1 \|\psi\|_2^2 - | \l \K\psi,d^* \psi \r_1 |
	-| \l d' \psi,\K\psi\r_1 | \\
& \ge (C_1 -2\ep) \| \psi\|_2^2 -C_\ep (\|d^*\psi\|_1^2
	+\|d\psi\|_1^2 +\|\psi\|_1^2 ) .
\end{align*}
Therefore the constant $C_0 =(C_1 -2\ep)/(1+C_\ep )$ does the job. 
Thus it remains to prove the claims. We begin with Claim 2.
\medskip\\
{\em Proof of Claim 2.\/}
Recall that, by formula (\ref{boundary-integral}),
$$
\l \K \theta,\phi \r_1 =\sum_I 
\left[ \sum_{k=0}^{N} \int_{b\Omega} 
(\nabla_{Y_k}(\theta\lfloor \vn))_I 
\overline{(\nabla_{Y_k}\phi)_I}
+\int_{b\Omega} (\theta\lfloor \vn)_I
\overline{\phi_I } \right] .
$$
Therefore, by applying the Schwarz inequality for Sobolev spaces on  
$b\Omega$, and recalling that the $Y_k$ are tangential vector fields,
we see that
\begin{align*}
|\l \K \psi,d^* \psi \r_1 | 
& \le \sum_I  \sum_{k=0}^{N} \int_{b\Omega}
|(\nabla_{Y_k}\psi\lfloor\vn)_I (\nabla_{Y_k} 
	(d^* \psi))_I | 
+ \| \psi\|_{L^2 (b\Omega)} \|d^* \psi\|_{L^2 (b\Omega)} \\
& \le c \sum_{kI} \|\nabla_{Y_k} (\psi \lfloor\vn)_I
					\|_{W^{1/2}(\Omega)}
\| \nabla_{Y_k} (d^* \psi)_I\|_{W^{-1/2}(\Omega)} 
+ \| \psi\|_{L^2 (b\Omega)} \|d^* \psi\|_{L^2 (b\Omega)} \\
& \le c\bigl( \|\psi\|_{W^{3/2}(b\Omega)} 
	\|d^* \psi\|_{W^{1/2}(b\Omega)} \bigr) \\ 
& \le c \|\psi\|_2 \|d^* \psi\|_1 \\
& \le \ep \| \psi\|_2^2 + C_{\ep} \| d^* \psi\|_1^2 .
\end{align*}
On the other hand, for $\psi\in{\cal D}$, we see that (using
``l.o.t.'' to denote lower order terms) 
\begin{align*}
\l d'\psi, \K\psi \r_1 
& = \sum_I \left[ \sum_{k=0}^{N} \int_{b\Omega}
	(\nabla_{Y_k} (d'\psi))_I
\overline{\bigl( (\nabla_{Y_k} \psi)\lfloor\vn\bigr)_I } \right] 
	+ \{ \text{l.o.t.} \} \\
& = \sum_{kI} \int_{b\Omega}
\bigl( \sum_{j\neq0,J} X_j Y_k \psi_J \e{J}{jI} \bigr)
\overline{(\nabla_{Y_k} \psi\lfloor\vn)_I } 
	+ \{ \text{l.o.t.} \} \\
& = -\sum_{kJ} \int_{b\Omega} (Y_k \psi_J)
\overline{ \bigl( \sum_{j\neq0,I} \e{J}{jI} X_j 
(\nabla_{Y_k} \psi\lfloor\vn)_I \bigr)} 
+ \{ \text{l.o.t.} \} \\
& = -\sum_{kJ} \int_{b\Omega} (Y_k \psi_J)
\overline{ Y_k
\bigl( \sum_{j\neq0,I} \e{J}{jI} X_j 
(\psi\lfloor\vn)_I \bigr)} 
+ \{ \text{l.o.t.} \} , 
\end{align*}
where we have used the fact that the term containing $X_0 \psi_J$
with $J\ni 0$ can be absorbed in the error terms because of the
equation (\ref{XXXX}).  Now notice that, for $\psi\in{\cal D}$,
we have that
\begin{align*}
X_0 (\psi\lfloor\vn)_I
& = \bigl( \nabla_{\vn} (\psi\lfloor\vn) \bigr) +
\{ 0\ \text{order terms}\} \\
& = \bigl( (\nabla_{\vn} \psi)\lfloor \vn\bigr)_I 
+\{ 0\ \text{order terms}\} \\
& = \{ 0\ \text{order terms}\} .
\end{align*}
Thus
\begin{multline*}
\sum_{kJ} \int_{b\Omega} (Y_k \psi_J)
\overline{ Y_k
\bigl( \sum_{j\neq0,I} \e{J}{jI} X_j 
(\psi\lfloor\vn)_I \bigr)} \\
=  \sum_{kJ} \int_{b\Omega} (Y_k \psi_J)
\overline{ \bigl( \nabla_{Y_k}
( d(\psi\lfloor\vn)) \bigr)_J} 
+ \{ \text{l.o.t.} \} .
\end{multline*}

Now we notice that, using equation
(\ref{XXXX}), it follows that for
$\psi\in{\cal D}$, the coefficients of $d(\psi\lfloor\vn)$ are all 
coefficients of $d\psi$, modulo lower 
order terms in the components of
$\psi$.  Notice also that the lower order terms that we produced in
the previous calculation are all 
${\cal O} (\|\psi\|_{W^1 (b\Omega)})$. 

Therefore
\begin{align*}
|\l d' \psi,\K\psi\r_1 |
& \le C \bigl(
\sum_{kJ} \| \nabla_{Y_k} \psi_J \|_{W^{1/2} (b\Omega)}
\| \nabla_{Y_k} d(\psi\lfloor\vn)\|_{W^{-1/2}(b\Omega)} 
+ \|\psi\|_{W^1 (b\Omega)}^2 \bigr)\\
& \le C \bigl( \sum_J \|\psi_J \|_{W^{3/2} (b\Omega)} 
\| d(\psi\lfloor\vn)_J \|_{W^{1/2}(b\Omega)} 
+ \|\psi\|_{W^1 (b\Omega)}^2 \bigr)\\
& \le C \bigl(
\|\psi\|_{W^{3/2}(b\Omega)} \|d\psi\|_{W^{1/2}(b\Omega)} 
+ \|\psi\|_{W^1 (b\Omega)}^2 \bigr)\\
& \le C ( \|\psi\|_2 \| d\psi\|_1 +\|\psi\|_{3/2}^2 ) \\
& \le \ep \|\psi\|_2^2 +C_\ep (\| d\psi\|_1^2 +\|\psi\|_1^2 ).
\end{align*}
This proves Claim 2.
\mbox{}\medskip\\
{\em Proof of Claim 1.\/}
Fix an open cover $\{ U_\ell \}_0^M$ of $\overline{\Omega}$ such
that $\overline{U_0} \ss \Omega$ and, for $\ell=1,\dots,M$,
on each $U_\ell$ we can find an orthonormal frame $\omega_0 =dx_0,
\omega_1, \dots, \omega_N$ for the space of $1$-forms.  Let $X_0=
\p/\p x_0, X_1, \dots, X_N$ be the dual frame of vector fields.
With 
respect to this basis the Laplacian is not diagonal on the space of 
$q$-forms, but it is diagonal in the top order terms.

On $U_\ell$ we write
$\psi=\sum_I \psi_I \omega^I$.  Let $\{ \eta_\ell^2 \}$ 
be a partition of unity subordinate to the cover $\{U_\ell \}$.
In this proof we denote by ${\cal E}$ any quantity which is 
${\cal O}(\| \psi\|_2 \| \psi\|_1 )$.  
Let $A$ for the moment denote either operator $d$ or $d'$.  Then
notice that
\begin{align*}
\l A\psi, A\psi\r_1
& = \sum_{jI} \l D_j (A\psi)_I, D_j (A\psi)_I \r_0 
			+ {\cal E} \\
& = \sum_{jI} \sum_\ell \int_\Omega \eta_\ell^2 
	(D_j A\psi)_I \overline{(D_j A\psi)_I} + {\cal E} \\
& = \sum_{jI} \sum_\ell \int_\Omega 
D_j \bigl( A(\eta_\ell \psi) \bigr)_I 
\overline{D_j \bigl( A(\eta_\ell \psi) \bigr)_I} + {\cal E} .
\end{align*}
Then we see that
it  suffices to consider forms with support in one of the
patches $U_\ell$.  When integrating by parts, this reduction
only produces error terms of the type ${\cal E}$ already considered.

Notice that the form $\eta_\ell \psi$ may not belong to ${\cal D}$.
Nonetheless $\eta_\ell \psi$ satisfies boundary equations of type
(\ref{XXXX}) and (\ref{YYYY}), and this is all we need.

Thus, without loss of generality, 
let $\psi\in{\cal D}$ have support in a small
open set on which we can write $\psi=\sum_I \psi_I
\omega^I$.  
We have
$$
d\psi=\sum_J \bigl( \sum_{j I} \e{J}{jI}  X_j \psi_I 
\bigr)\omega^J 
+ \{ 0\text{ order terms}\},
$$
and 
$$
d' \psi = \sum_L \bigl( \sum_{\ell I} \e{I}{\ell L}X'_\ell \psi_I
\bigr) \omega^L
+ \{ 0\text{ order terms}\},
$$
where $X'_\ell$ is the formal adjoint of $X_\ell$.  

Now notice that there exists a constant 
$C_0$ such that for $\theta\in W^1_q$, 
$\theta=\sum_I \theta_I \omega^I$, we have that
\begin{equation}\label{NORM}
\| \theta\|_1^2 \ge C_0 \sum_I \bigl( \sum_j \|X_j \theta_I \|_0^2 
+ \|\theta_I \|_0^2 \bigr) ,
\end{equation}
where the constant $C_0$ depends only on the choice of
$\omega_0,\dots,\omega_N$. 

Then,    
\begin{align}
\lefteqn{\sum_{jI} 
\bigl( \|X_j (d\psi)_I \|_0^2 +\|X_j (d'\psi)_I \|_0^2 \bigr)}
			\notag \\
& = \sum_J \sum_{pIqI'} \sum_{k=0}^{N} \int_\Omega 
	(\e{J}{pI}  X_p X_k \psi_I )  
	\overline{(\e{J}{q I'} X_{q} X_k \psi_{I'} )} \notag \\
& \qquad  +\sum_L \sum_{p Iq I'} 
	\sum_{k=0}^{N} \int_\Omega 
	(\e{I}{p L} X'_p X_k \psi_I) 
	\overline{(\e{I'}{q L} X'_{q} X_k \psi_{I'})}
			+{\cal E}  \notag \\ 
& = \sum_{kII'}\biggl[ \sum_{pq} \int_\Omega 
\bigl( \sum_J \e{J}{pI}\e{J}{qI'} \bigr) 
(X_p X_k \psi_I ) \overline{(X_q X_k \psi_{I'})} \notag \\
& \qquad 
+ \sum_{pq} \int_\Omega 
\bigl( \sum_L \e{I}{pL}\e{I'}{qL} \bigr) 
(X'_p X_k \psi_I ) \overline{(X'_q X_k \psi_{I'})}
		\biggr] +{\cal E} \notag \\ 
& = \sum_{kII'}\biggl[ \sum_{pq} \int_\Omega 
\bigl( \sum_J \e{J}{pI}\e{J}{qI'} \bigr) 
(X_p X_k \psi_I ) \overline{(X_q X_k \psi_{I'})} \notag \\
& \qquad 
+ \sum_{pq} \int_\Omega 
\bigl( \sum_L \e{I}{pL}\e{I'}{qL} \bigr) 
(X_p X_k \psi_I ) \overline{(X_q X_k \psi_{I'})}
		\biggr] +{\cal E} ,  \label{sum-of-d-d'}
\end{align}
where the last equality holds since
$X'_\ell =-X_\ell +L_0$, where $L_0$ is an
operator of order 0.  

Consider the term
\begin{multline}\label{---}
\int_\Omega  \e{J}{pI}\e{J}{qI'} 
(X_q X_k \psi_I ) \overline{(X_p X_k \psi_{I'})} 
+  \int_\Omega \e{I}{pL}\e{I'}{qL} 
(X_p X_k \psi_I ) \overline{(X_q X_k \psi_{I'})} \\
 \equiv I +I\!I .
\end{multline}
We need to distinguish two cases, according to whether $p=q$ or 
$p\neq q$.  Notice that
if $p=q$ in (\ref{---}), then
\begin{equation}\label{p=q}
I+I\!I = 2\| X_p X_k \psi_I \|_0^2 +{\cal E},
\end{equation}
since $\e{J}{qI}\e{J}{pI'}=\e{I}{pL}\e{I'}{qL} =0$ 
if $I\neq I'$. 

Suppose now that $p\neq q$.  We show in this case that
the term in (\ref{---}) is 
 of type ${\cal E}$.  
Notice that if $p\neq q$ then
\begin{equation}\label{EPSILON}
\e{J}{qI}\e{J}{pI}+\e{I}{pL}\e{I}{qL} =0,
\end{equation}
because, if $p<q$, then $\e{J}{qI}=-\e{I'}{qL}$ and
$\e{J}{pI'}=\e{I}{pL}$.  Then, if $p,q\neq0$, integration by parts
gives rise to no boundary term, and (\ref{EPSILON}) shows that
\begin{equation}\label{p,q-not0}
I+I\!I = {\cal E}
\end{equation}
when $p\neq q$, $p,q\neq0$.

If $p\neq0$ and $q=0$ we distinguish two cases according to
whether $k=0$ or $k\neq0$.

If $q=0$ and $k=0$, then we integrate by parts in $I\!I$.  
Recall that
$\psi\in\D$ implies that for $I'\ni0$ $\psi_{I'}$
satisfies the boundary equation (\ref{XXXX}).
Then we have
\begin{align*}
I\!I
& = -\int_\Omega \e{J}{0I}\e{J}{pI'}
(X_0 X_p X_0 \psi_I )\overline{X_0 \psi_{I'}}
+ \int_{b\Omega} \e{I}{pL}\e{I'}{0L}
(X_p X_0 \psi_I)\overline{X_0 \psi_{I'}} +{\cal E}\\
& = \int_\Omega \e{J}{0I}\e{J}{pI'}
(X_0 X_0 \psi_I) \overline{(X_p X_0 \psi_{I'})} 
+\int_{b\Omega} (X_p X_0 \psi_I) 
\overline{L_0 \psi_{I'}} +{\cal E} .
\end{align*}
Now using equation (\ref{EPSILON}) again, we see that 
\begin{align}
|I + I\!I|
& \le \big| \int_{b\Omega} 
(X_p X_0 \psi_I) L_0 \psi_{I'}\big| + {\cal E} \notag \\ 
& \le C \int |X_0 \psi_I X_p \psi_{I'} | +{\cal E}\notag \\
& \le C\|X_0 \psi_I \|_{W^{1/2}(b\Omega)}
			\|\psi_I\|_{W^{1/2}(b\Omega)}
+{\cal E} \notag \\
& = {\cal E} ,    \label{q=0,k=0}
\end{align}
that is
$I+I\!I={\cal E}$ in this case. 

Finally, suppose that $q=0$ and $k\neq0$. 
Notice that $\e{J}{0I}=\e{I'}{0L}=1$, and $\e{J}{pI'}=-\e{I}{pL}$
since $I'\ni0$.  Then
\begin{align*}
I+I\!I
& = -\int_\Omega \e{I}{pL} (X_0 X_k \psi_I)
\overline{(X_p X_k \psi_{I'})} 
+\int_\Omega \e{I}{pL} (X_p X_k \psi_I)
\overline{(X_0 X_k \psi_{I'})} + {\cal E} \\
& = \int_\Omega \e{I}{pL} (X_p X_0 X_k \psi_I)
\overline{X_k \psi_{I'}} 
+\int_\Omega \e{I}{pL} (X_p X_k \psi_I)
\overline{(X_0 X_k \psi_{I'})} + {\cal E} \\
& = -\int_\Omega \e{I}{pL} (X_p X_k \psi_I)
\overline{(X_0 X_k \psi_{I'})} 
+\int_{b\Omega} \e{I}{pL} (X_p X_k \psi_I)
\overline{X_k \psi_{I'}} \\
& \qquad \qquad \qquad 
+\int_\Omega \e{I}{pL} (X_p X_k \psi_I)
\overline{(X_0 X_k \psi_{I'})} + {\cal E} \\
& = -\int_{b\Omega} X_k \psi_I 
\overline{X_k (X_p \psi_{I'}\e{I}{pL})} +{\cal E}. 
\end{align*}
Hence, using the fact that $k\neq0$ we have that
\begin{align}
|\sum_{pI'} I+I\!I|
& \le \big| \int_{b\Omega} X_k \psi_I 
\overline{ X_k \bigl( \sum_{pI'} X_p \psi_{I'}\e{I}{pI'} \bigr)}
\big|	+{\cal E} \notag \\
& = \big| \int_{b\Omega} X_k \psi_I 
\overline{ X_k ( d\psi\lfloor\vn)_I}\big| +{\cal E} \notag \\
& \le \| X_k \psi_I \|_{W^{1/2}(b\Omega)} 
\|X_k (d\psi\lfloor\vn)_I \|_{W^{-1/2}(b\Omega)} +{\cal E} \notag \\ 
& \le C \| \psi \|_{W^{3/2}(b\Omega)} 
\| d\psi\lfloor\vn \|_{W^{1/2}(b\Omega)} +{\cal E} \notag \\ 
& \le C \|\psi\|_2 \|d\psi\|_1 +{\cal E} \notag \\ 
& \le \ep \|\psi\|_2^2 + C_\ep \|d\psi\|_1^2 +{\cal E} .  
				\label{q=0,knot0} 
\end{align}
Therefore, collecting (\ref{p=q}) - (\ref{q=0,knot0}) and
substituting them into 
(\ref{sum-of-d-d'}), and recalling (\ref{NORM}), we obtain that
\begin{align*}
\frac1{C_0} (\| d\psi\|_1^2 + \|d' \psi\|_1^2  )
& \ge \sum_{kI} 
\bigl( \|X_k (d\psi)_I \|_0^2 +\|X_k (d'\psi)_I \|_0^2 \bigr) \\
& = \sum_{kII'} \bigl[ \sum_{pq} I +I\!I \bigr] \\
& \ge  \sum_{kpI} 2\| X_p X_k \psi_I \|_0^2 
	- (\ep \|\psi\|_2^2 + C_\ep \|d\psi\|_1^2 )
	+{\cal E} \\
& = (2-\ep) \|\psi\|_2^2 - C_\ep \|d\psi\|_1^2 	+{\cal E} \\
& \ge (2-2\ep)\|\psi\|_2^2 - C_\ep (\|d\psi\|_1^2 +\|\psi\|_1^2) .
\end{align*}
From this Claim 1 follows easily.
This proves the theorem. 
\endpf

\section{The Regularity Theorem in the Case of Functions}
\setcounter{equation}{0}
We now turn to the question of regularity. In this
section we are going to 
prove the regularity result for the case $q = 0$.

\begin{theorem}\label{regularity}  \sl
Let $f \in C^\infty(\overline{\Omega})$ be orthogonal in
the $W^1$ inner product to the space of constant functions.
Then there is a unique function 
$u \in C^\infty(\overline{\Omega})$
that solves the boundary value problem
\begin{equation*}
\begin{cases}
d^* d\, u  =  f \qquad & \on \Omega \\
du \in \dom d^*  & 
\end{cases} 
\end{equation*}
In other words, $u$ solves the system
\begin{equation*}
\begin{cases}
\displaystyle{
(- \btu + G_\Omega) u  =  f }\qquad & \on \Omega \\ 
\displaystyle{
\pd{^2 u}{n^2}    =  0} 
& \on  b\Omega 
\end{cases}  
\end{equation*}
Moreover, the solution $u$ satisfies the 
desired (coercive) estimates
for each $s > 1/2$; i.e.\ there exists a $c_s > 0$ such that
\begin{equation}\label{Dagger}
\|u\|_{s+2} \leq c_s \bigl ( \|f\|_s + \|u\|_0 \bigr ) .
\end{equation}
\end{theorem}
\bgpf  
By the existence theorem, Theorem \ref{EXISTENCE}, we know that the
solution to the above \bvp\ exists for all $f$ orthogonal to the
harmonic space ${\cal H}_0$.  In the forthcoming Theorem \ref{H_0} we 
shall prove
that ${\cal H}_0$ reduces to the constants.  Now we turn to the
question of estimates.

Let $u$ be a solution of the boundary value problem,
with $f \in C^\infty(\overline{\Omega})$.  
Let $s > 1/2$.  We wish
to estimate $\|u\|_{s+2}$ in terms of $\|f\|_s$.  It is clear
that it suffices to estimate $\|\eta u \|_{s+2}$ for a given
cut-off function $\eta$ with small support.

We first suppose that $\supp \eta \cap b\Omega = \emptyset$.
[The second step will be to assume that 
$\supp \eta \cap b\Omega \ne \emptyset$.]
\subsection*{The interior estimate}  
Let $\eta \in C_0^\infty (\Omega)$.   
Recalling that $f=-\btu u +G_\Omega u$ 
we have that 
\begin{align}
\| \eta u\|_{s+2}
& \le c \| \btu (\eta u)\|_s \notag \\
& \le c \bigl( \| \btu(\eta u) +\eta f\|_s +\| \eta f\|_s 
\bigr) \notag \\  
& \le c\bigl( \|\btu (\eta u) -\eta\btu u \|_s 
+ \| \eta G_\Omega u\|_s  + \|\eta f\|_s \bigr) \notag\\ 
& \le c \bigl( \| \eta_1 u\|_{s+1} +\| \eta G_\Omega u \|_s 
+\|\eta f\|_s \bigr) ,  \label{yet1more}
\end{align}
where $\eta_1 \equiv 1$ on $\supp \eta$.  
Next we want to show that 
$\|\eta G_\Omega u\|_s \le c \| u\|_{s+1}$.  Recall that 
by Proposition \ref{esp-K-Omega} 
and equation (\ref{LuigiXIV})
$G_\Omega$ is the solution of the \bvp
$$
\begin{cases}
\displaystyle{ -\btu v+v =0 }\qquad & \on \Omega \\
\displaystyle{ 
\pd{v}{n} = \sum_{k=0}^{N} \nabla_{Y^*_k}\bigl(
(\nabla_{Y_k} du)\lfloor\vn \bigr) 
+\pd{u}{n}  } & \on b\Omega
\end{cases} .
$$
Now it is easy to see that
(using the Fourier transform for instance)
$$
\| \eta G_\Omega u \|_s \le c\|\eta_1 G_\Omega u\|_{s-1}
\le c \| G_\Omega u \|_{s-1} .
$$
(Here, for $t>0$,
$\| \, \cdot \, \|_{-t}$ is the norm in the 
Sobolev space $W_{-t}(\Omega) \equiv (\stackrel{\circ}
{W^t}(\Omega))^*$.)
Moreover, by \cite{TRI} Theorem 4.2.4
we have
\begin{align}
\| G_\Omega u \|_{s-1}
& \le c\|T_2 \pd{u}{n} -T'_2 u\|_{W^{s-5/2}(b\Omega)} \notag\\
& \le c\bigl( \big\| \pd{u}{n} \big\|_{W^{s-1/2}(b\Omega)} 
+\| u\|_{W^{s-1/2}(b\Omega)} \bigr) \notag \\
& \le c\|u\|_{s+1} . \label{yet2more}
\end{align}
Thus (\ref{yet1more}) and (\ref{yet2more}) give that
$$
\| \eta u\|_{s+2} \le c ( \| \eta f \|_s + \| u\|_{s+1} ).
$$
\subsection*{The boundary estimate} 
Fix a Fermi coordinate patch $$\bigl( U,\Psi=(x_0,\dots,x_N)\bigr),$$
and let $\eta\in C^\infty_0 (U)$. 

In what follows we denote by $G_\Omega$ both the operator defined 
on $\Omega$ and (when restricted to $U$) the same operator
expressed in the local chart and thus defined on $\Psi(U)$.
This technical ambiguity should cause no confusion.  We will
denote by $G_{\RR^{N+1}}$ the operator arising from
considering the adjoint of $d^*$ in the half space.
Finally, given a function $u \in C^\infty(U)$ 
we denote by
$\tilde{u}$ the function 
$$u \circ \Psi^{-1} \in C^\infty(\Psi(U)).$$ 

Now consider $\eta, \eta_1 \in C_0^\infty(U)$ with $\eta_1 = 1$
on the support of $\eta$.  These functions 
are chosen to be constant 
along the normal direction to $b\Omega$ near the boundary.
In our
local coordinates, $\eta u$ satisfies the following boundary
condition: 
$$
\frac{\partial^2}{\partial {x_0}^2} (\widetilde{\eta u}) 
=0 \qquad \qquad \on \{ x_0 =0\} .
$$
Recall that if $v$ satisfies the \bvp 
\begin{equation*}
\begin{cases}
\displaystyle{
- \btu v+ G_{\RR_+^{N+1}} v=f }& \qquad \on \RR_+^{N+1} \\  
\displaystyle{
\pd{^2 v}{x_0^2}   =  0 } &\qquad \on b\RR_+^{N+1} 
\end{cases} 
\end{equation*}
and if $v$ has compact support,
then
$$
\|v\|_{s+2} \leq C \bigl\{ \|f\|_s \ + \ \|v\|_{s+1} \bigr\} , 
\qquad \qquad  s > 1/2 .
$$

Then
\begin{equation}\label{SSTAR}
\|\widetilde{\eta u}\|_{s+2} \leq
C\bigl\{ \big\|(-\btu +G_{\rnp})(\widetilde{\eta u}) \big\|_s 
      + \|\widetilde{\eta u}\|_{s+1} \bigr\}. 
\end{equation}
Our goal is to replace $- \btu +\ G_{\rnp}$ with 
$-\btu + \ep \L_2  + G_\Omega$ modulo error terms that are 
controlled by lower order norms of $\eta u$.
Thus the estimate (\ref{SSTAR}) gives that
\begin{eqnarray}
\|\eta u\|_{s+2} 
& \leq & C \bigl\{ \|(-\btu + \ep \L_2  
	+ G_\Omega)(\widetilde{\eta u}) \|_s
  	+ \|\widetilde{\eta u}\|_{s+1} \nonumber \\
&      & \qquad + \ep \|\L_2  (\widetilde{\eta u}) \|_s 
+ \bigl\| G_{\rnp}(\widetilde{\eta u}) 
		\bigl[ \eta G_\Omega u\bigr]\til \bigr\|_s  
 	 \bigr\} \nonumber \\
& \leq & C \biggl\{ \bigl \| \tilde{\eta} 
	[(-\btu+\ep\L_2  +G_\Omega)(u)]\til \bigr\|_s 
+ \ep \|\eta u \|_{s+2} + \|\eta_1 u \|_{s+1} \nonumber \\ 
&      & \qquad + \bigl\| G_{\rnp} (\widetilde{\eta u})
  - [\eta G_\Omega u]\til \bigr\|_s
\biggr\} \nonumber \\
&   \le & C
\bigl\{ \big\| -\tilde\eta \btu \tilde u +G_{\RR^{N+1}} 
(\widetilde{\eta u}) \big\|_s + \|\eta_1 u\|_{s+1} \bigr\} 
\nonumber \\
& \le & C \bigl\{ \big\| \tilde\eta  (-\btu +\ep\L_2  +
G_\Omega) (\tilde u) \big\|_s +\ep \|\eta u\|_{s+2} \nonumber \\
&  & \qquad
+ \big\| G_{\RR^{N+1}} (\widetilde{\eta u}) 
-\bigl[ \eta G_\Omega u \bigr]\til  \big\|_s 
+ \|\eta_1 u \|_{s+1} \bigr\} \label{(1)}  \\
& =  &  C \bigl\{ \| \eta f \|_s +\ep \|\eta u\|_{s+2} 
+ \big\| G_{\RR^{N+1}} (\widetilde{\eta u}) 
-\bigl[ \eta G_\Omega u \bigr]\til  \big\|_s 
+ \|\eta_1 u \|_{s+1} \bigr\} .  \nonumber 
\end{eqnarray}      
The function $G_{\rnp}(\widetilde{\eta u})$ solves the
\bvp 
$$
\begin{cases}
\displaystyle{
(-\btu + I)\bigl( G_{\rnp}(\widetilde{\eta u}) \bigr)  =  0} 
			& \on \rnp \\
\displaystyle{
\pd{}{x_0} \bigl( G_{\rnp} (\widetilde{\eta u})\bigr) =  
\pd{}{x_0}(\widetilde{\eta u}) -\btu'_{\RR^N} 
\pd{}{x_0} (\widetilde{\eta u}) }
			& \on b\rnp 
\end{cases} .
$$
On the other hand, by (\ref{LuigiXIV}),
the function $[\eta G_\Omega u]\til$, 
solves the \bvp
$$
\begin{cases}
\displaystyle{
\bigl( -\btu+I +\epsilon \L_2  \bigr) 
([\eta G_\Omega u]\,\tilde{\text{}}\,)
 =  - [2 \nabla \eta \cdot \nabla G_\Omega u 
+\btu \eta \cdot G_\Omega u ]\,\tilde{\text{}}\,} & \on \rnp \\
\displaystyle{
\pd{}{x_0} \bigl( [\eta G_\Omega u]\til \bigr) 
 =  T_1 \pd{}{x_0} (\widetilde{\eta u}) -(\btu'_{\RR^{N+1}} 
- \ep \L_2 ) \bigl( \pd{}{x_0} (\widetilde{\eta u})  \bigr) 
+  T_2 (\widetilde{\eta u}) } & \on b\rnp 
\end{cases} .
$$
Let 
$w \equiv G_{\rnp}(\widetilde{\eta u})-[\eta G_\Omega u ]\til$. 
Then $w$ solves the system 
\begin{equation}\label{SSSTAR}
\begin{cases}
\displaystyle{
(\btu - I)w  =  \ep \L_2 \bigl( [\eta G_\Omega u ]\til
\bigr) + [2 \nabla \eta \cdot \nabla G_\Omega u 
      + (\btu \eta) G_\Omega u ]\til }
& \on \rnp \\ 
\displaystyle{
\pd{w}{x_0}   =  \tilde{\eta} \ep {\cal T}_2 
\pd{\tilde u}{x_0}
   + \tilde\eta T_1 \pd{\tilde u}{x_0} 
  + \tilde \eta  T_2 \tilde u }
& \on b\rnp
\end{cases} 
\end{equation}
The right hand side of the first equation in the system 
(\ref{SSSTAR}) can be 
rewritten as 
$$
\ep \eta' L_2 \widetilde{G_\Omega u} + \eta'' L_1
\widetilde{G_\Omega u} 
$$
where the order of $L_j$ is $j$ and $\supp \eta', \supp \eta''$ 
lie in $\supp \tilde{\eta}.$

Now we apply Theorem 4.2.4 in \cite{TRI}.
We obtain that, for any $s>1/2$,
\begin{eqnarray}
 \|w\|_s 
& \leq & C \biggl\{ \ep \bigl\| 
L_2 ([\eta G_\Omega u]\til ) \bigr\|_{F^{s-2}} 
+ \bigl\| \eta' L_1 [G_\Omega u]\til \bigr\|_{F^{s-2}} 
\nonumber \\
&   & \quad + \ep \bigl \|\eta {\cal T}_2  \pd{u}{n} 
\bigr \|_{W^{s-3/2}(b\Omega)}  + \bigl \| \eta T_1  
 \pd{u}{n} \bigr \|_{W^{s - 3/2}(b\Omega)}
+ \bigl \|\eta T_2 u \bigr \|_{W^{s-3/2}(b\Omega)}
\biggr\}  \nonumber \\
& \leq & C\biggl\{ \ep \bigl \| \eta G_\Omega u \bigr \|_{F^s}
  + \bigl \|\eta' [G_\Omega u] \bigr \|_{F^{s-1}}
   + \ep \left \| \eta \frac{\partial u}{\partial n} 
        \right\|_{W^{s + 1/2}(b\Omega)} \nonumber \\
 &   & \qquad + \left \| \eta \pd{u}{n} 
          \right \|_{W^{s-1/2}(b\Omega)}
      + \bigl\|\eta_1 u \bigr\|_{W^{s+1/2}(b\Omega)} 
      		\biggr\}\nonumber \\
& \stackrel{G_{\Omega} u = \btu u + f}{\leq} & 
	C \biggl\{ \ep \bigl \| \eta \btu u \bigr \|_s 
	+ \ep \bigl\| \eta f \bigr\|_s  
	+ \bigl \|\eta_1 \btu u \bigr \|_{F^{s-1}} \nonumber  \\ 
&   & \qquad  + \bigl \| \eta_1 f \bigr \|_s 
 	+ \ep \bigl\| \eta u \bigr\|_{s+2} 
 	+ \bigl\|\eta_1 u \bigr\|_{s+1} \biggr\} \nonumber \\
& \leq & C \biggl\{ \ep \bigl\| \eta u \bigr\|_{s+2} 
	+ \bigl\| \eta_1 f \bigr\|_s
        + \bigl\| \eta_1 u \bigr\|_{s+1} \biggr\} . \label{(2)} 
\end{eqnarray} 
Now estimates (\ref{(1)}) and (\ref{(2)}) yield that
$$
\|\eta u\|_{s+2} \leq C \bigl\{ 
\|\eta_1 f \|_s + \ep \|\eta u\|_{s+2} + \|\eta_1 u\|_{s+1}  
\bigr\} ,
$$
that is, for $s>1/2$, 
$$
\|\eta u \|_{s+2} 
\leq C \bigl\{ \|\eta_1 f \|_s + \|\eta_1 u \|_{s+1} \bigr\} . 
$$
This concludes the proof.
\endpf

\section{Estimates for $q$-Forms}
\setcounter{equation}{0}

\subsection*{The Regularity Theorem}
In this section we prove the estimate for the solution of the \bvp\
in the case of $q$-forms.  We follow the same outline as in the case 
of functions.  For higher degree forms some extra technicalities are
needed.  The differential
operators acting on forms that are involved are not (usually)
diagonal, and we need some extra care for the off-diagonal terms.

\begin{theorem}\label{REG-THM-Q}  \sl
Let $s>1/2$.  Let
$\alpha\in W^{s}_q (\Omega)$, 
and $\alpha$ orthogonal to the
kernel of the \bvp.  Let $\psi\in\bw^q$ be a solution of 
\begin{equation*}
\begin{cases}
\displaystyle{(dd^* +d^* d) \psi =\alpha} & \on \Omega \\
\psi\in \dom d^* & \\
d\psi\in \dom d^* &
\end{cases} 
\end{equation*}
There exists $c_s >0$, independent of $\psi$, such that
$$
\| \psi\|_{s+2} \le c_s (\|\alpha\|_s + \| \psi\|_{s+1} ).
$$
\end{theorem}
Thus we wish to estimate $\| \psi\|_{s+2}$ in terms of
$\|\alpha\|_{s}$,  
for $s>1/2$.  It suffices to estimates $\eta\psi$ for a cut-off
function $\eta$ with support contained in a small open set.  
If $\eta\in C^\infty_0 (\Omega)$ then the argument used in the case of
functions (see the proof of \ref{regularity}) 
applies with no substantial change. 

Thus we
assume that $\eta$ is a cut-off function whose support is contained 
in an open set on which there exists a Fermi coordinate chart.  
As before, we also suppose that 
$(\p/\p n)\eta=0$ in a 
tubular neighborhood of the boundary. 

As in the case of $q$-forms in the half space we write
\begin{align*}
(dd^* +d^* d) & = (dd' +d'd)+(d\K +\K d) \\
& \equiv  -\btu+G_\Omega .
\end{align*}
We now proceed as in the function case.  We need to estimate 
$\|\eta\psi\|_{s+2}$. 
We introduce Fermi local coordinates, and set the stage as in 
Section
\ref{SPECIAL}.
Recall that such coordinates have the property that the
normal direction is orthogonal to the remaining directions at all
points in the coordinate patch (that is, such coordinates ``flatten" 
the boundary.).  Given a form $\phi$ on $\Omega$
(i.e.\ written in global coordinates), we write $\tilde\phi$ to
indicate the form written in the Fermi coordinates.  We adopt the 
same notation for the operators.

We wish to apply the estimates proved in the half space for the
$q$-forms.  In order to do this we notice that 
$\widetilde{(\eta\psi)}$  is a form defined on 
the half space, and we wish to 
check if it satisfies a \bvp\ for
which we have favorable estimates.  Since $\psi$ and $d\psi$ belong 
to $\dom d^*$, they  satisfy equations (\ref{XXXX}) and (\ref{YYYY}).
 We need to examine those equations more closely.  In order to do
this we rewrite equations (\ref{XXXX}) and (\ref{YYYY}) in our local
coordinates.

\begin{proposition}   \sl
Let $\psi\in\bw^q (\overline{\Omega})$ be such that both 
$\psi,\, d\psi \in\dom d^*$. Let $\eta$ be a cut-off function as
before. Then the form $\widetilde{\eta\psi}$ satisfies the boundary  
equations
\begin{equation}\label{bdry-cond-0inI}
\pd{}{x_0} (\widetilde{\eta\psi})_I + \sum_{|J|=q \atop J\ni 0} 
\gamma_{IJ} (\widetilde{\eta\psi})_J =0 \qquad \text{if } I\ni 0, 
\end{equation}  
 and, 
\begin{equation}\label{bdry-cond-0not-inI}
\pd{^2}{x^2_0} (\widetilde{\eta\psi})_I + \sum_{|K|=q \atop K\ni 0}  
T_1 \bigl(  (\widetilde{\eta\psi})_K \bigr)
+\sum_{|L'|=q \atop L'\not\ni 0 } a_{K'L'} \pd{}{x_0}
(\widetilde{\eta\psi})_{L'} + {\cal E}
=0  \quad \text{if } I\not\ni 0 ,
\end{equation}
where ${\cal E}$ are $\{ 0\text{ order terms}\}$ in the 
components of
$\psi$.
\end{proposition}
We stress that we have improved equations (\ref{XXXX})
and (\ref{YYYY}) by separating the ``normal" and ``tangential''
components of the lower order terms in the boundary equations. 
\bgpf
By Lemma \ref{covariant-diff}, 
the condition $\psi\in\dom d^*$, i.e.\ $\nabla_{\vn}\psi\lfloor\vn=0$  
on $b\Omega$, becomes
$$
\pd{}{x_0} (\widetilde{\eta\psi})_I +\sum_{|J|=q\atop J\ni 0} 
\gamma_{IJ} (\widetilde{\eta\psi})_J =0 \qquad \on \RR^N
\text{\ if } I\ni0. 
$$
On the other hand, the second boundary condition gives that
$$
 \pd{}{x_0} (\widetilde{d\psi})_K +\sum_{|L|=q+1\atop L\ni 0}  
\gamma_{KL} (\widetilde{d\psi})_L =0 \qquad \on \RR^N
\text{\ if } K\ni0,
$$
i.e., $\on \RR^N$ for $K\ni 0$,
\begin{multline*}
\pd{}{x_0} \biggl( \pd{\tilde\psi_{K'}}{x_0} \e{K}{0K'} 
+ \sum_{K''\atop j=1,\dots,N} \pd{\tilde\psi_{K''}}{x_j} 
\e{K}{jK''} \biggr) +\pd{}{x_0}\bigl( \sum_J \psi_J d\omega^J
\bigr)_J \\
+ \sum_{|L|=q+1\atop L\ni0} \gamma_{KL} \bigl(
\pd{\tilde\psi_{L'}}{x_0} 
\e{L}{0L'} +\sum_{L''\atop j=1,\dots,N} 
 X_j \widetilde{\psi_{L''}}
\e{L}{jL''}  \bigr) +\{ 0\text{ order terms}\} =0 . 
\end{multline*}
Now we use the first boundary condition  to replace 
the terms of the form
$(\p/\p x_0) 
\tilde\psi_{K''}$, with $K'' \ni 0$, with $0$ order terms 
in $\psi$. Hence, the term 
$(\p/x_0)\bigl( \sum_J \psi_J d\omega^J \bigr)_J $ only contributes
normal derivatives of tangential components and $0$ order
terms. 

Therefore, for $K'\not\ni 0$, we have
$$
\pd{^2}{x^2_0} (\widetilde{\eta\psi})_{K'}
+ \sum_{K''} T_1
(\widetilde{\eta\psi})_{K''} \e{0K'}{jK''} + (0\text{ order terms})
+\sum_{|L'|=q\atop L' \not\ni0}  a_{K'L'} 
\pd{}{x_0} (\widetilde{\eta\psi})_{L'} =0.
$$
From this identity, equation (\ref{bdry-cond-0not-inI}) follows. 
\endpf
At this point we adapt the estimates in the half space 
for the \bvp s 
\begin{equation}\label{Pb1}
\begin{cases}  
\displaystyle{(-\btu+G_{\rnp})\phi=\alpha } & \on\rnp\\
	\displaystyle{ \pd{^2 \phi_I}{x_0} =0 }& 
			\on\RR^N \quad\text{if }I\not\ni 0 
\end{cases}\ ,
\end{equation}
and
\begin{equation}\label{Pb2}
\begin{cases}  
\displaystyle{(-\btu+G_{\rnp})\phi=\alpha } & \on\rnp\\
\displaystyle{ \pd{\phi_I}{x_0} = 0 }
&\on\RR^N \quad\text{if }I\ni 0 
\end{cases} \ .
\end{equation}
to the form $\widetilde{\eta\psi}$.

We wish to apply Theorems \ref{apr-est} and \ref{apriori-est-q-hs}
to the problems (\ref{Pb1}) and (\ref{Pb2}) respectively.  In both
cases the boundary data is estimated by $\|
\widetilde{\eta\psi}\|_{s+1}$.
Let $\eta_1 \in C^\infty_0 (\overline{\Omega})$,
$\eta_1 \equiv 1 \on \supp \eta$. We obtain that  
\begin{align}\label{I-II-above}
\| \widetilde{\eta\psi} \|_{s+2}
& \le c\left( 
	\| \eta_1 (-\btu+G_{\rnp}) (\widetilde{\eta\psi})\|_s
	+\| \widetilde{\eta\psi}\|_{s+1} \right) \notag  \\
& \le c\left( \big\| 
\bigl[ \eta(-\btu+G_\Omega)\psi\bigr]\til \big\|_s 
	+ \big\| \bigl[ -\eta\btu\psi \bigr]\til 
	+ \btu(\widetilde{\eta\psi}) \big\|_s \right.\notag  \\
& \qquad\qquad \left. + \big\| \eta_1 G_{\rnp}(\widetilde{\eta\psi})
	-\eta_1 \bigl[ \eta G_\Omega \psi\bigr]\til \big\|_s 
	+ \| \widetilde{\eta\psi} \|_{s+1} \right) \notag  \\
& \le c\left( \| \widetilde{\eta\alpha} \|_s 
+ \big\| G_{\rnp}(\widetilde{\eta\psi}) 
	-\bigl[ \eta G_\Omega \psi\bigr]\til \big\|_s 
	+\| \widetilde{\eta\psi}\|_{s+1} \right) .
\end{align}  
Now we need to estimate $\big\| \eta_1 G_{\rnp}(\widetilde{\eta\psi}) 
	-\eta_1 \bigl[ \eta G_\Omega \psi\bigr]\til \big\|_s $ 
in the equation
(\ref{I-II-above}) above.  
We have the following theorem.
\begin{theorem}\label{marco1}  \sl
Let $s>1/2$. Then, for any $\epsilon>0$, 
there exists a $C_\epsilon >0$ 
such that 
$$
\big\| \eta_1 G_{\rnp}(\widetilde{\eta\psi}) 
	-\eta_1 \bigl[ \eta G_\Omega \psi\bigr]\til \big\|_s 
 \le \epsilon \| \widetilde{\eta\psi} \|_{s+2} 
	+C_\epsilon \| \eta_1 \psi \|_s ,
$$
where $\eta_1 \in C^\infty_0$ and $\eta_1 \equiv 1$ on 
$\text{supp}\, \eta$. 
\end{theorem}

Assuming the Theorem for now, we finish the proof of the
estimate for 
$\| \widetilde{\eta\psi}\|_{s+2}$,
and therefore the proof of Theorem \ref{REG-THM-Q}.  
Using Theorem \ref{marco1} and 
(\ref{I-II-above}) above we see that
\begin{equation*}
\| \widetilde{\eta\psi}\|_{s+2}
\le c \bigl( \| \eta\alpha\|_s + \ep \| \eta_1 \psi \|_{s+2}
\| \eta_1 \psi\|_{s+1} \bigr) ,
\end{equation*}
from which we obtain that
$$
\| \widetilde{\eta\psi}\|_{s+2}
\le c \bigl( \| \eta\alpha\|_s + \| \eta_1 \psi\|_{s+1} \bigr) ,
$$
which is what we wished to prove. \qed

\subsection*{Proof of Theorem 11.3}
Now we turn to the proof of Theorem \ref{marco1}.
The proof will be broken up into several lemmas.

The
form $G_{\rnp}(\widetilde{\eta\psi})\equiv\theta
=\sum_I \theta_I \omega^I$ is such that its components $\theta_I$ 
satisfy the following \bvp s
\begin{equation}\label{G0}
\begin{cases}
\displaystyle{(-\btu+I)\theta_I =0} & \on\rnp\\
\displaystyle{ \pd{\theta_I }{x_0} =
(-\btu'+I)(d(\widetilde{\eta\psi})_{0I})  }
	&\on\RR^N \quad\text{if }I\not\ni 0 
\end{cases} 
\end{equation}
and
\begin{equation}\label{Gnot0}
\begin{cases}
\displaystyle{(-\btu+I)\theta_I =0} & \on\rnp\\
\displaystyle{ \theta_I = (-\btu'+I)(\widetilde{\eta\psi})_I } 
		 & \on\RR^N \quad\text{if }I\ni 0
\end{cases} 
\end{equation}
On the other hand we need to investigate what kind of \bvp\ the form 
 $[\eta G_\Omega \psi]\til$ satisfies.  Recall that $G_\Omega \equiv
(\K d+d\K)$, and that the operator $\K$ 
is the operator on forms that
gives the solution of the following \bvp\ for $\theta\in\bw^q$.  
Namely we have 
$$  
\begin{cases}
\displaystyle{ (-\btu+I)\theta =0} & \on\Omega \\
\displaystyle{ \nabla_{\vn}\theta = \sum_{i=0}^{N}
\nabla_{Y_i^*}(\nabla_{Y_i} \phi\lfloor \vn ) +\phi\lfloor \vn} 
			& \on b\Omega 
\end{cases}
$$
for $\phi\in\bw^{q+1}(\overline{\Omega})$.
\begin{lemma} \label{eta-K-Omega-phi}   \sl
Let $\phi\in\bw^q (\overline{\Omega})$.  The form 
$[\eta\K \phi]\til \equiv \theta$
is such that its $I$-component satisfies the following \bvp\ on
$\rnp$ 
$$ 
\begin{cases}
\displaystyle{ (-\btu-\ep{\cal L}_2 +I) \theta_I
= \bigl( \eta_1 L_1 \widetilde{\K \phi}\bigr)_I } & \on \rnp \\
\displaystyle{ \pd{ \theta_I}{x_0} = -\sum_{|K|=q-1} \tilde \eta  
\gamma_{IK} ( \widetilde{\K \phi})_K +(E_2 \phi)_I } & \on\RR^N 
\end{cases}
$$
where
$$
(E_2 \phi)_I = 
\begin{cases}
-\tilde \btu_T \phi_{0I} +\phi_{0I} +(T_1 \phi)_I 
	& \text{if } I\not\ni 0 \\
0 & \text{if }I\ni 0 
\end{cases}  
$$
${\cal L}_2$ is defined by equation (\ref{Delta-tilde}).  
Moreover, $L_1$ is a first order differential operator,
$T_1 $ is a first order differential operator, both sending
$(q+1)$-forms into $q$-forms.
\end{lemma}
\bgpf
It suffices to use Lemma \ref{covariant-diff} and the computation
leading to formula (\ref{LuigiXIV}).
\endpf
\begin{lemma}\label{marco2}   \sl
Let $\phi\in\bw^q (\overline{\Omega})$, and let $I\ni0$.  
Then for any $s>1/2$,
$$
\| [\eta\K \phi]_I \til \|_s \le c \|\phi \|_s .
$$
\end{lemma}
\bgpf
Recall (see Corollary \ref{K-ord-1}) 
that $\K$ is an operator of order 1.
Therefore we direct the reader's attention to the fact that 
 the content of this lemma is that, on the normal components 
of the form, the operator $\K$ is of order 0. (Recall that, in the
case of the half space, the operator $\K$ was identically zero on
normal components.)

It suffices to recall the estimates in the negative norms for
elliptic \bvp s, as in \cite{TRI}, Theorem 4.2.4.
Using the spaces $F^s_{2,2} (\Omega)$, which we will simply denote by
$F^s$, we have
\begin{equation*}
\| [\eta\K \phi]_I \til \|_s
\le c\bigl( \| \bigl( \eta_1 L_1 \widetilde{\K \phi} \bigr)_I
						\|_{F^{s-2}} 
+  \sum_{K\ni0} \| [\eta\K \phi]_K \til \|_{W^{s-3/2}(b\Omega)}
\bigr) .   
\end{equation*}
Let $\epsilon >0$.
By selecting $\eta$ with support suitably small, we can achieve 
$$
\| [\eta\K \phi]_K \til \|_{W^{s-3/2}(b\Omega)} \le \epsilon 
\| [\eta\K \phi]_K \til \|_{W^{s-1/2}(b\Omega)} \, .
$$
Therefore
\begin{align*}
\sum_{I\ni0} \| [\eta\K \phi]_I \til \|_s
& \le c\|\bigl( \eta_1 L_1 \widetilde{\K\phi} \bigr)_I \|_{F^{s-2}} 
+\ep \sum_K \|[\eta\K \phi]_K \til \|_{W^{s-1/2}(b\Omega)} \\
& \le c\|\bigl( \eta_1 L_1 \widetilde{\K\phi} \bigr)_I \|_{F^{s-2}}  
+ \ep\sum_K \|[\eta\K \phi]_K \til \|_{F^{s}} .
\end{align*}
Then, by absorbing the last term on the right hand side over on the left, 
using the continuity of differential operators on the spaces $F^s$,  
and invoking Corollary \ref{K-ord-1}, we
obtain that
\begin{align*}
\sum_{I\ni0} \| [\eta\K \phi]_I \til \|_s
& \le c\| [\eta\K \phi] \til \|_{F^{s-1}(\Omega)} \\
& \le c \| \eta_1 \phi \|_s  .
\end{align*}
This concludes the proof of the lemma.
\endpf 
\begin{lemma} \label{eta-G-Omega-psi-0inI}   \sl
Let $\psi\in\bw^q (\overline{\Omega})$.  
Then, for $I\ni0$, $[\eta G_\Omega \psi]_I \til$ is  
a solution of the following \bvp 
$$  
\begin{cases}
\displaystyle{ (-\btu-\ep{\cal L}_2 +I) w = \bigl[ \eta_1 L_1
(G_\Omega \psi) \bigr]_I \til }
				& \on\rnp \\
\displaystyle{ w=\tilde\eta (-\btu_T +I)
\widetilde{\psi_I} + \zeta}   & \on \RR^N 
\end{cases} 
$$
where ${\cal L}_2$ is defined by equation (\ref{Delta-tilde}), 
$\zeta$ is a function satisfying the estimate
$$
\| \zeta\|_{W^{s-1/2}(b\Omega)} \le c \| \eta_1 \psi\|_{s+1},
$$
$L_1$ is a differential operator of order 1, and
$\eta_1 \in C^\infty_0 $ with $\eta_1 \equiv 1$ on
$\text{supp}\, \eta$. 
\end{lemma}
\bgpf 
The equation on $\rnp$ is easily seen to be satisfied.  Concerning  
the boundary equation, recall that $G_\Omega =\K d+d\K$.  Then, for
$I\ni0$, $I=0I'$, 
the function $\bigl[ \eta\K d\psi\bigr]_I \til |_{\RR^N}$ is the 
restriction to  $b\rnp$ of the solution of the \bvp\ in Lemma
\ref{eta-K-Omega-phi} with $\phi=d\psi$.  Moreover, for $I=0I'$,
on $\RR^N$ we have
\begin{align*}  
\bigl[ \eta d\K \psi\bigr]_I \til 
& = \pd{}{x_0} \bigl( \eta\K \psi\bigr)_{I'} \til 
	+\sum_{i=1,\dots,N \atop |J|=q-1} X_i
	\bigl[ \eta \K \psi\bigr]_J \til \e{0I'}{jJ} 
			+ (0\text{ order terms})\\  
& = -\sum_{K\not\ni0} \tilde\eta\gamma_{I'K} 
	(\widetilde{\K \psi})_K +(-\btu_T +I)\psi_{0I'}
	+\sum_{J\ni0} T_1 \bigl[ \eta\K \psi\bigr]_J \til
+ (0\text{ o. t.'s}) \, ,
\end{align*}
where  $T_1$ is a tangential differential operator of order 1. 
Therefore $w$ satisfies the boundary equation
\begin{align*}
w & = \tilde\eta(-\btu_T  +I)\psi_I 
-\sum_{K\not\ni0} \tilde\eta\gamma_{I'K} (\widetilde{\K \psi})_K
	+ \bigl[ \eta\K d\psi \bigr]_I \til 
+\sum_{J\ni0} T_1 \bigl[ \eta\K \psi\bigr]_J \til \\  
& \equiv  \tilde\eta(-\tilde\btu' +I)\psi_I + \zeta .
\end{align*}
 Now the desired estimate follows from Lemma \ref{marco2}. 
Indeed, recalling that $I\ni0$, we see that
\begin{align*}
\| \zeta\|_{W^{s-1/2}(b\Omega)}
& \le c \left( \| \widetilde{\K \psi} \|_{W^{s-1/2}(b\Omega)} 
+ \| \bigl[ \eta\K d\psi \bigr]_I \til \|_{W^{s-1/2}(b\Omega)}  
\right. \\
& \qquad \qquad \left.
+\sum_{J\ni0} \| T_1 \bigl[ \eta\K \psi\bigr]_J \til 
			\|_{W^{s-1/2}(b\Omega)} \right) \\  
& \le c \left( \| \widetilde{\K \psi} \|_s  
+ \| \bigl[ \eta\K d\psi \bigr]_I \|_s
+\sum_{J\ni0} \| \bigl[ \eta\K \psi\bigr]_J \|_{s+1}
\right) \\ 
& \le c\|\eta_1 \psi\|_{s+1} . \qed
\end{align*}
\renewcommand{\qed}{}\endpf
Before estimating the term $[\eta G_\Omega \psi]_I \til$ for
$I\not\ni0$ we need an extra lemma.
\begin{lemma}\label{LAST?}  \sl
Let $w_1=\eta(\K d\psi)_K \til -\eta(d\K \psi)_K \til$,
with $K\not\ni0$.  Then $w_1$
solves 
$$
\begin{cases}
\displaystyle{
(-\btu+I)w_1 = L_1 (\K d\psi) +L_2 (\K \psi) } 
	\qquad & \on \Omega \\
\displaystyle{
\pd{w_1}{n} = T_2 \psi + L_0 (\K \psi)} & \on \RR^N
\end{cases} \ ,
$$
where the $L_j$'s 
are operators of order $j$, and $T_2$ is a tangential
differential operator of order $2$.  

Moreover, there exists $c>0$ independent of $\psi$ such that 
$$
\|w_1\|_s \le c\| \psi\|_{s+1}.
$$
\end{lemma}
\bgpf
Recall that $(\K \phi)_K$ satisfies the \bvp
$$
\begin{cases}
\displaystyle{
(-\btu+I)\theta =0 }& \\
\displaystyle{
\pd{\theta}{n} = Y_k^* \bigl( (\nabla_{Y_k}\phi \lfloor) \vn \bigr)_K
+ (\phi\lfloor\vn)_K } & 
\end{cases} 
$$
Therefore,
\begin{align*}
\lefteqn{
(-\btu+I)\bigl( \eta\bigl[ (\K d\psi)_K -(d\K\psi)_K \bigr] \bigr)
}\\ 
& = L_1 (\K d\psi)_K +\eta\bigl[ \btu (d\K \psi)_K -(d\K \psi)_K
\bigr] + L'_1 (d\K\psi) \\
& = L_1 (\K d\psi)_K + L'_1 (d\K\psi) +L_2 (\K\psi)
	+\psi d\bigl[ (\btu \K\psi)_K -(\K\psi)_K \bigr] \\
& = L_1 (\K d\psi)_K +L'_1 (d\K\psi)_K +L'_2 (\K\psi) \\
& = L_1 (\K d\psi)_K +L_2 (\K\psi) .
\end{align*}
Next we analyze the boundary equation.  We have
$$
\pd{}{n} \bigl[ \eta(\K d\psi)_K -\eta(d\K\psi)_K \bigr]
= \eta\sum_{k=0}^{N} Y^*_k \bigl[ (\nabla_{Y_k} d\psi)\lfloor\vn 
\bigr]_K -\eta\pd{}{n} (d\K\psi)_K \, .
$$
Now we want to commute the normal derivative and the operator $d$ in
the far rightmost term in the equation above.  We have
\begin{align*}
\pd{}{n} (d\K\psi)_K 
& = \pd{}{n} \bigl( \sum_{i=1,\dots,N \atop I} 
X_i (\K\psi)_I \e{K}{iI} \bigr) + \pd{}{n}L_0 (\K\psi)  \\ 
& =  \sum_{i=1,\dots,N \atop I} 
X_i \bigl( \pd{}{n} (\K\psi)_I \bigr) \e{K}{iI} 
+ L'_1 (\K\psi)+ \pd{}{n}L_0 (\K\psi)  \\ 
& =  \sum_{i=1,\dots,N \atop I} X_i 
\biggl( \sum_{k=0}^{N} Y^*_k \bigl[
(\nabla_{Y_k}\psi)\lfloor\vn\bigr]_I \biggr) \e{K}{iI}
 + L'_1 (\K\psi)+ \pd{}{n}L_0 (\K\psi)  \\ 
& =  \sum_{i=1,\dots,N \atop I} X_i 
\biggl( \sum_{k=0}^{N} Y^*_k 
Y_k \psi_{0I}  \biggr) \e{K}{iI} +T_2 \psi 
 + L'_1 (\K\psi)+ \pd{}{n}L_0 (\K\psi)  \\ 
& = \sum_{k=0}^{N} Y^*_k 
Y_k  \biggl( 
 \sum_{i=1,\dots,N \atop I}  
X_i \psi_{0I}\e{K}{iI}   \biggr) 
+T'_2 \psi 
 + L'_1 (\K\psi) \\
& \qquad \qquad \qquad + L_0 \pd{}{n} (\K\psi) + L'_0 (\K\psi) \\ 
& = \sum_{k=0}^{N} Y^*_k \bigl[ (\nabla_{Y_k} d\psi)\lfloor\vn 
\bigr]_K 
+ T_2 \psi +L_0 (\K\psi) ,
\end{align*}
where
we have used the fact that $(\p/\p n) \K \psi$ on $b\Omega$
equals a second order tangential operator in $\psi$. 
Hence,
$$
\pd{}{n} \bigl[ \eta(\K d\psi)_K -\eta(d\K\psi)_K \bigr]
= T_2 \psi +L_0 (\K\psi) .
$$
Finally we prove the estimate.  Using \cite{TRI} Theorem 4.2.4 and
Corollary \ref{K-ord-1}, 
we have
\begin{align*}
\| w_1\|_s 
& \le c \bigl( \| L_1 \K d\psi \|_{F^{s-2}} 
+ \| L_2 \K \psi \|_{F^{s-2}} + \| T_2 \psi\|_{W^{s-3/2}(\Omega)} 
+ \| T_1 \K\psi\|_{W^{s-3/2}(\Omega)} \bigr) \\
& \le c \bigl( \| \K d\psi\|_{F^{s-1}} + \|\K\psi\|_s 
+\|\psi\|_{s+1} +\|\K\psi\|_s \bigr)\\
& \le c \| \psi\|_{s+1} . \qed
\end{align*}
\renewcommand{\qed}{}\endpf
\begin{lemma} \label{eta-G-Omega-psi-0not-inI}   \sl
Let $\psi\in\bw^q$.  Then for $I\not\ni0$, 
the function $[\eta G_\Omega \psi]_I \til$ is a
solution of the following \bvp: 
$$
\begin{cases}
\displaystyle{ (-\btu+\ep{\cal L}_2 +I) w = \bigl[ \eta_1 L_1
( G_\Omega \psi) \bigr]_I \til }
				& \on\rnp \\
\displaystyle{ \pd{w}{x_0}  =(-\btu_T +I)\pd{
\widetilde{\psi_{0I}}}{x_0}  + \zeta}	& \on\RR^N 
 \end{cases} 
$$
Here 
$\zeta$ is a function satisfying the estimate
$$
\| \zeta\|_{W^{s-3/2}(b\Omega)} \le c \| \eta_1 \psi\|_{s+1} ,
$$
and $L_1$ and  $\eta_1$ are as in the previous lemma.
\end{lemma}
\bgpf
As in the  proof of Lemma \ref{eta-G-Omega-psi-0inI}
we need only check that
the boundary equation is satisfied, and that the desired estimate
holds. 

Notice that, by Lemma \ref{eta-K-Omega-phi}, 
on the set
$\{ x_0 =0 \}$, for $I\not\ni0$, we have
$$
\pd{}{x_0} \bigl[ \eta\K d\psi \bigr]_I \til
= -\sum_{K\not\ni0} \tilde\eta\tilde\gamma_{IK} 
	[\K d\psi ]_K \til
	+(E_2 \widetilde{d\psi})_I \, ,
$$
where $E_2$ is defined as in Lemma \ref{eta-K-Omega-phi}.

On the other hand, if $I\not\ni 0$, on $\RR^N$ we have that 
\begin{align*}
\pd{}{x_0} \bigl[ \eta d\K \psi \bigr]_I \til 
& = \tilde\eta\sum_{J\not\ni0 \atop i=1,\dots,N} \pd{}{x_0}X_i
	(\widetilde{\K \psi})_I 
+\pd{}{x_0} \bigl[ \sum_J \bigl( \eta\K \psi\bigr)_J \til
d\omega^J \bigr] \\
& = - \sum_{K\not\ni 0\atop i=0,\dots,N} \tilde\eta a_{IK} X_i \bigl(
\widetilde{\K \psi} \bigr)_K 
+[L_0 (\widetilde{\K \psi})]_I +(T_2 \widetilde{\eta\psi} )_I \, ,
\end{align*}
where $T_2$ is a second order tangential differential operator.  Here
we have used the fact that $(\p/\p x_0) \K \psi$ on the set $\{
x_0 =0 \}$ equals a second order tangential operator on $\psi$. 
Therefore, for $I\not\ni 0$,
\begin{align*}
\left. \pd{w}{x_0}\right|_{\RR^N}
&  = (-\btu_T +I) \pd{\tilde\psi_{0I}}{x_0} \\
& \qquad - \sum_{K\not\ni0} 
\tilde\eta \gamma_{IK}  \bigl[\K d\psi\bigr]_K \til
-\sum_{i=1,\dots,N \atop K\not\ni0} \eta a_{IK} X_i \bigl(
\widetilde{\K \psi}\bigr)_K  
+[L_0 (\widetilde{\K \psi})]_I 
+(T_2 \widetilde{\eta\psi})_I  \\
& \equiv (-\btu_T +I) \pd{\tilde\psi_{0I}}{x_0} +\zeta .
\end{align*}
Hence we need only check the estimate.  

Let $w_1$ be given by Lemma \ref{LAST?}.  Then
$$
\zeta= -\sum_{K\not\ni0} \gamma_{IK} \tilde\eta\bigl[ 
d\K\psi\bigr]_K \til +w_1 
-\sum_{i=1,\dots,N \atop K\not\ni0} \eta a_{IK} X_i \bigl(
\widetilde{\K \psi}\bigr)_K  
+[L_0 (\widetilde{\K \psi})]_I 
+(T_2 \widetilde{\eta\psi})_I  .
$$
Using the above equality, Lemma \ref{marco2}, Lemma \ref{LAST?}, 
and Corollary \ref{K-ord-1}
we see that 
\begin{align*}
\| \zeta\|_{W^{s-3/2}(b\Omega)} 
&  \le c\biggl( \sum_{K\not\ni0}
\|\tilde \eta_1 \bigl( d\K \psi \bigr)_K \til 
\|_{W^{s-3/2}(b\Omega)}  +
\| w_1 \|_{W^{s-3/2}(b\Omega)} \\
& \qquad \qquad 
+ \|\eta_1 \K \psi\|_{W^{s-1/2}(b\Omega)}      
+ \| \widetilde{\eta\psi} \|_{W^{s+1/2}(b\Omega)} \biggr) \\
& \le c\bigl( 
\| \tilde\eta_1  \widetilde{\K \psi}  \|_{W^{s-1/2}(b\Omega)} 
+\| w_1 \|_{W^{s-1/2}(b\Omega)}  
+ \|\eta_1\psi\|_{s+1}  \bigr) \\
& \le c\| \eta_1 \psi\|_{s+1}  .
\end{align*}  
This proves the estimates and concludes the proof.
\endpf

\begin{pf*}{End of the Proof of \ref{marco1}}
We are finally able
to compare $\bigl[ \eta G_\Omega \psi \bigr]_I \til$ and
$\bigl[ G_{\rnp}(\widetilde{\eta\psi})\bigr]_I$.  
We begin with the case $I\ni0$. 
By (\ref{Gnot0}) and Lemma \ref{eta-G-Omega-psi-0not-inI}
it is easy to see that the function $w$ solves the \bvp
$$
\begin{cases}
\displaystyle{
(-\btu+I)w = (-\ep{\cal L}_2) \bigl[ \eta G_\Omega \psi
\bigr]_I \til 
+ \eta_1 L_1 \bigl[ G_\Omega \psi \bigr] \til } \qquad & \on \rnp \\ 
\displaystyle{
w = (-\ep{\cal T}_2) (\widetilde{\eta\psi})_I 
+\eta_1 (T_1 \psi)_I 
+\zeta}
& \on \RR^N 
\end{cases} 
$$
where $\zeta$ is as in Lemma \ref{eta-G-Omega-psi-0not-inI}.
Therefore, using 
Theorem 4.2.4 in \cite{TRI}, Lemma 
\ref{eta-G-Omega-psi-0not-inI}, and Corollary \ref{K-ord-1} 
we see that 
\begin{align*}
\| w\|_s 
& \le c\bigl( \ep \| {\cal L}_2(\eta G_\Omega \psi)_I \|_{F^{s-2}} 
+  \| \eta_1 L_1 (  G_\Omega \psi) \|_{F^{s-2}} 
+\|\ep {\cal T}_2 (\widetilde{\eta\psi}) +\eta T_1 \psi
\|_{W^{s-1/2} (b\Omega)} \\
& \qquad \qquad \qquad + \| \zeta\|_{W^{s-1/2}(b\Omega)} \bigr) \\
& \le c\bigl( \ep \| \eta G_\Omega \psi \|_{s} 
+\|  \eta G_\Omega \psi \|_{F^{s-1}}  +\| G_{\rnp}
(\eta\psi)\|_{F^{s-1}} +\ep \| \eta\psi\|_{s+2} 
+ \|\eta_1 \psi \|_{s+1} \bigr)  \\  
& \le c \bigl( \ep \|\eta\psi\|_{s+2} +\|\eta\psi\|_{s+1} \bigr) \, ,
\end{align*}  
where we also use the estimate 
$$
\|\eta G_\Omega \psi -G_\Omega (\eta\psi)\|_s 
\le c \|\eta_1 \psi\|_{s+1} .
$$
This follows by writing $G_\Omega =\K d+\K d$, and noticing that the
commutator between $\K$ and the multiplication by $\eta$ is a
$0$ order operator, since $\eta\K \theta -\K (\eta\theta)$ solves the \bvp
$$
\begin{cases}
\displaystyle{
(\btu -I)w = L_1 (\K \theta) }\qquad & \on \Omega \\
\displaystyle{ 
\pd{w}{n} = T_1 \theta} & \on b\Omega 
\end{cases} 
$$
This is a problem to which we apply Theorem 4.2.4 in \cite{TRI} again.

Now let $I\not\ni0$.  Then the function $w$ solves the \bvp
$$
\begin{cases}
\displaystyle{
(-\btu+I)w = (-\ep{\cal L}_2 ) \bigl[ \eta G_\Omega \psi
\bigr]_I \til 
+ \eta_1 L_1 \bigl[ G_\Omega \psi \bigr] \til }& \on \rnp \\
\displaystyle{
\pd{w}{x_0} = (\ep{\cal T}_2 ) 
\pd{(\widetilde{\eta\psi})_{0I}}{x_0}  +\zeta }
& \on \RR^N 
\end{cases} 
$$
where $\zeta$ is given by Lemma \ref{eta-G-Omega-psi-0not-inI}. 
Therefore, by Corollary \ref{K-ord-1},
\begin{align*}
\| w\|_s
& \le    \ep \| \eta\psi \|_{s+2} + c\left( 
\| \eta_1 G_\Omega \psi \|_{F^{s-1}(\Omega)}
+\| \zeta \|_{W^{s-3/2}(b\Omega)}  \right) \\
& \le \ep \| \eta\psi \|_{s+2} + c \| \eta_1  \psi \|_{s+1} .\qed 
\end{align*}
\renewcommand{\qed}{}\endpf

\section{The Decomposition Theorem and Conclusions}
\setcounter{equation}{0}


\subsection*{Proof of the Main Result} We are now in a position to
finish the proof of Theorem \ref{MAIN-THM-DMN}.

Let $X$ be the subspace of $W^{s+2}_q (\Omega)$, ($s>1/2$) consisting
of the
forms $\phi$ satisfying the boundary conditions $\phi,\, d\phi\,
\in\dom d^*$.  Observe that $X$ is a closed subspace of $W^{s+2}_q
(\Omega)$, (this  is the reason why $s>1/2$).  By the Regularity
Theorem \ref{REG-THM-Q} we know that, on $X$, the norm
$$
\| \psi\|_{s+1} + \|(-\btu+G_\Omega )\psi\|_s 
$$ 
is equivalent to the norm $\| \psi\|_{s+2}$.  Let ${\cal P}$ be the
orthogonal projection of $X$ onto the kernel of $-\btu +G_\Omega$.
Then, by standard functional analysis arguments, (see \cite{ZIE} p.178),
$$
\| \psi-{\cal P}\psi\|_{s+1} \le C \| -\btu \psi+G_\Omega \psi\|_s .
$$
Therefore, if $\psi$ is orthogonal to the kernel of $-\btu+G_\Omega$,
the regularity theorem 
Theorem \ref{REG-THM-Q} says that
$$
\|\psi\|_{s+2} \le C_s \|-\btu\psi+G_\Omega \psi\|_s .
$$
We now complete the proof of the existence theorem---Theorem 
\ref{EXISTENCE}.
For $\alpha\in W^1_q \cap {\cal H}_q^\perp$ fixed, we approximate
$\alpha$ by $\alpha_m \in Z\equiv (F-I)(\D)$.  For such $\alpha_m$'s
we can find $\phi_m \in\D$ orthogonal to ${\cal H}_q$, solutions of
the \bvp, and such that, for $\alpha\ge1$,
$$
\| \phi\|_{s+2} \le c_s \|\alpha\|_s \, . 
$$
Thus the $\phi_m$'s converge in $W^{s+2}_q$ to a certain $\phi$.
Clearly $\phi$ solves the \bvp\ with data $\alpha$, since the
boundary conditions are preserved in the $W^{s+2}$ topology for
$s>1/2$. 

Hence, if 
$\alpha\in W^s (\Omega)$, ($s\ge1$) and $\alpha$ 
is orthogonal to
$\ker(-\btu+G_\Omega)$ in the $W^1$ inner product, 
then there exists
a unique solution to the problem
(\ref{bvp-domain})---$\phi\in W^{s+2} (\Omega)$ 
orthogonal to $\ker(-\btu+G_\Omega)$ and such that
$$
\|\phi\|_{s+2} \le C_s \|\alpha\|_s .
$$
If $\alpha\in W^s(\Omega)$ with $1/2 <s<1$, then we still have existence
and regularity by using a density argument, since $\alpha\in W^s_q
(\Omega)$ can be approximated by $\alpha_m \in W^1_q (\Omega)$ in the
$W^s$-topology.  This concludes the proof of Theorem
\ref{MAIN-THM-DMN}.  \qed

\subsection*{Decomposition of $W^1_q$}
The aim of this part is to prove the following result:
\begin{theorem}\label{decomp-thm}   \sl
Let $\Omega$ be a smoothly bounded domain in $\RR^{N+1}$.  Let $W^1_q
(\Omega)$ denote the 1-Sobolev space of $q$-forms.  Then we have the
strong orthogonal decomposition
$$
W^1_q = dd^* (W^1_q) \bigoplus d^* d(W^1_q) \bigoplus {\cal H}_q\, ,
$$
where ${\cal H}_q$ is a finite dimensional subspace, and 
 $d^*$ denotes the $W^1_q$-Hilbert space adjoint of $d$. 
\end{theorem}
\bgpf
All of this is standard.  If $\alpha$ is orthogonal to 
${\cal H}_q$ then $\alpha$ belongs to
$(dd^* +d^* d)W^1_q$.  The fact that $d(W^1_q)$ and $d^* (W^1_q)$ are
orthogonal subspaces is also clear since
$$
\l d\phi,d^*\psi\r_1 = \l d^2 \phi ,\psi\r_1 = 0,
$$
for all $\phi,\psi\in C^\infty (\overline{\Omega})$.
\endpf
In the case $q=0$ we are able to determine
the harmonic space ${\cal H}_0$.  The rest of the section is devoted
to this end.
\begin{theorem}\label{H_0}   \sl
We have that ${\cal H}_0 = \{ \text{\rm constants}\}$.   
More precisely,
the \bvp
\begin{equation}
\begin{cases}\label{our-bvp}
\displaystyle{
d^* d u  =   f  } & \qquad \on \Omega  \\
du  \in  \dom d^*  & 
\end{cases}
\end{equation}
has a unique solution $u$ orthogonal (in the $W^1 (\Omega)$ inner
product) to the constant functions, for each $f\in W^1 (\Omega)$
which is also orthogonal to the constants, i.e.
$$
\int_\Omega f =0 .
$$
\end{theorem}
This theorem is a consequence of the next theorem.  
By equation (\ref{Jan3-1}) we can
rewrite the \bvp\  (\ref{our-bvp})  as 
\begin{equation}\label{BP'} 
\begin{cases}
\displaystyle{ 
- \btu u + G_\Omega u  = f } \qquad & \on \Omega \\
\displaystyle{
\sum_{j=0}^{N} n_j
\pd{}{n} \left (\pd{u}{x_j} \right ) =0} &  \on b\Omega 
\end{cases}  \ ,
\end{equation}
where $\vn =(n_0,\dots,n_N)$ is the normal vector field.
By our choice of $\vn$, and by
setting $G_\Omega u = v$ we obtain the new system
$$
\begin{cases}
\displaystyle{
- \btu u + v  =  f   } \qquad & \on \Omega  \\
\displaystyle{\pd{^2 u}{n^2} = 0 } &  \on b\Omega  \\
\displaystyle{ 
\btu v - v  =  0  } & \on \Omega  \\
\displaystyle{
\pd{v}{n} =  \pd{u}{n} + 
\sum_{k=0}^{N} Y_k^* \bigl[ ( \nabla_{Y_k} du)\lfloor\vn \bigr] } 
& \on b\Omega 
\end{cases} 
$$
Now assume $f$ to be suffiently regular.
By solving the first equation for $v$ and substituting
we find that (\ref{BP'}) finally becomes 
\begin{equation}\label{bvp-we-solve}
\begin{cases}
\displaystyle{
\btu^2 u - \btu u  =  f - \btu f \qquad }& \on \Omega \\
\displaystyle{
 \pd{^2 u}{n^2} =  0 } & \on b\Omega  \\
\displaystyle{
\pd{}{n} \btu u 
-\sum_{k=0}^{N} Y_k^* \bigl[ ( \nabla_{Y_k} du)\lfloor\vn \bigr] 
- \pd{u}{n} =  - \pd{f}{n} } & \on b\Omega  
\end{cases} 
\end{equation}
This last is the problem that we are
going to study.
\begin{theorem}   \sl
The boundary value problem (\ref{bvp-we-solve})
is elliptic.  It has index 0.  The kernel of the system is given by
$\{ u = \text{constant}\}$ and a (unique) solution exists
for each piece of data $f$ that satisfies
$$
\int_\Omega f = 0 \ , \qquad \qquad f
\in W^r , \ r \ge 2.
$$
\end{theorem}
\begin{remark}  \rm
The condition $f \in W^r$, $r \ge2$ is
necessary in order to
guarantee that $\btu f \in L^2(\Omega).$
In this case the solution $u$ satisfies the estimates
$$
\|u\|_{W^{s+4}(\Omega)}
 \leq C \biggl ( \|f - \btu f \|_{W^s(\Omega)}
      + \left \| \frac{\partial f}{\partial n}
  \right \|_{W^{1/2+ s}(b\Omega)} + \|u\|_s \biggr ),
$$
that is, for $ s\ge 0$,
$$
 \| u \|_{W^{s+4}(\Omega)}
\leq C \bigl ( \|f\|_{s+2} + \|u\|_s \bigr ) .
$$
This result, when applied to problem (\ref{our-bvp}),
 is not optimal.  This is the reason why we had to go
through the more complicated computations in the previous sections.
\end{remark}
\bgpf
In order to prove that the \bvp\ is elliptic it suffices to verify
the Lopatinski condition as in \cite{HOR1}.  This is standard, and
since we do not use this fact in what follows, we leave the details
to the reader. 
Thus the \bvp\ (\ref{bvp-we-solve})
admits
solutions with appropriate estimates modulo
a finite dimensional kernel
and a finite dimensional cokernel.  Our next
job is to explicitly
determine this kernel and cokernel.
\medskip \\
{\bf Determination of the Kernel:}
Suppose that the function
$u$, defined on $\Omega$, satisfies
$$
\begin{cases}
\displaystyle{
\btu (\btu u) - \btu u  =  0 }\qquad & \on \Omega \\
\displaystyle{
\pd{^2 u}{n^2} 
  =  0 } &  \on b\Omega \\
\displaystyle{
\pd{}{n} \btu u 
-\sum_{k=0}^{N} Y_k^* \bigl[ ( \nabla_{Y_k} du)\lfloor\vn \bigr] 
- \pd{u}{n}  =  0 } & \on b\Omega
\end{cases}
$$
We intend to show that $u$ must therefore be
constant.  This will
then imply that the kernel is the one
dimensional space of constant
functions.  Without loss of generality we suppose that
all functions are {\em real
valued}.  
 
Using Green's theorem we see that
\begin{multline*}
0 = \int_\Omega u  \btu \bigl[ \btu u - u\bigr]
 \, dV  \\
=  \int_\Omega |\btu u|^2 + \int_{b\Omega}
  u\pd{}{n} (\btu u) 
 - \int_{b\Omega} \btu u \pd{u}{n} + \int_\Omega |\grad  u|^2 -
\int_{b\Omega} u \pd{u}{n} .
\end{multline*}
Thus we find that
\begin{eqnarray}
\lefteqn{
\int_\Omega |\btu u|^2 + |\grad  u|^2} \nonumber \\
& = & \int_{b\Omega} \pd{u}{n} \bigl( \btu u +u \bigr)
 - u\pd{}{n} (\btu u)  \nonumber \\ 
 & \stackrel{(2^{\text{nd}}\ \text{bdry. cond.})}{=} &
 \int_{b\Omega} \pd{u}{n} \bigl( \btu u +u \bigr) 
- \int_{b\Omega} u\left[ 
\sum_{k=0}^{N} Y_k^* \bigl[ ( \nabla_{Y_k} du)\lfloor\vn \bigr] 
+ \pd{u}{n} \right]  \nonumber \\ 
  & = & \int_{b\Omega} \pd{u}{n}\btu u -
  \int_{b\Omega} u    
\sum_{k=0}^{N} Y_k^* \bigl[ ( \nabla_{Y_k} du)\lfloor\vn \bigr]  
.  \label{ABOVE}
\end{eqnarray}
Now we consider the second integral on the right hand side 
of equation (\ref{ABOVE}).
We have   
\begin{eqnarray*}
\int_{b\Omega} 
\sum_{k=0}^{N} Y_k^* \bigl[ ( \nabla_{Y_k} du)\lfloor\vn \bigr]  u 
&  =&
\sum_k \int_{b\Omega} Y_k u
\bigl[ (\nabla_{Y_k}du\lfloor)\vn\bigr]   \\
&  \stackrel{(1^{\text{st}}\ \text{bdry. cond.})}{=} &
\sum_k  \int_{b\Omega} Y_k u \bigl[ \pd{(Y_k u)}{n} \bigr] \\
& = & - \frac{1}{2} \int_{b\Omega}  \pd{}{n} |\grad  u|^2 ,
\end{eqnarray*}
where we use the first boundary condition, and 
the simple
facts that $(\nabla_{Y_k}du)$ $\lfloor\vn = (\p/\p n)Y_k u$, 
and  that
$|\grad u|^2 =\sum_k |Y_k u|^2 +|(\p/\p n)u|^2$.  

Substituting this last equality
into (\ref{ABOVE}), and using Green's 
theorem again, 
we find that
\begin{eqnarray*}
\int_\Omega |\btu u|^2 + |\grad  u|^2
& = & \int_{b\Omega} \btu u 
 \pd{u}{n} - \frac{1}{2}
 \pd{}{n} |\grad  u|^2  \\
& =  & \int_\Omega |\btu u|^2 
	+ \grad  u \cdot \grad (\btu u) 
	- \frac{1}{2} \btu (|\grad  u|^2) .
\end{eqnarray*}
Notice that
\begin{eqnarray*}
\btu |\grad  u|^2 & = & 2 \sum_j \left( \btu
\pd{u}{x_j} \right) \pd{u}{x_j}
+ 2 \grad  \pd{u}{x_j} \cdot \grad  \pd{u}{x_j} \\
 & = & 2 \grad  (\btu u) \cdot \grad  u + 2 \sum_j
 \grad  \pd{u}{x_j} \cdot \grad  \pd{u}{x_j} .
\end{eqnarray*}
As a result,
$$
\int_\Omega |\btu u|^2 +
|\grad  u|^2 = \int_\Omega |\btu u|^2
  - \sum_j \grad  \pd{u}{x_j} \cdot
     \grad  \pd{u}{x_j} .
$$
That is,
$$
\int_\Omega |\grad  u|^2 + \sum_j
\left | \grad  \pd{u}{x_j}
\right |^2 = 0 .
$$
This equality implies that $u$ is a constant.
\medskip \\
{\bf Determination of the Cokernel:}
We want to show that if $f$ is
such that the system (\ref{bvp-we-solve}) admits
a solution then
$$
\int_\Omega f = 0 .
$$
Notice that
this would imply that $\text{Cokernel}
\supseteq \bigl (
\text{ker} (f \mapsto \int_\Omega f) \bigr )^\perp.$
 
We have
\begin{eqnarray*}
\int_\Omega \btu(\btu u) - \btu u & = & \int_{b\Omega}
	\pd{}{n} (\btu u -u)\\
& \stackrel{(2^{\text{nd}}\ \text{bdry. cond.})}{=} & 
\int_{b\Omega}  \sum_k Y_k^* \bigl[ (\nabla_{Y_k}du\lfloor)\vn\bigr] 
	- \pd{f}{n} \\
&\stackrel{(\text{Green on }b\Omega)}{=} &
 	-\int_{b\Omega} \pd{f}{n}  \\
& = & - \int_{\Omega} \btu f .
\end{eqnarray*}
Now the left hand side in these last
equations equals $\int_{\Omega} (f - \btu f)$.
Therefore
$$
\int_\Omega f = 0 .
$$
 Now we examine the dimension of the cokernel.  Suppose that
$$
(g; v_1, v_2) \in C^\infty(\overline{\Omega})
\times C^\infty(b\Omega)
      \times C^\infty(b\Omega)
$$
and 
$$
(g; v_1, v_2) \perp \text{Range of} \, (\ref{bvp-we-solve}) ,
$$
i.e.\
\begin{multline}
\int_\Omega \bigl ( \btu^2 u - \btu u \bigr ) g +
 \int_{b\Omega} \left [ \sum_{j=0}^{N}
n_j \pd{}{n}
 \left ( \pd{u}{x_j}
\right ) \right ] v_1  \\
+ \int_{b\Omega}
\left [ \pd{}{n}
\btu u + \sum_{j=1}^N \div^T \left ( n_j
\nabla^T \pd{u}{n} \right )
- \frac{\partial u}{\partial n} \right]
v_2 = 0 \label{cokernel} 
\end{multline}
for all $u \in C^\infty.$
Then we want to show that
$$
(g; v_1, v_2) \ \text{lies in some one dimensional subspace.}
$$
This will finish the proof.

Thus suppose that (\ref{cokernel}) 
holds for all
$u \in C^\infty(\overline{\Omega}).$
By taking $u \in C^\infty_0(\Omega)$ and
integrating by parts we obtain
\begin{equation*}
0  =  \int_\Omega g(\btu^2 u - \btu u) 
 =  \int_\Omega u \bigl(
\btu^2 g - \btu g \bigr) .
\end{equation*}
This equality implies that $\btu^2 g -
\btu g = 0$ on $\Omega.$
 
Next, assume (\ref{cokernel}) and apply Green's
theorem to the integral on
$\Omega.$  It equals
\begin{eqnarray*}
\lefteqn{\int_\Omega \btu g(\btu - I) u  
+ \int_{b\Omega} g\pd{}{n}
(\btu - I) u - \int_{b\Omega} \pd{g}{n}(\btu - I)u} \\ 
 & = & \int_\Omega u (\btu^2 g - \btu g) +
\int_{b\Omega} \pd{u}{n}  \btu g
 - \int_{b\Omega} u \pd{}{n} \btu g  \\
 &   & \qquad + \int_{b\Omega} g\pd{}{n}
 \btu u  - \int_{b\Omega} g\pd{u}{n} - \int_{b\Omega} \btu u
\pd{g}{n}
 + \int_{b\Omega} u \pd{g}{n} \, .
\end{eqnarray*}
Substituting this last identity into (\ref{cokernel}),
and recalling that
$\btu^2 g - \btu g = 0$ on $\Omega$ we find that
\begin{eqnarray}
0 
& = & \int_{b\Omega}  \pd{u}{n} \btu g
	- \int_{b\Omega} u \pd{}{n} \btu g
	+ \int_{b\Omega} g\pd{}{n} \btu u 
	- \int_{b\Omega} g\pd{u}{n}  \nonumber\\ 
&   & \quad - \int_{b\Omega} \btu u \pd{g}{n}
	+ \int_{b\Omega} u \pd{g}{n}
	+ \int_{b\Omega} \pd{^2 u}{n^2}  v_1    \nonumber\\
&   & \quad + \int_{b\Omega} v_2 \pd{}{n}\btu u 
-\sum_k \int_{b\Omega} Y_k^* \bigl[ (\nabla_{Y_k}du)\lfloor\vn\bigr]  
 	v_2 - \int_{b\Omega} v_2 \pd{u}{n}  \nonumber\\
& = & \int_{b\Omega} u \left[ \pd{g}{n} - \pd{}{n} \btu g \right] 
	+ \int_{b\Omega} \pd{u}{n}
	\bigl[ \btu g - g - v_2 \bigr] \nonumber\\
&   & + \int_{b\Omega} \pd{^2 u}{n^2} v_1
	-\sum_k \int_{b\Omega} 
	Y_k^* \bigl[ (\nabla_{Y_k}du)\lfloor\vn\bigr]  v_2 
	- \int_{b\Omega} \btu u \pd{g}{n} \nonumber\\
&   &  + \int_{b\Omega} (g+v_2 )\pd{}{n} \btu u
 . \label{bdry-eq} 
\end{eqnarray}
Therefore we have
\begin{multline}
0=\int_{b\Omega} u \left[ \pd{g}{n} - \pd{}{n} \btu g \right] 
+ \int_{b\Omega} \pd{u}{n} \bigl[ \btu g - g - v_2 \bigr] 
 + \int_{b\Omega} \pd{^2 u}{n^2} v_1 \\
-\sum_k \int_{b\Omega} \bigl[ \pd{}{n}Y_k u]  Y_k v_2 
	-\int_{b\Omega} \btu u \pd{g}{n} 
	+\int_{b\Omega}\pd{}{n} \btu u(g+v_2 ) . \label{bdry-eq-2}
\end{multline}
As in the half space case,
we write $x'$ to indicate the
variable in $b\Omega$.  Given any
function $v$ on $b\Omega$, we let
$$
\tilde{v}(x) = v(\pi(x)) \equiv v(x') ,
$$
where $\pi$ is the normal projection
defined in a neighborhood
of $b\Omega.$  The functions $\tilde v$ are defined in a 
(fixed)
neighborhood of $b\Omega.$  Then we can extend them to all of
$\rnp$ by multiplying by a (fixed) cut-off function.  
 
Now we proceed in several
steps to calculate the cokernel:
\medskip \\
{\bf Step 1:}  Take $u(x)$ to be of the form
$$
u(x) = \rho(x)^3   \tilde{u}_1 (x) \, ,
$$
where $\tilde{u}_1 \in C^\infty(b\Omega)$
is generic.  Such a choice
yields that
$$
g + v_2 = 0 \qquad \qquad \text{on} \ b\Omega .
$$
{\em Proof of Step 1:}  For such a choice of $u$,
equation (\ref{bdry-eq-2}) becomes
$$
0 = \int_{b\Omega} (g+v_2 ) \pd{}{n} \btu \bigl(
\rho^3 \tilde{u}_1 \bigr)
$$
Now,
\begin{eqnarray*}
\pd{}{n} \btu \bigl[ \rho^3 \tilde{u}_1
\bigr] \biggr |_{b\Omega}
  & = & \pd{}{n} \bigl[ \tilde{u}_1 \btu \rho^3   
  + 2 \bigl( \grad  \rho^3 \cdot \grad  \tilde{u}_1 \bigr)
  + \rho^3 \btu \tilde{u}_1 \bigr] \biggr|_{b\Omega} \\
 & = & \left ( \pd {} n \btu \rho^3 \right )   \tilde{u}_1
               \biggr |_{b\Omega} \\
                & \equiv & 6 u_1 .
\end{eqnarray*}
[It is a standard fact that we may assume that $(\partial \rho/\partial n)
= |\grad\, \rho| = 1$ on $b\Omega.$]
This implies that $g + v_2 = 0$ on $b\Omega.$
\medskip \\

\noindent {\bf Step 2:}  Take $u = \rho^2(x) \tilde{u}_1(x)$
with $u_1 \in C^\infty(b\Omega).$
This implies that
$$
v_1 = \pd g n \qquad \qquad \text{\rm on} \ b\Omega.
$$
{\em Proof of Step 2:}  Using Step 1, we see that this
particular choice of $u$ 
gives
\begin{eqnarray*}
   \btu u \biggr |_{b\Omega} & = & \btu \bigl
                         ( \rho^2 \ut_1 \bigr )
                        \biggr |_{b\Omega} \\
        & = & \bigl( \ut_1\btu \rho^2 
               + 2 \grad  \rho^2 \cdot \grad  \ut_1
               + \rho^2 \btu \ut_1 \bigr )
\biggr |_{b\Omega} \\
        & = & 2 |\grad  \rho |^2 \ut_1 \biggr |_{b\Omega} \\
        & = & 2 u_1 ,
\end{eqnarray*}
while $(\p^2 \rho^2 /\p n)=2\ \on b\Omega$.  
Therefore, 
$$
v_1 = \pd{g}{n} \qquad \qquad \on b\Omega .
$$

Hence, putting Steps 1 and 2 together, equation 
(\ref{bdry-eq-2})
becomes
\begin{multline}
0= 
\int_{b\Omega} u \left[ \pd{g}{n} - \pd{}{n} \btu g \right]
     + \int_{b\Omega} \pd{u}{n} \btu g 
     + \int_{b\Omega} \pd{^2 u}{n^2} \pd{g}{n} \\
+ \sum_k \int_{b\Omega} \bigl[ \pd{}{n}Y_k u]  Y_k \pd{g}{n} 
	-\int_{b\Omega} \btu u \pd{g}{n}  \label{bdry-eq-3}
\end{multline}
for all $u \in C^\infty(\overline{\Omega})$.
\medskip \\
{\bf Step 3:}  We have that 
\begin{multline}
\int_{b\Omega} \pd{u}{n} \btu g 
+\int_{b\Omega} \pd{^2 u}{n^2} \pd{g}{n} 
\sum_k +\int_{b\Omega} \bigl[ \pd{}{n}Y_k u]  Y_k \pd{g}{n} 
	-\int_{b\Omega} \btu u \pd{g}{n} \\
= \int_{b\Omega} u\biggl[ Y_k^* \pd{}{n} D_k g \biggr] 
+ \int_{b\Omega} \pd{u}{n} \pd{^2 g}{n}  \label{bdry-eq-4}
\end{multline}
{\em Proof of Step 3:}  
Recall that for $k=0,\dots,N$, 
$Y_k \equiv D_k -n_k (\p/\p n)$, and that $\vn = \sum_j n_j D_j$.
Therefore, it is easy to see that
$$
\pd{^2 u}{n^2}\pd{g}{n} +\sum_k (Y_k g)\bigl( \pd{}{n} Y_k u\bigr) 
=\sum_k D_k g \pd{}{n} (D_k u) .
$$
Thus a repeated
application of Green's
theorem gives that 
\begin{align*}
\lefteqn{ 
\int_{b\Omega} \pd{u}{n} \btu g 
+\int_{b\Omega} \pd{^2 u}{n^2} \pd{g}{n} 
+\sum_k \int_{b\Omega} \bigl[ \pd{}{n}Y_k u]  Y_k \pd{g}{n} 
	-\int_{b\Omega} \btu u \pd{g}{n} }\\
& = \int_{b\Omega} \pd{u}{n} \btu g 
+ \sum_k \int_{b\Omega} \bigl[ \pd{}{n}Y_k u]  Y_k \pd{g}{n} 
-\int_{b\Omega} \btu u \pd{g}{n} \\
& = \int_{b\Omega} \pd{u}{n} \btu g 
+ \sum_k \biggl[ \int_\Omega \grad D_k u \cdot \grad D_k u 
\biggr] - \int_\Omega \btu u \btu g \\
& = \int_{b\Omega} \pd{u}{n} \btu g 
+ \sum_k \int_{b\Omega} \bigl( \pd{u}{n} \pd{D_k g}{n} 
-\int_{b\Omega} \pd{u}{n}\btu g \\
& = \sum_k \int_{b\Omega} u \bigl( Y_k^* \pd{}{n} D_k g \bigr)
+ \sum_k \int_{b\Omega} \pd{u}{n} \bigl( n_k \pd{}{n}D_k g \bigr) \\ 
& = \sum_k \int_{b\Omega} u\bigl( Y_k^* \pd{}{n}D_k g\bigr) 
+ \int_{b\Omega} \pd{u}{n} \pd{^2 g}{n^2} .
\end{align*}
This proves Claim 3.

Now because of Step 3, equality (\ref{bdry-eq-3}) becomes
\begin{equation} \label{bdry-eq-5}
\int_{b\Omega} u \left[ \pd{g}{n} - \pd{}{n} \btu g 
+\sum_k \bigl( Y_k^* \pd{}{n}D_k g \bigr) \right]
     + \int_{b\Omega} \pd{u}{n} \pd{^2 g}{n} .
\end{equation}
Thus $g$ must satisfy the system
$$
\begin{cases}
\displaystyle{
\btu^2 g - \btu g  =  0} & \on \Omega \\
\displaystyle{
\sum_j   \pd{^2 g}{n^2}  = 0 } & \on b\Omega \\
\displaystyle{
\pd{g}{n}- \pd{}{n} \btu g  + \sum_k Y_k^* \pd{}{n}D_k g
= 0} \qquad & \on  b\Omega
\end{cases} 
$$
Using the first boundary equation we can rewrite the second one 
as in the homogeneous problem that we have studied earlier
when we calculated the kernel.
We conclude that $g \equiv C$.  Therefore
the cokernel of the problem
is given by
$$
(g; v_1, v_2) = (c; 0, - c) ,
$$
which is obviously a one dimensional space.

All of these calculations can be applied to our original problem provided that 
$f\in W^2 (\Omega)$.  Using the regularity result we can extend it to
all $f\in W^s$, $s>1/2$ and conclude that has a unique solution 
$u$ with $\int_\Omega u =0$, 
for
all $f\in W^s (\Omega)$ with $\int_\Omega f=0$.
\endpf 
 
\section{Final Remarks}
 \setcounter{equation}{0}

In the present paper we have laid the foundations for
Hodge theory of the standard exterior differential operator $d$
in the Sobolev topology acting on $q$-forms on a domain in $\RR^N$.
The theory that have presented here is complete.
However, there is much work that remains to be done.  For maximum
applicability, the Hodge theory on a (compact) manifold
needs to be developed.  Also the case of all $s$ and all $q$ should
be treated in a unified manner.
 
The ultimate goal of this program is to develop the Hodge theory
of the $\overline{\partial}$-Neumann problem in the Sobolev topology
on a strongly pseudoconvex domain.  We will turn to this
task in a subsequent sequence of papers. 
\mbox{}\vspace{1cm}\\

\end{document}